\documentclass[]{article}
\usepackage{graphicx}
\usepackage{amsfonts}
\usepackage{amsmath}
\usepackage{amssymb}
\usepackage{fancyhdr}
\usepackage{titlesec}
\usepackage{indentfirst}
\usepackage{booktabs}
\usepackage{verbatim}
\usepackage{color}
\usepackage{amsthm}
\usepackage{caption}
\usepackage{subcaption}
\usepackage{physics}
\usepackage{mathtools}
\usepackage{tikz}
\usepackage[ruled,vlined,linesnumbered]{algorithm2e}
\usepackage{algpseudocode}

\usepackage{hyperref}
\usepackage{cleveref}
\usepackage[page,header]{appendix}
\usepackage{titletoc}

\newcommand{\rd}{\,\mathrm{d}}

\numberwithin{equation}{section}

\newcommand{\rom}[1]{\text{\MakeUppercase{\romannumeral #1}}}
\newcommand{\rombracket}[1]{(\text{\MakeUppercase{\romannumeral #1}})}

\newcommand{\jh}[1]{\textcolor{blue}{JH: #1}}
\newcommand{\ki}[1]{\textcolor{magenta}{[gw: #1]}}

\newcommand{\bA}{\mathbf{A}}
\newcommand{\bB}{\mathbf{B}}
\newcommand{\bE}{\mathbf{E}}
\renewcommand{\bf}{\mathbf{f}}
\newcommand{\bF}{\mathbf{F}}
\newcommand{\bG}{\mathbf{G}}
\newcommand{\bM}{\mathbf{M}}
\newcommand{\bH}{\mathbf{H}}
\newcommand{\bJ}{\mathbf{J}}
\newcommand{\bu}{\mathbf{u}}
\newcommand{\bx}{\mathbf{x}}
\newcommand{\bv}{\mathbf{v}}

\newcommand{\dt}{\Delta t}
\renewcommand{\dv}{\Delta v}
\newcommand{\dx}{\Delta x}
\newcommand{\e}{\mathrm{e}}
\newcommand{\pdt}{\partial_t}
\newcommand{\pdx}{\partial_x}
\newcommand{\gradx}{\nabla_\bx}
\newcommand{\gradv}{\nabla_\bv}
\newcommand{\T}{\mathrm{T}}

\begin{document}

\title{Dynamical Tensor Train Approximation for Kinetic Equations}

\author{
Geshuo Wang\footnote{Department of Applied Mathematics, University of Washington, Seattle, WA 98195 (geshuo@uw.edu).} \
\ and \ 
Jingwei Hu\footnote{Department of Applied Mathematics, University of Washington, Seattle, WA 98195 (hujw@uw.edu). Corresponding author.}
}

\maketitle

\begin{abstract}
The numerical solution of kinetic equations is challenging due to the high dimensionality of the underlying phase space.
In this paper, we develop a dynamical low-rank method based on the projector-splitting integrator in tensor-train (TT) format.
The key idea is to discretize the three-dimensional velocity variable using tensor trains while treating the spatial variable as a parameter, thereby exploiting the low-rank structure of the distribution function in velocity space.
In contrast to the standard step-and-truncate approach, this method updates the tensor cores through a sweeping procedure, allowing the use of relatively small TT-ranks and leading to substantial reductions in memory usage and computational cost.
We demonstrate the effectiveness of the proposed approach on several representative kinetic equations.
\end{abstract}

{\small 

}

\section{Introduction}
Kinetic equations are fundamental partial differential equations (PDEs) that describe the statistical evolution of large particle systems in phase space. 
The unknown is the one-particle distribution function $f(t,\bx,\bv)$, which depends on time $t$, position $\bx\in\mathbb{R}^d$, and velocity $\bv\in\mathbb{R}^d$ (with $d$ denoting the dimension, typically ranging from 1 to 3). 
Different kinetic models play central roles in various fields.
For instance, the Boltzmann equation is fundamental in rarefied gas dynamics \cite{cercignani}, while the Vlasov-Fokker-Planck equation governs the behavior of charged particles in plasma physics \cite{Villani02}.
Due to their importance, the numerical simulation of kinetic equations has been an active research area in scientific computing. However, the high dimensionality of the phase space makes direct discretization expensive in terms of both memory and computational cost, posing the main challenge for numerical simulations.

To mitigate this issue, researchers have exploited the fact that solutions of many kinetic equations often lie close to low-dimensional manifolds in the high-dimensional solution space. Low-rank methods leverage this observation by approximating the distribution function using separable representations, thereby reducing both storage requirements and computational complexity. While kinetic equations do not generally admit exact low-rank solutions, in many practical situations their dynamics remain near low-rank structures, making such approximations highly effective.

For problems involving two variables, the singular value decomposition (SVD) provides a standard low-rank tool. In higher dimensions, tensor decompositions such as the Tucker format \cite{tucker1966some}, the canonical polyadic decomposition \cite{hitchcock1927expression}, and the tensor-train (TT) format \cite{oseledets2011tensor} are widely employed.
Tensor methods generalize matrix decompositions to capture multi-way interactions, alleviating exponential storage scaling.
In particular, the TT format, also known in quantum physics as the matrix product state (MPS) \cite{white1992density,white1993density,schollwock2005density,schollwock2011density}, has emerged as a powerful representation. It has been successfully applied in a variety of areas, notably in quantum dynamics \cite{haegeman2016unifying,jorgensen2019exploiting,ye2021constructing,sun2024efficient}.

Within the low-rank framework for solving time-dependent PDEs, two prominent strategies have been developed. The step-and-truncate (SAT) approach \cite{dolgov2014low,kormann2015semi,ehrlacher2017dynamical} integrates the full high-dimensional equation over a short time step, after which the solution is truncated back to the chosen low-rank manifold, commonly via truncated SVD or tensor truncation. In contrast, the dynamical low-rank (DLR) approximation \cite{koch2007dynamical,einkemmer2018low} evolves the solution directly on the low-rank manifold by projecting the governing equation onto its tangent space. While direct implementations of DLR can suffer from numerical instabilities related to matrix inversions, these issues are mitigated by the projector-splitting integrator, which has been extended from matrices \cite{lubich2014projector} to tensor trains  \cite{lubich2015time}. DLR methods typically use a fixed rank throughout the simulation, providing predictable computational cost. Hybrid approaches that combine these principles have also been proposed \cite{CKL22,kieri2019projection,nakao2025reduced}, further broadening the methodological landscape. Many other variants of low-rank integrators exist; without attempting to be exhaustive, we refer the reader to the recent review \cite{EKKMQ25} for a discussion in the context of solving kinetic equations.

To balance tractability and realism, many studies focus on reduced kinetic models, such as one spatial dimension with one to three velocity components (1D1V, 1D2V, 1D3V) or two spatial and two velocity dimensions (2D2V). 
These simplified settings retain essential kinetic features while allowing the development and testing of new algorithms. 
Another key modeling decision lies in the choice of decomposition. A common approach is to separate the spatial variable $\bx$ and the velocity variable $\bv$ \cite{einkemmer2018low,HW22,coughlin2024robust,einkemmer2025asymptotic}.
More refined decompositions split all components of $(\bx,\bv)$ into higher-order tensor structures \cite{dolgov2012fast,kormann2015semi,guo2024conservative,YL24}, enabling a richer but more computationally demanding representation.

In this work, we propose a TT-based method tailored for kinetic equations. Unlike previous works, our method treats the velocity and spatial variables differently: the velocity variable is represented as a three-dimensional tensor in TT format, while the spatial variable is treated as a parameter. Rather than approximating the full six-dimensional distribution function with a single global tensor train, we employ a collection of TT representations, one at each spatial grid point. The tensor cores are evolved using the projector-splitting integrator introduced in \cite{lubich2015time}, together with carefully designed temporal and spatial discretizations to ensure efficiency.
This localized strategy avoids the rapid TT-rank growth that typically arises from coupling spatial and velocity variables, while still exploiting the strong compressibility of the solution in velocity space. This design is further motivated by kinetic theory: the local equilibrium distribution (the Maxwellian) is fully separable in velocity space and corresponds to a rank-one tensor. Consequently, when the system evolves near equilibrium, the distribution function is expected to remain low rank in the proposed representation. Extensive numerical experiments confirm that the proposed method achieves substantial reductions in memory usage and computational cost while maintaining high accuracy. Finally, we note that a similar decomposition strategy was recently applied in \cite{galindoolarte2025nodal} to a 1D2V setting, using matrix SVD tools and a SAT approach. On the other hand, the dynamical TT approach has been explored and further developed in \cite{DV21,DRV21}, where the Fokker-Planck equation was considered. In the terminology of this paper, this corresponds to the (spatially homogeneous) linear Fokker-Planck equation, without interaction with physical space.

The rest of this paper is organized as follows. In \Cref{sec_TT_operators}, we introduce the basic operations of tensor trains that are relevant to the subsequent discussion.
In \Cref{sec_time_integration_of_tensor_trains}, we describe in detail the projector-splitting method adapted to the TT representation proposed in this work. 
In \Cref{sec_examples}, we apply the proposed method to a series of kinetic equations, including both spatially homogeneous and inhomogeneous cases. In particular, we consider the Vlasov-Amp\`ere-Fokker-Planck equation, a kinetic model widely used to describe plasma dynamics. The paper is concluded in \Cref{sec_future_work}.

\section{Operations of tensor trains}
\label{sec_TT_operators}
In this section, we briefly review some basic operations in the TT format, with a particular focus on the TT representation of three-way tensors.
These operations will be utilized in the subsequent construction of the numerical methods.
While we present the three-dimensional case for simplicity, the underlying principles directly extend to higher-dimensional tensor trains.
Most of these operations are well established in the literature \cite{oseledets2011tensor,lee2018fundamental}, and the remaining ones are straightforward to interpret.

We assume the following two three-way tensors $\bA, \bB\in\mathbb{R}^{N_1\times N_2\times N_3}$ given in TT form:
\begin{equation}
\label{eq_TT_bA} 
    \bA_{k_1k_2k_3} = \sum_{\alpha_1=1}^{r_1} \sum_{\alpha_2=1}^{r_2} 
    A_{k_1\alpha_1}^{(1)}  
    A_{\alpha_1k_2\alpha_2}^{(2)}
    A_{\alpha_2k_3}^{(3)},\quad
    \bB_{k_1k_2k_3} = \sum_{\beta_1=1}^{s_1} \sum_{\beta_2=1}^{s_2} 
    B_{k_1\beta_1}^{(1)}  
    B_{\beta_1k_2\beta_2}^{(2)}
    B_{\beta_2k_3}^{(3)}.
\end{equation}
Here, $A^{(1)}\in\mathbb{R}^{N_1\times r_1},A^{(2)}\in\mathbb{R}^{r_1\times N_2\times r_2},A^{(3)}\in\mathbb{R}^{r_2\times N_3}$ are called the cores of $\bA$.
The TT-rank or bond dimension of $\bA$ is $(r_1,r_2)$. The terms for the tensor $\bB$ are defined analogously.
In principle, any tensor can be represented in TT form as long as the TT-rank is chosen  sufficiently large \cite{oseledets2011tensor}.
However, in many physical applications, the TT-rank can be kept low because of the intrinsic low-rank nature of the underlying physical quantities \cite{khoromskij2018tensor,bachmayr2023low}. Storing the three-way tensor $\bA$ in full form requires $\mathcal{O}(N^3)$ memory with $N = \max\{N_1,N_2,N_3\}$.
In contrast, its TT representation requires only $\mathcal{O}(Nr^2)$ memory, where $r = \max\{r_1, r_2\}$ denotes the maximal TT-rank.
When the TT-rank $(r_1, r_2)$ remains moderate, the TT format yields a substantial reduction in storage cost.
\begin{enumerate}
    \item \textbf{Scaling of a tensor train.} Given a constant $c$, $c\bA$ is a tensor train given by
    \begin{equation}
    \begin{aligned}
        (c\bA)_{k_1k_2k_3} &= \sum_{\alpha_1=1}^{r_1}
        \sum_{\alpha_2=1}^{r_2} \left(cA_{k_1\alpha_1}^{(1)}\right) A_{\alpha_1k_2\alpha_2}^{(2)} A_{\alpha_2k_3}^{(3)} \\
        &=\sum_{\alpha_1=1}^{r_1}
        \sum_{\alpha_2=1}^{r_2} A_{k_1\alpha_1}^{(1)}\left(cA_{\alpha_1k_2\alpha_2}^{(2)}\right)  A_{\alpha_2k_3}^{(3)} =\sum_{\alpha_1=1}^{r_1}
        \sum_{\alpha_2=1}^{r_2} A_{k_1\alpha_1}^{(1)}  A_{\alpha_1k_2\alpha_2}^{(2)}
        \left(cA_{\alpha_2k_3}^{(3)}\right).
    \end{aligned}
    \end{equation}
    That is, multiplying a tensor train by a constant only requires scaling one of its cores.
    \item \textbf{Summation of two general tensor trains.} The summation of $\bA$ and $\bB$ is a tensor train $\mathbf{C}$ with cores given by
    \begin{equation}
    \label{eq_tensor_train_summation}
    \begin{aligned}
        &C_{k_1\gamma_1}^{(1)} = 
        \begin{dcases}
            A_{k_1\gamma_1}^{(1)}, 
            & \text{if~} 1 \leqslant \gamma_1 \leqslant r_1,\\
            B_{k_1(\gamma_1-r_1)}^{(1)},
            & \text{if~} r_1+1 \leqslant \gamma_1 \leqslant r_1+s_1,
        \end{dcases} \\
        &C_{\gamma_1 k_2 \gamma_2}^{(2)} = 
        \begin{dcases}
            A_{\gamma_1 k_2\gamma_2}^{(2)}, 
            & \text{if~} 1 \leqslant \gamma_1 \leqslant r_1, 1\leqslant \gamma_2 \leqslant r_2,\\
            B_{(\gamma_1-r_1)k_2(\gamma_2-r_2)}^{(2)},
            & \text{if~} r_1+1 \leqslant \gamma_1 \leqslant r_1+s_1, r_2+1\leqslant \gamma_2 \leqslant r_2+s_2, \\
            0, & \text{else},
        \end{dcases} \\
        &C_{\gamma_2k_3}^{(3)} = \begin{dcases}
            A_{\gamma_2k_3}^{(3)}, &\text{if~} 1\leqslant \gamma_2 \leqslant r_2, \\
            B_{(\gamma_2-r_2)k_3}^{(3)}, &\text{if~} r_2+1\leqslant \gamma_2 \leqslant r_2+s_2.
        \end{dcases}
    \end{aligned}
    \end{equation}
    That is,  $\mathbf{C}=\bA+\bB$ is a tensor train with TT-rank $(r_1+s_1,r_2+s_2)$.
    \item \textbf{Summation of two tensor trains with a single varying core.} If $r_1=s_1$, $r_2=s_2$, $A_{k_1\alpha_1}^{(1)} = B_{k_1\alpha_1}^{(1)}$ and $A_{\alpha_1k_2\alpha_2}^{(2)} = B_{\alpha_1k_2\alpha_2}^{(2)}$, i.e., $\bA$ and $\bB$ have the same TT-rank and differ only in the third core, then the summation of $\bA$ and $\bB$ can be simply defined as
    \begin{equation}
        (\bA+\bB)_{k_1k_2k_3} = 
        \sum_{\alpha_1=1}^{r_1}
        \sum_{\alpha_2=1}^{r_2}
        A_{k_1\alpha_1}^{(1)} A_{\alpha_1k_2\alpha_2}^{(2)} \left(A_{\alpha_2k_3}^{(3)}+B_{\alpha_2k_3}^{(3)}\right).
    \end{equation}
    The summation can be defined similarly if the two tensor trains have the same TT-rank but differ only in the second or first core.  Compared with the summation of two general tensor trains, summing two tensor trains with a single varying core does not increase the TT-rank so is more efficient.
    \item \textbf{Elementwise inverse of a rank-(1,1) tensor train.} 
    If $r_1=r_2=1$, $\bA$ is a rank-(1,1) tensor train.
    If we further assume that all elements of $\bA$ are non-zero, we can define its elementwise inverse $\frac{1}{\bA}$ in TT form by 
    \begin{equation}
        \left(\frac{1}{\bA}\right)_{k_1k_2k_3}
        = \frac{1}{A_{k_1 1}^{(1)}}
        \frac{1}{A_{1 k_2 1}^{(2)}}
        \frac{1}{A_{1 k_3}^{(3)}},
    \end{equation}
    which is also a rank-(1,1) tensor train.
    The elementwise inverse of a general tensor train does not admit such a simple form; readers may refer to \cite{oseledets2010ttcross} for an approximation algorithm.
    \item \textbf{Hadamard product of a rank-(1,1) tensor train and a general tensor train.}
    If $r_1=r_2=1$, $\bA$ is a rank-(1,1) tensor train, and we define the Hadamard product (elementwise product) $\bA \odot \bB$ of $\bA$ and $\bB$ by 
    \begin{equation}
        (\bA\odot \bB)_{k_1k_2k_3} = \sum_{\beta_1=1}^{s_1}
        \sum_{\beta_2=1}^{s_2}
        \left(A_{k_11}^{(1)}B_{k_1\beta_1}^{(1)}\right)
        \left(A_{1k_21}^{(2)}B_{\beta_1k_2\beta_2}^{(2)}\right)
        \left(A_{1k_3}^{(3)}B_{\beta_2k_3}^{(3)}\right).
    \end{equation}
    The Hadamard product of two general tensor trains also admits a simple representation \cite{lee2018fundamental,sun2024hatt}.
    Since it does not appear in our numerical methods, we omit its discussion here.
    \item \textbf{Summation of all elements.} For a tensor train $\bA$ given in \eqref{eq_TT_bA}, we can sum all its elements by
    \begin{equation}
    \label{eq_TT_sum_all_terms}
        \sum_{k_1=1}^{N_1}
        \sum_{k_2=1}^{N_2}
        \sum_{k_3=1}^{N_3}
        \bA_{k_1k_2k_3} = \sum_{\alpha_1=1}^{r_1}
        \sum_{\alpha_2=1}^{r_2}
        \left(
            \sum_{k_1=1}^{N_1} A_{k_1\alpha_1}^{(1)}
        \right)
        \left(
            \sum_{k_2=1}^{N_2} A_{\alpha_1 k_2 \alpha_2}^{(2)}
        \right)
        \left(
            \sum_{k_3=1}^{N_3} A_{\alpha_2 k_3}^{(3)}
        \right).
    \end{equation}
    While the order of summation does not affect the final result, it can have a significant impact on the computational cost.
    We refer the reader to previous work \cite{oseledets2011tensor} for further discussion.
    In what follows, we apply this method to compute integrals of functions whose values at discrete grid points are represented in the TT format.
     \item \textbf{Orthonormalization I.} The TT representation of a tensor is generally not unique. The orthonormalization process helps construct tensor trains that satisfy orthogonality conditions, serving as a generalization of the QR decomposition for matrices. We start from a three-way tensor train $\bA$ in the form \eqref{eq_TT_bA}.
     The first step is to apply the QR decomposition to $A^{(3)}$ such that
     \begin{equation}
         A_{\alpha_2 k_3}^{(3)} = \sum_{\alpha_2' = 1}^{r_2} R_{\alpha_2\alpha_2'}^{(2)} Q_{\alpha_2'k_3}^{(3)},
         \quad
          \sum_{k_3=1}^{N_3} Q_{\alpha_2k_3}^{(3)}Q_{\alpha_2'k_3}^{(3)}=\delta_{\alpha_2\alpha_2'}.
     \end{equation}
     We then multiply $A^{(2)}$ with $R^{(2)}$ and apply the QR decomposition again:
     \begin{equation}
\sum_{\alpha_2=1}^{r_2}A_{\alpha_1k_2\alpha_2}^{(2)}
         R_{\alpha_2\alpha_2'}^{(2)}
         = \sum_{\alpha_1'=1}^{r_1} R^{(1)}_{\alpha_1\alpha_1'} Q_{\alpha_1'k_2\alpha_2'}^{(2)},
         \quad \sum_{k_2=1}^{N_2} \sum_{\alpha_2'=1}^{r_2}
         Q_{\alpha_1k_2\alpha_2'}^{(2)}
         Q_{\alpha_1'k_2\alpha_2'}^{(2)} = \delta_{\alpha_1\alpha_1'}.
     \end{equation}
     The last step is to multiply $A^{(1)}$ with $R^{(1)}$, yielding a tensor train of the form
     \begin{equation} \label{form1}
\sum_{\alpha_1=1}^{r_1} A^{(1)}_{k_1\alpha_1}R^{(1)}_{\alpha_1 \alpha'_1}:=C_{k_1\alpha'_1}^{(1)}, \quad
         \bA_{k_1k_2k_3} 
         = \sum_{\alpha'_1=1}^{r_1}
         \sum_{\alpha'_2=1}^{r_2}
         C_{k_1\alpha'_1}^{(1)}
         Q_{\alpha'_1k_2\alpha'_2}^{(2)}
         Q_{\alpha'_2k_3}^{(3)}.
     \end{equation}
     Here, $Q^{(2)},Q^{(3)}$ satisfy the orthonormality conditions, and $C^{(1)}$ is the non-orthonormal part of the tensor train. 
     
     In the study of tensors and tensor trains, the cumbersome subscripts and summations can be represented more intuitively using tensor diagram notation \cite{penrose1971applications,schollwock2011density,lubich2015time,paeckel2019time}. Diagrammatically, the above process can be expressed as
\begin{equation}
\begin{tikzpicture}[baseline=0]
\draw (0,0) -- (1.6,0);
\draw (1.6,0) -- (3.2,0);
\draw (0,0) -- (0,0.8);
\draw (1.6,0) -- (1.6,0.8);
\draw (3.2,0) -- (3.2,0.8);
\draw[fill=white,line width = 1.2] (0,0) circle (0.2);
\draw[fill=white,line width = 1.2] (1.6,0) circle (0.2);
\draw[fill=white,line width = 1.2] (3.2,0) circle (0.2);
\node at (0,-0.5) {$A^{(1)}$};
\node at (1.6,-0.5) {$A^{(2)}$};
\node at (3.2,-0.5) {$A^{(3)}$};
\node at (0.8,0.2) {$\alpha_1$};
\node at (2.4,0.2) {$\alpha_2$};
\node at (0,1) {$k_1$};
\node at (1.6,1) {$k_2$};
\node at (3.2,1) {$k_3$};
\end{tikzpicture}
=
\begin{tikzpicture}[baseline=0]
\draw (0,0) -- (1.6,0);
\draw (1.6,0) -- (3.2,0);
\draw (0,0) -- (0,0.8);
\draw (1.6,0) -- (1.6,0.8);
\draw (3.2,0) -- (3.2,0.8);
\draw[fill=white,line width = 1.2] (0,0) circle (0.2);
\draw[fill=white,line width = 1.2] (1.4,0) -- (1.7,-0.2) -- (1.7,0.2) -- (1.4,0);
\draw[fill=white,line width = 1.2] (3.0,0) -- (3.3,-0.2) -- (3.3,0.2) -- (3.0,0);
\node at (0,-0.5) {$C^{(1)}$};
\node at (1.6,-0.5) {$Q^{(2)}$};
\node at (3.2,-0.5) {$Q^{(3)}$};
\node at (0.8,0.2) {$\alpha_1$};
\node at (2.4,0.2) {$\alpha_2$};
\node at (0,1) {$k_1$};
\node at (1.6,1) {$k_2$};
\node at (3.2,1) {$k_3$};
\end{tikzpicture}
\end{equation}
In the left diagram, the first circle has two legs, $k_1$ and $\alpha_1$, indicating that it is a two-way tensor (matrix). The second circle has three legs and represents a three-way tensor. The last circle again has two legs, representing a matrix. The $\alpha_1$ leg of $A^{(1)}$ is connected to the $\alpha_1$ leg of $A^{(2)}$, meaning that these two indices are contracted. Mathematically, it corresponds to the summation over $\alpha_1$ in \eqref{eq_TT_bA}. The contraction of index $\alpha_2$ is similar. After contractions, the resulting diagram has three free legs $k_1,k_2,k_3$, indicating that the final object is a three-way tensor with these physical indices. The right diagram is interpreted in the same way, except that we use a circle to represent a general tensor core and a triangle to represent a core satisfying an orthonormality condition. The orientation of the triangle intuitively indicates the direction of the orthonormality property.

The orthonormalization procedure can also be applied to full tensors, in which case it is commonly referred to as the TT-SVD algorithm \cite{oseledets2011tensor}. Since in all our numerical examples the initial conditions can be represented directly as tensor trains without forming the full tensors, we omit the details of the TT-SVD algorithm.

\item \textbf{Orthonormalization II.} 
    If we assume a tensor train is already in the form given on the right-hand side of \eqref{form1}, we can perform another orthonormalization procedure by first applying the QR decomposition to $C^{(1)}$ such that
    \begin{equation}
         C_{k_1 \alpha_1'}^{(1)} = \sum_{\alpha_1 = 1}^{r_1} P_{k_1\alpha_1}^{(1)}
    S_{\alpha_1\alpha'_1}^{(1)},
         \quad
          \sum_{k_1=1}^{N_1} 
    P_{k_1\alpha_1}^{(1)}
    P_{k_1\alpha_1'}^{(1)} 
    = \delta_{\alpha_1\alpha_1'}.
     \end{equation}
     We then multiply $S^{(1)}$ with $Q^{(2)}$ and apply the QR decomposition again:
     \begin{equation}
\sum_{\alpha'_1=1}^{r_1}S_{\alpha_1\alpha_1'}^{(1)}
    Q_{\alpha'_1k_2\alpha'_2}^{(2)}
         :=C_{\alpha_1 k_2\alpha_2'}^{(2)} = \sum_{\alpha_2=1}^{r_2} P_{\alpha_1 k_2 \alpha_2}^{(2)}
    S_{\alpha_2 \alpha_2'}^{(2)} ,
         \quad \sum_{k_2=1}^{N_2} 
    \sum_{\alpha_1=1}^{r_1}
    P_{\alpha_1 k_2 \alpha_2}^{(2)}
    P_{\alpha_1 k_2 \alpha_2'}^{(2)}
    = \delta_{\alpha_2\alpha_2'}.
     \end{equation}
     Finally, multiplying $S^{(2)}$ with $Q^{(3)}$ yields a tensor train of the form
      \begin{equation} 
\sum_{\alpha_2'=1}^{r_2} S_{\alpha_2 \alpha_2'}^{(2)} 
    Q_{\alpha_2' k_3}^{(3)}:=C_{\alpha_2 k_3}^{(3)}, \quad \sum_{\alpha'_1=1}^{r_1}
         \sum_{\alpha'_2=1}^{r_2}
         C_{k_1\alpha'_1}^{(1)}
         Q_{\alpha'_1k_2\alpha'_2}^{(2)}
         Q_{\alpha'_2k_3}^{(3)}=\sum_{\alpha_1=1}^{r_1}\sum_{\alpha_2=1}^{r_2} P^{(1)}_{k_1\alpha_1}P^{(2)}_{\alpha_1k_2\alpha_2}C^{(3)}_{\alpha_2k_3}.
     \end{equation}
          Here, $P^{(1)},P^{(2)}$ satisfy the orthonormality conditions, and $C^{(3)}$ is the non-orthonormal part of the tensor train. Diagrammatically, the above process can be expressed as
\begin{equation}
\begin{tikzpicture}[baseline=0]
\draw (0,0) -- (1.6,0);
\draw (1.6,0) -- (3.2,0);
\draw (0,0) -- (0,0.8);
\draw (1.6,0) -- (1.6,0.8);
\draw (3.2,0) -- (3.2,0.8);
\draw[fill=white,line width = 1.2] (0,0) circle (0.2);
\draw[fill=white,line width = 1.2] (1.4,0) -- (1.7,-0.2) -- (1.7,0.2) -- (1.4,0);
\draw[fill=white,line width = 1.2] (3.0,0) -- (3.3,-0.2) -- (3.3,0.2) -- (3.0,0);
\node at (0,-0.5) {$C^{(1)}$};
\node at (1.6,-0.5) {$Q^{(2)}$};
\node at (3.2,-0.5) {$Q^{(3)}$};
\node at (0.8,0.2) {$\alpha_1$};
\node at (2.4,0.2) {$\alpha_2$};
\node at (0,1) {$k_1$};
\node at (1.6,1) {$k_2$};
\node at (3.2,1) {$k_3$};
\end{tikzpicture}
=
\begin{tikzpicture}[baseline=0]
\draw (0,0) -- (1.6,0);
\draw (1.6,0) -- (3.2,0);
\draw (0,0) -- (0,0.8);
\draw (1.6,0) -- (1.6,0.8);
\draw (3.2,0) -- (3.2,0.8);
\draw[fill=white,line width = 1.2] (0.2,0) -- (-0.1,-0.2) -- (-0.1,0.2) -- (0.2,0);
\draw[fill=white,line width = 1.2] (1.8,0) -- (1.5,-0.2) -- (1.5,0.2) -- (1.8,0);
\draw[fill=white,line width = 1.2] (3.2,0) circle (0.2);
\node at (0,-0.5) {$P^{(1)}$};
\node at (1.6,-0.5) {$P^{(2)}$};
\node at (3.2,-0.5) {$C^{(3)}$};
\node at (0.8,0.2) {$\alpha_1$};
\node at (2.4,0.2) {$\alpha_2$};
\node at (0,1) {$k_1$};
\node at (1.6,1) {$k_2$};
\node at (3.2,1) {$k_3$};
\end{tikzpicture}
\end{equation}
The orthonormalization in the opposite direction can be done similarly.
\end{enumerate}

\section{Time integration of tensor trains}
\label{sec_time_integration_of_tensor_trains}
In this section, we present our method for a general kinetic equation of the form
\begin{equation}
\label{eq_original_equation}
    \partial_t f(t,x,\bv) = F(f), \quad x \in \Omega_x\subset\mathbb{R}, \quad \bv=(v^{(1)},v^{(2)},v^{(3)}) \in \mathbb{R}^3.
\end{equation}
The right-hand-side function $F(f)$ is problem-dependent and will be specified in the following sections. Note that for the velocity space, we consider the (physically relevant) three-dimensional case, whereas for the physical space we restrict to one dimension for simplicity. The presentation can be extended to multiple spatial dimensions in a straightforward manner.

Our method is based on the projector-splitting integrator introduced in \cite{lubich2015time}, implemented using a discretize-then-project (DtP) approach. For the kinetic equation \eqref{eq_original_equation}, this means that we first discretize both the physical space and the velocity space, and then project the resulting ODE system onto the tangent space of the low-rank solution manifold. Although recent studies \cite{kusch2023stability,zhang2025stability} have reported restrictive time-step conditions for the DtP approach compared to the project-then-discretize (PtD) approach -- where the original equation is first projected onto the low-rank space and the resulting PDE is then discretized -- the DtP approach offers a significant advantage in avoiding complicated projections of the governing equation, particularly in high-dimensional settings.

We first introduce the spatial discretization by employing the finite difference method. In the physical space, we assume grid points $\{x_j\}_{j=1}^{N_x}\in \Omega_x$, whereas in the velocity space we choose grid points $\{v_{k_1},v_{k_2},v_{k_3}\}_{k_1,k_2,k_3=1}^{N_v}\in \mathbb{R}^3$. For simplicity, we assume uniform grids and the same number of grid points in each velocity dimension, although our method can be extended to more general settings. Additionally, we introduce the temporal discretization $t_n = n\dt$ with $n=0,1,2,\cdots$, and define the discrete values of $f$ on the grid points at time $t_n$ by
\begin{equation}
    f_{k_1k_2k_3}^{j,n} \approx f(t_n,x_j,v_{k_1},v_{k_2},v_{k_3}).
\end{equation}
In our framework, for each fixed time $t_n$ and spatial grid point $x_j$, we represent the array $f_{k_1k_2k_3}^{j,n}$ as a three-way tensor train in the velocity space, denoted $\bf^{j,n}$.

With these notations, assuming an explicit time-stepping scheme, equation \eqref{eq_original_equation} is discretized as
\begin{equation}
    \bf^{j,n+1} = \mathcal{F}_{\Delta t}^{j}(\bf^{1,n},\cdots,\bf^{N_x,n}), \quad j=1,\dots,N_x.
\end{equation}
A natural idea for evolving the above system is to evaluate the right-hand side $\mathcal{F}_{\dt}^j(\bf^{1,n},\cdots,\bf^{N_x,n})$ in TT format and assign the resulting tensor train to $\bf^{j,n+1}$. However, this process generally increases the TT-rank, causing the memory cost to grow rapidly.
When the cost becomes unaffordable, TT-rounding \cite{oseledets2011tensor} must be applied to reduce the memory footprint. This forms the basis of the SAT approach. Despite its simplicity, the computational cost and memory requirement of this approach can still be quite high.

In contrast to the SAT approach, our method is based on the projector-splitting integrator \cite{lubich2014projector,lubich2015time}. Although more intrusive, this method fixes the rank of the solution throughout the time evolution (with rank adaptivity easily incorporated if needed), and is therefore more computationally efficient and memory friendly. At each time step, the method involves several forward and backward sub-projections to update the tensor cores of the solution.

In the following, we present a first-order projector-splitting scheme for equation \eqref{eq_original_equation} in complete detail. Our presentation is at the fully discrete level and employs a notation that is very different from (and, we hope, more intuitive than) that in \cite{lubich2015time}. To facilitate the exposition, we assume that the update in the forward sub-projection step follows the scheme
\begin{equation}
\label{eq_distribution_propagation}
    \bf^{j,*} = \mathcal{G}_{\Delta t}^{j}(\bf^{1},\cdots,\bf^{N_x}), \quad j=1,\dots,N_x,
\end{equation}
and that the update in the backward sub-projection step follows the scheme
\begin{equation}
\label{eq_distribution_backward_propagation}
    \bf^{j,*} = \mathcal{H}_{\Delta t}^j(\bf^{1},\cdots,\bf^{N_x}), \quad j=1,\dots,N_x.
\end{equation}
The mappings $\mathcal{G}_{\Delta t}^{j}$ and $\mathcal{H}_{\Delta t}^{j}$ arise from suitable first-order time-stepping schemes for \eqref{eq_original_equation}, which will be specified in the following section when we consider concrete equations. We will also need the following five forms of the TT representation:
\begin{gather}
\begin{tikzpicture}[baseline=0]
\draw (0,0) -- (1.6,0);
\draw (1.6,0) -- (3.2,0);
\draw (0,0) -- (0,0.8);
\draw (1.6,0) -- (1.6,0.8);
\draw (3.2,0) -- (3.2,0.8);
\draw[fill=white,line width = 1.2] (0,0) circle (0.2);
\draw[fill=white,line width = 1.2] (1.4,0) -- (1.7,-0.2) -- (1.7,0.2) -- (1.4,0);
\draw[fill=white,line width = 1.2] (3.0,0) -- (3.3,-0.2) -- (3.3,0.2) -- (3.0,0);
\node at (0,-0.5) {$C^{j,(1)}$};
\node at (1.6,-0.5) {$Q^{j,(2)}$};
\node at (3.2,-0.5) {$Q^{j,(3)}$};
\node at (0.8,0.2) {$\alpha_1$};
\node at (2.4,0.2) {$\alpha_2$};
\node at (0,1) {$k_1$};
\node at (1.6,1) {$k_2$};
\node at (3.2,1) {$k_3$};
\node at (1.6,-1) {\rombracket{1}};
\end{tikzpicture}
\begin{tikzpicture}[baseline=0]
\draw (0,0) -- (1.6,0);
\draw (1.6,0) -- (3.2,0);
\draw (0,0) -- (0,0.8);
\draw (1.6,0) -- (1.6,0.8);
\draw (3.2,0) -- (3.2,0.8);
\draw[fill=white,line width = 1.2] (0.2,0) -- (-0.1,-0.2) -- (-0.1,0.2) -- (0.2,0);
\draw[fill=white,line width = 1.2] (0.8,0) circle (0.2);
\draw[fill=white,line width = 1.2] (1.4,0) -- (1.7,-0.2) -- (1.7,0.2) -- (1.4,0);
\draw[fill=white,line width = 1.2] (3.0,0) -- (3.3,-0.2) -- (3.3,0.2) -- (3.0,0);
\node at (-0.2,-0.5) {$P^{j,(1)}$};
\node at (0.8,-0.5) {$S^{j,(1)}$};
\node at (1.8,-0.5) {$Q^{j,(2)}$};
\node at (3.2,-0.5) {$Q^{j,(3)}$};
\node at (0.4,0.2) {$\alpha_1$};
\node at (1.2,0.2) {$\alpha_1'$};
\node at (2.4,0.2) {$\alpha_2$};
\node at (0,1) {$k_1$};
\node at (1.6,1) {$k_2$};
\node at (3.2,1) {$k_3$};
\node at (1.6,-1) {\rombracket{2}};
\end{tikzpicture} \nonumber\\
\begin{tikzpicture}[baseline=0]
\draw (0,0) -- (1.6,0);
\draw (1.6,0) -- (3.2,0);
\draw (0,0) -- (0,0.8);
\draw (1.6,0) -- (1.6,0.8);
\draw (3.2,0) -- (3.2,0.8);
\draw[fill=white,line width = 1.2] (0.2,0) -- (-0.1,-0.2) -- (-0.1,0.2) -- (0.2,0);
\draw[fill=white,line width = 1.2] (1.6,0) circle (0.2);
\draw[fill=white,line width = 1.2] (3.0,0) -- (3.3,-0.2) -- (3.3,0.2) -- (3.0,0);
\node at (0,-0.5) {$P^{j,(1)}$};
\node at (1.6,-0.5) {$C^{j,(2)}$};
\node at (3.2,-0.5) {$Q^{j,(3)}$};
\node at (0.8,0.2) {$\alpha_1$};
\node at (2.4,0.2) {$\alpha_2$};
\node at (0,1) {$k_1$};
\node at (1.6,1) {$k_2$};
\node at (3.2,1) {$k_3$};
\node at (1.6,-1) {\rombracket{3}};
\end{tikzpicture} 
\begin{tikzpicture}[baseline=0]
\draw (0,0) -- (1.6,0);
\draw (1.6,0) -- (3.2,0);
\draw (0,0) -- (0,0.8);
\draw (1.6,0) -- (1.6,0.8);
\draw (3.2,0) -- (3.2,0.8);
\draw[fill=white,line width = 1.2] (0.2,0) -- (-0.1,-0.2) -- (-0.1,0.2) -- (0.2,0);
\draw[fill=white,line width = 1.2] (1.8,0) -- (1.5,-0.2) -- (1.5,0.2) -- (1.8,0);
\draw[fill=white,line width = 1.2] (2.4,0) circle (0.2);
\draw[fill=white,line width = 1.2] (3.0,0) -- (3.3,-0.2) -- (3.3,0.2) -- (3.0,0);
\node at (0,-0.5) {$P^{j,(1)}$};
\node at (1.4,-0.5) {$P^{j,(2)}$};
\node at (2.4,-0.5) {$S^{j,(2)}$};
\node at (3.4,-0.5) {$Q^{j,(3)}$};
\node at (0.8,0.2) {$\alpha_1$};
\node at (2.0,0.2) {$\alpha_2'$};
\node at (2.8,0.2) {$\alpha_2$};
\node at (0,1) {$k_1$};
\node at (1.6,1) {$k_2$};
\node at (3.2,1) {$k_3$};
\node at (1.6,-1) {\rombracket{4}};
\end{tikzpicture}
\begin{tikzpicture}[baseline=0]
\draw (0,0) -- (1.6,0);
\draw (1.6,0) -- (3.2,0);
\draw (0,0) -- (0,0.8);
\draw (1.6,0) -- (1.6,0.8);
\draw (3.2,0) -- (3.2,0.8);
\draw[fill=white,line width = 1.2] (0.2,0) -- (-0.1,-0.2) -- (-0.1,0.2) -- (0.2,0);
\draw[fill=white,line width = 1.2] (1.8,0) -- (1.5,-0.2) -- (1.5,0.2) -- (1.8,0);
\draw[fill=white,line width = 1.2] (3.2,0) circle (0.2);
\node at (0,-0.5) {$P^{j,(1)}$};
\node at (1.6,-0.5) {$P^{j,(2)}$};
\node at (3.2,-0.5) {$C^{j,(3)}$};
\node at (0.8,0.2) {$\alpha_1$};
\node at (2.4,0.2) {$\alpha_2$};
\node at (0,1) {$k_1$};
\node at (1.6,1) {$k_2$};
\node at (3.2,1) {$k_3$};
\node at (1.6,-1) {\rombracket{5}};
\end{tikzpicture}
\label{eq_TT_decomposition_forms1}
\end{gather}
The algorithm proceeds as follows:
\begin{itemize}
    \item \textbf{Step 0:} Based on the initial condition, at each spatial point $x_j$, form the initial tensor train with TT-rank $(r_1,r_2)$, and convert it into form \rombracket{1} in \eqref{eq_TT_decomposition_forms1}. Denote the resulting tensor train by $\bf_\rom{1}^{j,0}$. This step follows the \textbf{operation 7} described in Section~\ref{sec_TT_operators}.
    Ideally, one should avoid constructing full tensors and then converting them into tensor trains, as storing full tensors is computationally expensive. 
     \end{itemize}
    Now assume that we have $\{\bf^{j,n}_\rom{1}\}_{j=1}^{N_x}$ in form \rombracket{1} at time $t_n$. The following steps aim to compute $\{\bf_\rom{1}^{j,n+1}\}_{j=1}^{N_x}$ in form \rombracket{1} at $t_{n+1}$.
    The update is split into five steps, where \textbf{steps 1, 3, and 5} are forward steps and \textbf{steps 2 and 4} are backward steps. 
    At each spatial point $x_j$,
    \begin{itemize}
    \item \textbf{Step 1a:} Compute $\bG^{j,(1)}:= \mathcal{G}_{\Delta t}^{j}(\bf_\rom{1}^{1,n},\cdots,\bf_\rom{1}^{N_x,n})$ using the scheme \eqref{eq_distribution_propagation}. 
    \item \textbf{Step 1b:}
    Compute the projection 
    $\tilde{C}^{j,(1)}: = \left\langle 
        \bG^{j,(1)}, Q^{j,(2)}, Q^{j,(3)}      \right\rangle_\rom{1}$. Componentwise, this is
    \begin{equation}
        \tilde{C}_{k_1\alpha_1}^{j,(1)}
        = \sum_{k_2=1}^{N_v}
        \sum_{k_3=1}^{N_v}
        \sum_{\alpha_2=1}^{r_2}
        \bG_{k_1k_2k_3}^{j,(1)}
        Q_{\alpha_1k_2\alpha_2}^{j,(2)}
        Q_{\alpha_2k_3}^{j,(3)},
    \end{equation}
    or diagrammatically,
    \begin{equation}
    \label{eq_substep_1a_diagram}
    \begin{tikzpicture}[baseline=0]
        \draw (0.7,0) -- (0,0) -- (0,0.7);
        \draw[fill=white,line width = 1.2] (0,0) circle (0.2);
        \node at (0,-0.5) {$\tilde{C}^{j,(1)}$};
        \node at (0,0.9) {$k_1$};
        \node at (0.6,0.2) {$\alpha_1$};
    \end{tikzpicture}
    =
    \begin{tikzpicture}[baseline=0]
        \draw (0,0) -- (0, 1.4);
        \draw (1.6,0) -- (1.6, 1.4);
        \draw (3.2,0) -- (3.2, 1.4);
        \draw (0,0) -- (3.2,0);
        \draw (0,1.4) -- (3.2,1.4);
        \draw[rounded corners=2pt, dashed, gray] (-0.4,1) rectangle (3.6,1.8);
        \draw[fill=white,line width = 1.2] (0,1.4) circle (0.2);
        \draw[fill=white,line width = 1.2] (1.6,1.4) circle (0.2);
        \draw[fill=white,line width = 1.2] (3.2,1.4) circle (0.2);
        \draw[fill=white,line width = 1.2] (1.4,0) -- (1.7,-0.2) -- (1.7,0.2) -- (1.4,0);
        \draw[fill=white,line width = 1.2] (3.0,0) -- (3.3,-0.2) -- (3.3,0.2) -- (3.0,0);
        \draw[fill=white,line width = 1.2, draw=white] (-0.5,-0.5) -- (-0.5,0.5) -- (0.5,0.5) -- (0.5,-0.5) -- (-0.5,-0.5);
        \node at (1.6,2.1) {$\bG^{j,(1)}$};
        \node at (1.6,-0.5) {$Q^{j,(2)}$};
        \node at (3.2,-0.5) {$Q^{j,(3)}$};
        \node at (0.2,0.7) {$k_1$};
        \node at (1.8,0.7) {$k_2$};
        \node at (3.4,0.7) {$k_3$};
        \node at (1,0.2) {$\alpha_1$};
        \node at (2.4,0.2) {$\alpha_2$};
    \end{tikzpicture}
    \end{equation}
    \item \textbf{Step 1c:} Replace $C^{j,(1)}$ in $\bf_\rom{1}^{j,n}$ by $\tilde{C}^{j,(1)}$ to obtain $\tilde{\bf}_\rom{1}^j$:
    \begin{equation}
        \begin{tikzpicture}[baseline=0]
        \draw (0,0) -- (1.6,0);
        \draw (1.6,0) -- (3.2,0);
        \draw (0,0) -- (0,0.8);
        \draw (1.6,0) -- (1.6,0.8);
        \draw (3.2,0) -- (3.2,0.8);
        \draw[fill=white,line width = 1.2] (0,0) circle (0.2);
        \draw[fill=white,line width = 1.2] (1.4,0) -- (1.7,-0.2) -- (1.7,0.2) -- (1.4,0);
        \draw[fill=white,line width = 1.2] (3.0,0) -- (3.3,-0.2) -- (3.3,0.2) -- (3.0,0);
        \node at (0,-0.5) {$C^{j,(1)}$};
        \node at (1.6,-0.5) {$Q^{j,(2)}$};
        \node at (3.2,-0.5) {$Q^{j,(3)}$};
        \node at (0.8,0.2) {$\alpha_1$};
        \node at (2.4,0.2) {$\alpha_2$};
        \node at (0,1) {$k_1$};
        \node at (1.6,1) {$k_2$};
        \node at (3.2,1) {$k_3$};
        \node at (1.6,-1) {$\bf_\rom{1}^{j,n}$};
        \end{tikzpicture}
        \rightarrow
        \begin{tikzpicture}[baseline=0]
        \draw (0,0) -- (1.6,0);
        \draw (1.6,0) -- (3.2,0);
        \draw (0,0) -- (0,0.8);
        \draw (1.6,0) -- (1.6,0.8);
        \draw (3.2,0) -- (3.2,0.8);
        \draw[fill=white,line width = 1.2] (0,0) circle (0.2);
        \draw[fill=white,line width = 1.2] (1.4,0) -- (1.7,-0.2) -- (1.7,0.2) -- (1.4,0);
        \draw[fill=white,line width = 1.2] (3.0,0) -- (3.3,-0.2) -- (3.3,0.2) -- (3.0,0);
        \node at (0,-0.5) {\textcolor{red}{$\tilde{C}^{j,(1)}$}};
        \node at (1.6,-0.5) {$Q^{j,(2)}$};
        \node at (3.2,-0.5) {$Q^{j,(3)}$};
        \node at (0.8,0.2) {$\alpha_1$};
        \node at (2.4,0.2) {$\alpha_2$};
        \node at (0,1) {$k_1$};
        \node at (1.6,1) {$k_2$};
        \node at (3.2,1) {$k_3$};
        \node at (1.6,-1) {$\tilde{\bf}_\rom{1}^j$};
        \end{tikzpicture}
    \end{equation}
    \item \textbf{Step 1d:} Write $\tilde{\bf}_\rom{1}^{j}$ in form \rombracket{2} by performing the QR decomposition on $\tilde{C}^{j,(1)}$, and denote the resulting tensor train by $\bf_\rom{2}^j$:
    \begin{equation}
    \begin{tikzpicture}[baseline=0]
    \draw (0,0) -- (1.6,0);
    \draw (1.6,0) -- (3.2,0);
    \draw (0,0) -- (0,0.8);
    \draw (1.6,0) -- (1.6,0.8);
    \draw (3.2,0) -- (3.2,0.8);
    \draw[fill=white,line width = 1.2] (0,0) circle (0.2);
    \draw[fill=white,line width = 1.2] (1.4,0) -- (1.7,-0.2) -- (1.7,0.2) -- (1.4,0);
    \draw[fill=white,line width = 1.2] (3.0,0) -- (3.3,-0.2) -- (3.3,0.2) -- (3.0,0);
    \node at (0,-0.5) {$\tilde{C}^{j,(1)}$};
    \node at (1.6,-0.5) {$Q^{j,(2)}$};
    \node at (3.2,-0.5) {$Q^{j,(3)}$};
    \node at (0.8,0.2) {$\alpha_1$};
    \node at (2.4,0.2) {$\alpha_2$};
    \node at (0,1) {$k_1$};
    \node at (1.6,1) {$k_2$};
    \node at (3.2,1) {$k_3$};
    \node at (1.6,-1) {$\tilde{\bf}_\rom{1}^j$};
    \end{tikzpicture}
    =
    \begin{tikzpicture}[baseline=0]
    \draw (0,0) -- (1.6,0);
    \draw (1.6,0) -- (3.2,0);
    \draw (0,0) -- (0,0.8);
    \draw (1.6,0) -- (1.6,0.8);
    \draw (3.2,0) -- (3.2,0.8);
    \draw[fill=white,line width = 1.2] (0.2,0) -- (-0.1,-0.2) -- (-0.1,0.2) -- (0.2,0);
    \draw[fill=white,line width = 1.2] (0.8,0) circle (0.2);
    \draw[fill=white,line width = 1.2] (1.4,0) -- (1.7,-0.2) -- (1.7,0.2) -- (1.4,0);
    \draw[fill=white,line width = 1.2] (3.0,0) -- (3.3,-0.2) -- (3.3,0.2) -- (3.0,0);
    \node at (-0.2,-0.5) {$P^{j,(1)}$};
    \node at (0.8,-0.5) {$S^{j,(1)}$};
    \node at (1.8,-0.5) {$Q^{j,(2)}$};
    \node at (3.2,-0.5) {$Q^{j,(3)}$};
    \node at (0.4,0.2) {$\alpha_1$};
    \node at (1.2,0.2) {$\alpha_1'$};
    \node at (2.4,0.2) {$\alpha_2$};
    \node at (0,1) {$k_1$};
    \node at (1.6,1) {$k_2$};
    \node at (3.2,1) {$k_3$};
    \node at (1.6,-1) {${\bf}_\rom{2}^j$};
    \end{tikzpicture}
    \end{equation}
    
    \item \textbf{Step 2a:} Compute $\bH^{j,(1)}: = \mathcal{H}_{\Delta t}^j(\tilde{\bf}_\rom{1}^1,\cdots,\tilde{\bf}_\rom{1}^{N_x})$ using the scheme \eqref{eq_distribution_backward_propagation}.
    \item \textbf{Step 2b:} Compute the projection 
    $
        \tilde{S}^{j,(1)}
        := \left\langle
            P^{j,(1)},
            \bH^{j,(1)},
            Q^{j,(2)},
            Q^{j,(3)}
        \right\rangle_\rom{2}
    $.
    Componentwise, this is
    \begin{equation}
    \label{eq_step2b}
        \tilde{S}_{\alpha_1\alpha_1'}^{j,(1)}
        =\sum_{k_1=1}^{N_v}\sum_{k_2=1}^{N_v} \sum_{k_3=1}^{N_v} \sum_{\alpha_2=1}^{r_2} 
        P^{j,(1)}_{k_1\alpha_1}
        \bH_{k_1k_2k_3}^{j,(1)} Q_{\alpha_1'k_2\alpha_2}^{j,(2)} Q_{\alpha_2k_3}^{j,(3)},
    \end{equation}
    or diagrammatically,
    \begin{equation}
    \begin{tikzpicture}[baseline=0]
        \draw (-0.7,0) -- (0.7,0);
        \draw[fill=white,line width = 1.2] (0,0) circle (0.2);
        \node at (0,-0.5) {$\tilde{S}^{j,(1)}$};
        \node at (-0.6,0.2) {$\alpha_1$};
        \node at (0.6,0.2) {$\alpha_1'$};
    \end{tikzpicture}
    =
    \begin{tikzpicture}[baseline=0]
        \draw (0,0) -- (0, 1.4);
        \draw (1.6,0) -- (1.6, 1.4);
        \draw (3.2,0) -- (3.2, 1.4);
        \draw (0,0) -- (3.2,0);
        \draw (0,1.4) -- (3.2,1.4);
        \draw[rounded corners=2pt, dashed, gray] (-0.4,1) rectangle (3.6,1.8);
        \draw[fill=white,line width = 1.2] (0,1.4) circle (0.2);
        \draw[fill=white,line width = 1.2] (1.6,1.4) circle (0.2);
        \draw[fill=white,line width = 1.2] (3.2,1.4) circle (0.2);
        \draw[fill=white,line width = 1.2] (0.2,0) -- (-0.1,-0.2) -- (-0.1,0.2) -- (0.2,0);
        \draw[fill=white,line width = 1.2] (1.4,0) -- (1.7,-0.2) -- (1.7,0.2) -- (1.4,0);
        \draw[fill=white,line width = 1.2] (3.0,0) -- (3.3,-0.2) -- (3.3,0.2) -- (3.0,0);
        \draw[fill=white,line width = 1.2, draw=white] (0.6,-0.5) -- (0.6,0.5) -- (1.0,0.5) -- (1.0,-0.5) -- (0.6,-0.5);
        \node at (1.6,2.1) {$\bH^{j,(1)}$};
        \node at (0,-0.5) {$P^{j,(1)}$};
        \node at (1.6,-0.5) {$Q^{j,(2)}$};
        \node at (3.2,-0.5) {$Q^{j,(3)}$};
        \node at (0.2,0.7) {$k_1$};
        \node at (1.8,0.7) {$k_2$};
        \node at (3.4,0.7) {$k_3$};
        \node at (0.4,0.2) {$\alpha_1$};
        \node at (1.2,0.2) {$\alpha_1'$};
        \node at (2.4,0.2) {$\alpha_2$};
    \end{tikzpicture}
    \end{equation}
    \item \textbf{Step 2c:} Replace $S^{j,(1)}$ in $\bf_\rom{2}^j$ by $\tilde{S}^{j,(1)}$ to obtain $\tilde{\bf}_\rom{2}^j$:
    \begin{equation}
    \begin{tikzpicture}[baseline=0]
    \draw (0,0) -- (1.6,0);
    \draw (1.6,0) -- (3.2,0);
    \draw (0,0) -- (0,0.8);
    \draw (1.6,0) -- (1.6,0.8);
    \draw (3.2,0) -- (3.2,0.8);
    \draw[fill=white,line width = 1.2] (0.2,0) -- (-0.1,-0.2) -- (-0.1,0.2) -- (0.2,0);
    \draw[fill=white,line width = 1.2] (0.8,0) circle (0.2);
    \draw[fill=white,line width = 1.2] (1.4,0) -- (1.7,-0.2) -- (1.7,0.2) -- (1.4,0);
    \draw[fill=white,line width = 1.2] (3.0,0) -- (3.3,-0.2) -- (3.3,0.2) -- (3.0,0);
    \node at (-0.2,-0.5) {$P^{j,(1)}$};
    \node at (0.8,-0.5) {$S^{j,(1)}$};
    \node at (1.8,-0.5) {$Q^{j,(2)}$};
    \node at (3.2,-0.5) {$Q^{j,(3)}$};
    \node at (0.4,0.2) {$\alpha_1$};
    \node at (1.2,0.2) {$\alpha_1'$};
    \node at (2.4,0.2) {$\alpha_2$};
    \node at (0,1) {$k_1$};
    \node at (1.6,1) {$k_2$};
    \node at (3.2,1) {$k_3$};
    \node at (1.6,-1) {${\bf}_\rom{2}^j$};
    \end{tikzpicture}
    \rightarrow
    \begin{tikzpicture}[baseline=0]
    \draw (0,0) -- (1.6,0);
    \draw (1.6,0) -- (3.2,0);
    \draw (0,0) -- (0,0.8);
    \draw (1.6,0) -- (1.6,0.8);
    \draw (3.2,0) -- (3.2,0.8);
    \draw[fill=white,line width = 1.2] (0.2,0) -- (-0.1,-0.2) -- (-0.1,0.2) -- (0.2,0);
    \draw[fill=white,line width = 1.2] (0.8,0) circle (0.2);
    \draw[fill=white,line width = 1.2] (1.4,0) -- (1.7,-0.2) -- (1.7,0.2) -- (1.4,0);
    \draw[fill=white,line width = 1.2] (3.0,0) -- (3.3,-0.2) -- (3.3,0.2) -- (3.0,0);
    \node at (-0.2,-0.5) {$P^{j,(1)}$};
    \node at (0.8,-0.5) {\textcolor{red}{$\tilde{S}^{j,(1)}$}};
    \node at (1.8,-0.5) {$Q^{j,(2)}$};
    \node at (3.2,-0.5) {$Q^{j,(3)}$};
    \node at (0.4,0.2) {$\alpha_1$};
    \node at (1.2,0.2) {$\alpha_1'$};
    \node at (2.4,0.2) {$\alpha_2$};
    \node at (0,1) {$k_1$};
    \node at (1.6,1) {$k_2$};
    \node at (3.2,1) {$k_3$};
    \node at (1.6,-1) {$\tilde{\bf}_\rom{2}^j$};
    \end{tikzpicture}
    \end{equation}
    \item \textbf{Step 2d:} Write $\tilde{\bf}_\rom{2}^j$ in form \rombracket{3} by multiplying $\tilde{S}^{j,(1)}$ with $Q^{j,(2)}$ , and denote the resulting tensor train by $\bf_\rom{3}^j$:
    \begin{equation}
    \begin{tikzpicture}[baseline=0]
    \draw (0,0) -- (1.6,0);
    \draw (1.6,0) -- (3.2,0);
    \draw (0,0) -- (0,0.8);
    \draw (1.6,0) -- (1.6,0.8);
    \draw (3.2,0) -- (3.2,0.8);
    \draw[fill=white,line width = 1.2] (0.2,0) -- (-0.1,-0.2) -- (-0.1,0.2) -- (0.2,0);
    \draw[fill=white,line width = 1.2] (0.8,0) circle (0.2);
    \draw[fill=white,line width = 1.2] (1.4,0) -- (1.7,-0.2) -- (1.7,0.2) -- (1.4,0);
    \draw[fill=white,line width = 1.2] (3.0,0) -- (3.3,-0.2) -- (3.3,0.2) -- (3.0,0);
    \node at (-0.2,-0.5) {$P^{j,(1)}$};
    \node at (0.8,-0.5) {$\tilde{S}^{j,(1)}$};
    \node at (1.8,-0.5) {$Q^{j,(2)}$};
    \node at (3.2,-0.5) {$Q^{j,(3)}$};
    \node at (0.4,0.2) {$\alpha_1$};
    \node at (1.2,0.2) {$\alpha_1'$};
    \node at (2.4,0.2) {$\alpha_2$};
    \node at (0,1) {$k_1$};
    \node at (1.6,1) {$k_2$};
    \node at (3.2,1) {$k_3$};
    \node at (1.6,-1) {$\tilde{\bf}_\rom{2}^j$};
    \end{tikzpicture}
    = 
    \begin{tikzpicture}[baseline=0]
    \draw (0,0) -- (1.6,0);
    \draw (1.6,0) -- (3.2,0);
    \draw (0,0) -- (0,0.8);
    \draw (1.6,0) -- (1.6,0.8);
    \draw (3.2,0) -- (3.2,0.8);
    \draw[fill=white,line width = 1.2] (0.2,0) -- (-0.1,-0.2) -- (-0.1,0.2) -- (0.2,0);
    \draw[fill=white,line width = 1.2] (1.6,0) circle (0.2);
    \draw[fill=white,line width = 1.2] (3.0,0) -- (3.3,-0.2) -- (3.3,0.2) -- (3.0,0);
    \node at (0,-0.5) {$P^{j,(1)}$};
    \node at (1.6,-0.5) {$C^{j,(2)}$};
    \node at (3.2,-0.5) {$Q^{j,(3)}$};
    \node at (0.8,0.2) {$\alpha_1$};
    \node at (2.4,0.2) {$\alpha_2$};
    \node at (0,1) {$k_1$};
    \node at (1.6,1) {$k_2$};
    \node at (3.2,1) {$k_3$};
    \node at (1.6,-1) {$\bf_\rom{3}^j$};
    \end{tikzpicture}
    \end{equation}
    
    \item \textbf{Step 3a:}
    Compute $\bG^{j,(2)}:= \mathcal{G}_{\Delta t}^{j}(\bf_\rom{3}^{1},\cdots,\bf_\rom{3}^{N_x})$ using the scheme \eqref{eq_distribution_propagation}. 
    \item \textbf{Step 3b:} Compute the projection  
    $\tilde{C}^{j,(2)}
        := \left\langle
            P^{j,(1)},
            \bG^{j,(2)},
            Q^{j,(3)}
        \right\rangle_{\rom{3}}$. Componentwise, this is
    \begin{equation}
        \tilde{C}_{\alpha_1k_2\alpha_2}^{j,(2)}
        = \sum_{k_1=1}^{N_v} \sum_{k_3=1}^{N_v}
        P_{k_1\alpha_1}^{j,(1)} \bG_{k_1k_2k_3}^{j,(2)} Q_{\alpha_2k_3}^{j,(3)},
    \end{equation}
    or diagrammatically,
    \begin{equation}
    \begin{tikzpicture}[baseline=0]
        \draw (0,0) -- (0,0.7);
        \draw (-0.7,0) -- (0.7,0);
        \draw[fill=white,line width = 1.2] (0,0) circle (0.2);
        \node at (0,-0.5) {$\tilde{C}^{j,(2)}$};
        \node at (0,0.9) {$k_2$};
        \node at (-0.6,0.2) {$\alpha_1$};
        \node at (0.6,0.2) {$\alpha_2$};
    \end{tikzpicture}
    =
    \begin{tikzpicture}[baseline=0]
        \draw (0,0) -- (0, 1.4);
        \draw (1.6,0) -- (1.6, 1.4);
        \draw (3.2,0) -- (3.2, 1.4);
        \draw (0,0) -- (3.2,0);
        \draw (0,1.4) -- (3.2,1.4);
        \draw[rounded corners=2pt, dashed, gray] (-0.4,1) rectangle (3.6,1.8);
        \draw[fill=white,line width = 1.2] (0,1.4) circle (0.2);
        \draw[fill=white,line width = 1.2] (1.6,1.4) circle (0.2);
        \draw[fill=white,line width = 1.2] (3.2,1.4) circle (0.2);
        \draw[fill=white,line width = 1.2] (0.2,0) -- (-0.1,-0.2) -- (-0.1,0.2) -- (0.2,0);
        \draw[fill=white,line width = 1.2] (3.0,0) -- (3.3,-0.2) -- (3.3,0.2) -- (3.0,0);
        \draw[fill=white,line width = 1.2, draw=white] (1.1,-0.5) -- (1.1,0.5) -- (2.1,0.5) -- (2.1,-0.5) -- (1.1,-0.5);
        \node at (1.6,2.1) {$\bG^{j,(2)}$};
        \node at (0,-0.5) {$P^{j,(1)}$};
        \node at (3.2,-0.5) {$Q^{j,(3)}$};
        \node at (0.2,0.7) {$k_1$};
        \node at (1.8,0.7) {$k_2$};
        \node at (3.4,0.7) {$k_3$};
        \node at (1,0.2) {$\alpha_1$};
        \node at (2.4,0.2) {$\alpha_2$};
    \end{tikzpicture}
    \end{equation}
    \item \textbf{Step 3c:} Replace $C^{j,(2)}$ in $\bf_\rom{3}^j$ by $\tilde{C}^{j,(2)}$ to obtain $\tilde{\bf}_\rom{3}^j$:
    \begin{equation}
        \begin{tikzpicture}[baseline=0]
        \draw (0,0) -- (1.6,0);
        \draw (1.6,0) -- (3.2,0);
        \draw (0,0) -- (0,0.8);
        \draw (1.6,0) -- (1.6,0.8);
        \draw (3.2,0) -- (3.2,0.8);
        \draw[fill=white,line width = 1.2] (0.2,0) -- (-0.1,-0.2) -- (-0.1,0.2) -- (0.2,0);
        \draw[fill=white,line width = 1.2] (1.6,0) circle (0.2);
        \draw[fill=white,line width = 1.2] (3.0,0) -- (3.3,-0.2) -- (3.3,0.2) -- (3.0,0);
        \node at (0,-0.5) {$P^{j,(1)}$};
        \node at (1.6,-0.5) {$C^{j,(2)}$};
        \node at (3.2,-0.5) {$Q^{j,(3)}$};
        \node at (0.8,0.2) {$\alpha_1$};
        \node at (2.4,0.2) {$\alpha_2$};
        \node at (0,1) {$k_1$};
        \node at (1.6,1) {$k_2$};
        \node at (3.2,1) {$k_3$};
        \node at (1.6,-1) {$\bf_\rom{3}^j$};
        \end{tikzpicture} 
        \rightarrow
        \begin{tikzpicture}[baseline=0]
    \draw (0,0) -- (1.6,0);
    \draw (1.6,0) -- (3.2,0);
    \draw (0,0) -- (0,0.8);
    \draw (1.6,0) -- (1.6,0.8);
    \draw (3.2,0) -- (3.2,0.8);
    \draw[fill=white,line width = 1.2] (0.2,0) -- (-0.1,-0.2) -- (-0.1,0.2) -- (0.2,0);
    \draw[fill=white,line width = 1.2] (1.6,0) circle (0.2);
    \draw[fill=white,line width = 1.2] (3.0,0) -- (3.3,-0.2) -- (3.3,0.2) -- (3.0,0);
    \node at (0,-0.5) {$P^{j,(1)}$};
    \node at (1.6,-0.5) {\textcolor{red}{$\tilde{C}^{j,(2)}$}};
    \node at (3.2,-0.5) {$Q^{j,(3)}$};
    \node at (0.8,0.2) {$\alpha_1$};
    \node at (2.4,0.2) {$\alpha_2$};
    \node at (0,1) {$k_1$};
    \node at (1.6,1) {$k_2$};
    \node at (3.2,1) {$k_3$};
    \node at (1.6,-1) {$\tilde{\bf}_\rom{3}^j$};
    \end{tikzpicture}
    \end{equation}
    \item \textbf{Step 3d:} Write $\tilde{\bf}_\rom{3}^j$ in form \rombracket{4} by performing the QR decomposition on $\tilde{C}^{j,(2)}$, and denote the resulting tensor train by $\bf_\rom{4}^j$:
    \begin{equation}
    \begin{tikzpicture}[baseline=0]
    \draw (0,0) -- (1.6,0);
    \draw (1.6,0) -- (3.2,0);
    \draw (0,0) -- (0,0.8);
    \draw (1.6,0) -- (1.6,0.8);
    \draw (3.2,0) -- (3.2,0.8);
    \draw[fill=white,line width = 1.2] (0.2,0) -- (-0.1,-0.2) -- (-0.1,0.2) -- (0.2,0);
    \draw[fill=white,line width = 1.2] (1.6,0) circle (0.2);
    \draw[fill=white,line width = 1.2] (3.0,0) -- (3.3,-0.2) -- (3.3,0.2) -- (3.0,0);
    \node at (0,-0.5) {$P^{j,(1)}$};
    \node at (1.6,-0.5) {$\tilde{C}^{j,(2)}$};
    \node at (3.2,-0.5) {$Q^{j,(3)}$};
    \node at (0.8,0.2) {$\alpha_1$};
    \node at (2.4,0.2) {$\alpha_2$};
    \node at (0,1) {$k_1$};
    \node at (1.6,1) {$k_2$};
    \node at (3.2,1) {$k_3$};
    \node at (1.6,-1) {$\tilde{\bf}_\rom{3}^j$};
    \end{tikzpicture}
    =
    \begin{tikzpicture}[baseline=0]
    \draw (0,0) -- (1.6,0);
    \draw (1.6,0) -- (3.2,0);
    \draw (0,0) -- (0,0.8);
    \draw (1.6,0) -- (1.6,0.8);
    \draw (3.2,0) -- (3.2,0.8);
    \draw[fill=white,line width = 1.2] (0.2,0) -- (-0.1,-0.2) -- (-0.1,0.2) -- (0.2,0);
    \draw[fill=white,line width = 1.2] (1.8,0) -- (1.5,-0.2) -- (1.5,0.2) -- (1.8,0);
    \draw[fill=white,line width = 1.2] (2.4,0) circle (0.2);
    \draw[fill=white,line width = 1.2] (3.0,0) -- (3.3,-0.2) -- (3.3,0.2) -- (3.0,0);
    \node at (0,-0.5) {$P^{j,(1)}$};
    \node at (1.4,-0.5) {$P^{j,(2)}$};
    \node at (2.4,-0.5) {$S^{j,(2)}$};
    \node at (3.4,-0.5) {$Q^{j,(3)}$};
    \node at (0.8,0.2) {$\alpha_1$};
    \node at (2.0,0.2) {$\alpha_2'$};
    \node at (2.8,0.2) {$\alpha_2$};
    \node at (0,1) {$k_1$};
    \node at (1.6,1) {$k_2$};
    \node at (3.2,1) {$k_3$};
    \node at (1.6,-1) {$\bf_\rom{4}^j$};
    \end{tikzpicture}
    \end{equation}
    \item \textbf{Step 4a:} 
    Compute $\bH^{j,(2)}: = \mathcal{H}_{\Delta t}^j(\tilde{\bf}_\rom{3}^1,\cdots,\tilde{\bf}_\rom{3}^{N_x})$ using the scheme \eqref{eq_distribution_backward_propagation}.
    \item \textbf{Step 4b:} Compute the projection $\tilde{S}^{j,(2)} := \left\langle
            P^{j,(1)},
            P^{j,(2)},
            \bH^{j,(2)},
            Q^{j,(3)}
        \right\rangle_\rom{4}$. Componentwise, this is
    \begin{equation}
        \tilde{S}^{j,(2)}_{\alpha_2'\alpha_2}
        = \sum_{k_1=1}^{N_v}\sum_{k_2=1}^{N_v} \sum_{k_3=1}^{N_v} \sum_{\alpha_1=1}^{r_1} 
        P^{j,(1)}_{k_1\alpha_1} P^{j,(2)}_{\alpha_1 k_2 \alpha_2'}
        \bH_{k_1k_2k_3}^{j,(2)} Q_{\alpha_2k_3}^{j,(3)},
    \end{equation}
    or diagrammatically,
    \begin{equation}
    \begin{tikzpicture}[baseline=0]
        \draw (-0.7,0) -- (0.7,0);
        \draw[fill=white,line width = 1.2] (0,0) circle (0.2);
        \node at (0,-0.5) {$\tilde{S}^{j,(2)}$};
        \node at (-0.6,0.2) {$\alpha_2'$};
        \node at (0.6,0.2) {$\alpha_2$};
    \end{tikzpicture}
    =
    \begin{tikzpicture}[baseline=0]
        \draw (0,0) -- (0, 1.4);
        \draw (1.6,0) -- (1.6, 1.4);
        \draw (3.2,0) -- (3.2, 1.4);
        \draw (0,0) -- (3.2,0);
        \draw (0,1.4) -- (3.2,1.4);
        \draw[rounded corners=2pt, dashed, gray] (-0.4,1) rectangle (3.6,1.8);
        \draw[fill=white,line width = 1.2] (0,1.4) circle (0.2);
        \draw[fill=white,line width = 1.2] (1.6,1.4) circle (0.2);
        \draw[fill=white,line width = 1.2] (3.2,1.4) circle (0.2);
        \draw[fill=white,line width = 1.2] (0.2,0) -- (-0.1,-0.2) -- (-0.1,0.2) -- (0.2,0);
        \draw[fill=white,line width = 1.2] (1.8,0) -- (1.5,-0.2) -- (1.5,0.2) -- (1.8,0);
        \draw[fill=white,line width = 1.2] (3.0,0) -- (3.3,-0.2) -- (3.3,0.2) -- (3.0,0);
        \draw[fill=white,line width = 1.2, draw=white] (2.2,-0.5) -- (2.2,0.5) -- (2.6,0.5) -- (2.6,-0.5) -- (2.2,-0.5);
        \node at (1.6,2.1) {$\bH^{j,(2)}$};
        \node at (0,-0.5) {$P^{j,(1)}$};
        \node at (1.6,-0.5) {$P^{j,(2)}$};
        \node at (3.2,-0.5) {$Q^{j,(3)}$};
        \node at (0.2,0.7) {$k_1$};
        \node at (1.8,0.7) {$k_2$};
        \node at (3.4,0.7) {$k_3$};
        \node at (0.8,0.2) {$\alpha_1$};
        \node at (2.0,0.2) {$\alpha_2'$};
        \node at (2.8,0.2) {$\alpha_2$};
    \end{tikzpicture}
    \end{equation}
    \item \textbf{Step 4c:} Replace $S^{j,(2)}$ in $\bf_\rom{4}^j$ by $\tilde{S}^{j,(2)}$ to obtain $\tilde{\bf}_\rom{4}^j$:
    \begin{equation}
        \begin{tikzpicture}[baseline=0]
    \draw (0,0) -- (1.6,0);
    \draw (1.6,0) -- (3.2,0);
    \draw (0,0) -- (0,0.8);
    \draw (1.6,0) -- (1.6,0.8);
    \draw (3.2,0) -- (3.2,0.8);
    \draw[fill=white,line width = 1.2] (0.2,0) -- (-0.1,-0.2) -- (-0.1,0.2) -- (0.2,0);
    \draw[fill=white,line width = 1.2] (1.8,0) -- (1.5,-0.2) -- (1.5,0.2) -- (1.8,0);
    \draw[fill=white,line width = 1.2] (2.4,0) circle (0.2);
    \draw[fill=white,line width = 1.2] (3.0,0) -- (3.3,-0.2) -- (3.3,0.2) -- (3.0,0);
    \node at (0,-0.5) {$P^{j,(1)}$};
    \node at (1.4,-0.5) {$P^{j,(2)}$};
    \node at (2.4,-0.5) {$S^{j,(2)}$};
    \node at (3.4,-0.5) {$Q^{j,(3)}$};
    \node at (0.8,0.2) {$\alpha_1$};
    \node at (2.0,0.2) {$\alpha_2'$};
    \node at (2.8,0.2) {$\alpha_2$};
    \node at (0,1) {$k_1$};
    \node at (1.6,1) {$k_2$};
    \node at (3.2,1) {$k_3$};
    \node at (1.6,-1) {$\bf_\rom{4}^j$};
    \end{tikzpicture}
    \rightarrow 
    \begin{tikzpicture}[baseline=0]
    \draw (0,0) -- (1.6,0);
    \draw (1.6,0) -- (3.2,0);
    \draw (0,0) -- (0,0.8);
    \draw (1.6,0) -- (1.6,0.8);
    \draw (3.2,0) -- (3.2,0.8);
    \draw[fill=white,line width = 1.2] (0.2,0) -- (-0.1,-0.2) -- (-0.1,0.2) -- (0.2,0);
    \draw[fill=white,line width = 1.2] (1.8,0) -- (1.5,-0.2) -- (1.5,0.2) -- (1.8,0);
    \draw[fill=white,line width = 1.2] (2.4,0) circle (0.2);
    \draw[fill=white,line width = 1.2] (3.0,0) -- (3.3,-0.2) -- (3.3,0.2) -- (3.0,0);
    \node at (0,-0.5) {$P^{j,(1)}$};
    \node at (1.4,-0.5) {$P^{j,(2)}$};
    \node at (2.4,-0.5) {\textcolor{red}{$\tilde{S}^{j,(2)}$}};
    \node at (3.4,-0.5) {$Q^{j,(3)}$};
    \node at (0.8,0.2) {$\alpha_1$};
    \node at (2.0,0.2) {$\alpha_2'$};
    \node at (2.8,0.2) {$\alpha_2$};
    \node at (0,1) {$k_1$};
    \node at (1.6,1) {$k_2$};
    \node at (3.2,1) {$k_3$};
    \node at (1.6,-1) {$\tilde{\bf}_\rom{4}^j$};
    \end{tikzpicture}.
    \end{equation}
    \item \textbf{Step 4d:} Write $\tilde{\bf}_\rom{4}^j$ in form \rombracket{5} by multiplying $\tilde{S}^{j,(2)}$ with $Q^{j,(3)}$, and denote the resulting tensor train by $\bf_\rom{5}^j$:
    \begin{equation}
        \begin{tikzpicture}[baseline=0]
    \draw (0,0) -- (1.6,0);
    \draw (1.6,0) -- (3.2,0);
    \draw (0,0) -- (0,0.8);
    \draw (1.6,0) -- (1.6,0.8);
    \draw (3.2,0) -- (3.2,0.8);
    \draw[fill=white,line width = 1.2] (0.2,0) -- (-0.1,-0.2) -- (-0.1,0.2) -- (0.2,0);
    \draw[fill=white,line width = 1.2] (1.8,0) -- (1.5,-0.2) -- (1.5,0.2) -- (1.8,0);
    \draw[fill=white,line width = 1.2] (2.4,0) circle (0.2);
    \draw[fill=white,line width = 1.2] (3.0,0) -- (3.3,-0.2) -- (3.3,0.2) -- (3.0,0);
    \node at (0,-0.5) {$P^{j,(1)}$};
    \node at (1.4,-0.5) {$P^{j,(2)}$};
    \node at (2.4,-0.5) {$\tilde{S}^{j,(2)}$};
    \node at (3.4,-0.5) {$Q^{j,(3)}$};
    \node at (0.8,0.2) {$\alpha_1$};
    \node at (2.0,0.2) {$\alpha_2'$};
    \node at (2.8,0.2) {$\alpha_2$};
    \node at (0,1) {$k_1$};
    \node at (1.6,1) {$k_2$};
    \node at (3.2,1) {$k_3$};
    \node at (1.6,-1) {$\tilde{\bf}_\rom{4}^j$};
    \end{tikzpicture} 
    =
    \begin{tikzpicture}[baseline=0]
    \draw (0,0) -- (1.6,0);
    \draw (1.6,0) -- (3.2,0);
    \draw (0,0) -- (0,0.8);
    \draw (1.6,0) -- (1.6,0.8);
    \draw (3.2,0) -- (3.2,0.8);
    \draw[fill=white,line width = 1.2] (0.2,0) -- (-0.1,-0.2) -- (-0.1,0.2) -- (0.2,0);
    \draw[fill=white,line width = 1.2] (1.8,0) -- (1.5,-0.2) -- (1.5,0.2) -- (1.8,0);
    \draw[fill=white,line width = 1.2] (3.2,0) circle (0.2);
    \node at (0,-0.5) {$P^{j,(1)}$};
    \node at (1.6,-0.5) {$P^{j,(2)}$};
    \node at (3.2,-0.5) {$C^{j,(3)}$};
    \node at (0.8,0.2) {$\alpha_1$};
    \node at (2.4,0.2) {$\alpha_2$};
    \node at (0,1) {$k_1$};
    \node at (1.6,1) {$k_2$};
    \node at (3.2,1) {$k_3$};
    \node at (1.6,-1) {$\bf_\rom{5}^{j}$};
    \end{tikzpicture}
    \end{equation}
    \item \textbf{Step 5a:} 
    Compute $\bG^{j,(3)}:= \mathcal{G}_{\Delta t}^{j}(\bf_\rom{5}^{1},\cdots,\bf_\rom{5}^{N_x})$ using the scheme \eqref{eq_distribution_propagation}. 
    \item \textbf{Step 5b:} Compute the projection $\tilde{C}^{j,(3)}
        := \left\langle
            P^{j,(1)},
            P^{j,(2)},
            \bG^{j,(3)}
        \right\rangle_\rom{5}$. Componentwise, this is
    \begin{equation}
        \tilde{C}^{j,(3)}_{\alpha_2k_3}
        = \sum_{k_1=1}^{N_v}\sum_{k_2=1}^{N_v}\sum_{\alpha_1=1}^{r_1} 
    P^{j,(1)}_{k_1\alpha_1} P^{j,(2)}_{\alpha_1 k_2 \alpha_2}
    \bG_{k_1k_2k_3}^{j,(3)},
    \end{equation}
    or diagrammatically,
    \begin{equation}
    \begin{tikzpicture}[baseline=0]
        \draw (-0.7,0) -- (0,0) -- (0,0.7);
        \draw[fill=white,line width = 1.2] (0,0) circle (0.2);
        \node at (0,-0.5) {$\tilde{C}^{j,(3)}$};
        \node at (-0.6,0.2) {$\alpha_2$};
        \node at (0,1) {$k_3$};
    \end{tikzpicture}
    =
    \begin{tikzpicture}[baseline=0]
        \draw (0,0) -- (0, 1.4);
        \draw (1.6,0) -- (1.6, 1.4);
        \draw (3.2,0) -- (3.2, 1.4);
        \draw (0,0) -- (3.2,0);
        \draw (0,1.4) -- (3.2,1.4);
        \draw[rounded corners=2pt, dashed, gray] (-0.4,1) rectangle (3.6,1.8);
        \draw[fill=white,line width = 1.2] (0,1.4) circle (0.2);
        \draw[fill=white,line width = 1.2] (1.6,1.4) circle (0.2);
        \draw[fill=white,line width = 1.2] (3.2,1.4) circle (0.2);
        \draw[fill=white,line width = 1.2] (0.2,0) -- (-0.1,-0.2) -- (-0.1,0.2) -- (0.2,0);
        \draw[fill=white,line width = 1.2] (1.8,0) -- (1.5,-0.2) -- (1.5,0.2) -- (1.8,0);
        \draw[fill=white,line width = 1.2] (3.0,0) -- (3.3,-0.2) -- (3.3,0.2) -- (3.0,0);
        \draw[fill=white,line width = 1.2, draw=white] (2.7,-0.5) -- (2.7,0.5) -- (3.7,0.5) -- (3.7,-0.5) -- (2.7,-0.5);
        \node at (1.6,2.1) {$\bG^{j,(3)}$};
        \node at (0,-0.5) {$P^{j,(1)}$};
        \node at (1.6,-0.5) {$P^{j,(2)}$};
        \node at (0.2,0.7) {$k_1$};
        \node at (1.8,0.7) {$k_2$};
        \node at (3.4,0.7) {$k_3$};
        \node at (0.8,0.2) {$\alpha_1$};
        \node at (2.4,0.2) {$\alpha_2$};
    \end{tikzpicture}
    \end{equation}
    \item \textbf{Step 5c:} Replace $C^{j,(3)}$ in $\bf_\rom{5}^j$ by $\tilde{C}^{j,(3)}$ to obtain $\tilde{\bf}_{\rom{5}}^j$:
    \begin{equation}
    \begin{aligned}
        \begin{tikzpicture}[baseline=0]
        \draw (0,0) -- (1.6,0);
        \draw (1.6,0) -- (3.2,0);
        \draw (0,0) -- (0,0.8);
        \draw (1.6,0) -- (1.6,0.8);
        \draw (3.2,0) -- (3.2,0.8);
        \draw[fill=white,line width = 1.2] (0.2,0) -- (-0.1,-0.2) -- (-0.1,0.2) -- (0.2,0);
        \draw[fill=white,line width = 1.2] (1.8,0) -- (1.5,-0.2) -- (1.5,0.2) -- (1.8,0);
        \draw[fill=white,line width = 1.2] (3.2,0) circle (0.2);
        \node at (0,-0.5) {$P^{j,(1)}$};
        \node at (1.6,-0.5) {$P^{j,(2)}$};
        \node at (3.2,-0.5) {$C^{j,(3)}$};
        \node at (0.8,0.2) {$\alpha_1$};
        \node at (2.4,0.2) {$\alpha_2$};
        \node at (0,1) {$k_1$};
        \node at (1.6,1) {$k_2$};
        \node at (3.2,1) {$k_3$};
        \node at (1.6,-1) {${\bf}_\rom{5}^j$};
        \end{tikzpicture}
        \rightarrow
        \begin{tikzpicture}[baseline=0]
        \draw (0,0) -- (1.6,0);
        \draw (1.6,0) -- (3.2,0);
        \draw (0,0) -- (0,0.8);
        \draw (1.6,0) -- (1.6,0.8);
        \draw (3.2,0) -- (3.2,0.8);
        \draw[fill=white,line width = 1.2] (0.2,0) -- (-0.1,-0.2) -- (-0.1,0.2) -- (0.2,0);
        \draw[fill=white,line width = 1.2] (1.8,0) -- (1.5,-0.2) -- (1.5,0.2) -- (1.8,0);
        \draw[fill=white,line width = 1.2] (3.2,0) circle (0.2);
        \node at (0,-0.5) {$P^{j,(1)}$};
        \node at (1.6,-0.5) {$P^{j,(2)}$};
        \node at (3.2,-0.5) {\textcolor{red}{$\tilde{C}^{j,(3)}$}};
        \node at (0.8,0.2) {$\alpha_1$};
        \node at (2.4,0.2) {$\alpha_2$};
        \node at (0,1) {$k_1$};
        \node at (1.6,1) {$k_2$};
        \node at (3.2,1) {$k_3$};
        \node at (1.6,-1) {$\tilde{\bf}_\rom{5}^j$};
        \end{tikzpicture}
        \end{aligned}
    \end{equation}
    \item \textbf{Step 5d:} Write $\tilde{\bf}_\rom{5}^j$ in form \rombracket{1} by reversing the  \textbf{operation 8} in Section~\ref{sec_TT_operators}, and denote the resulting tensor train by $\bf_\rom{1}^{j,n+1}$:
    \begin{equation}
        \begin{tikzpicture}[baseline=0]
        \draw (0,0) -- (1.6,0);
        \draw (1.6,0) -- (3.2,0);
        \draw (0,0) -- (0,0.8);
        \draw (1.6,0) -- (1.6,0.8);
        \draw (3.2,0) -- (3.2,0.8);
        \draw[fill=white,line width = 1.2] (0.2,0) -- (-0.1,-0.2) -- (-0.1,0.2) -- (0.2,0);
        \draw[fill=white,line width = 1.2] (1.8,0) -- (1.5,-0.2) -- (1.5,0.2) -- (1.8,0);
        \draw[fill=white,line width = 1.2] (3.2,0) circle (0.2);
        \node at (0,-0.5) {$P^{j,(1)}$};
        \node at (1.6,-0.5) {$P^{j,(2)}$};
        \node at (3.2,-0.5) {$\tilde{C}^{j,(3)}$};
        \node at (0.8,0.2) {$\alpha_1$};
        \node at (2.4,0.2) {$\alpha_2$};
        \node at (0,1) {$k_1$};
        \node at (1.6,1) {$k_2$};
        \node at (3.2,1) {$k_3$};
        \node at (1.6,-1) {$\tilde{\bf}_\rom{5}^j$};
        \end{tikzpicture}
        = \begin{tikzpicture}[baseline=0]
        \draw (0,0) -- (1.6,0);
        \draw (1.6,0) -- (3.2,0);
        \draw (0,0) -- (0,0.8);
        \draw (1.6,0) -- (1.6,0.8);
        \draw (3.2,0) -- (3.2,0.8);
        \draw[fill=white,line width = 1.2] (0,0) circle (0.2);
        \draw[fill=white,line width = 1.2] (1.4,0) -- (1.7,-0.2) -- (1.7,0.2) -- (1.4,0);
        \draw[fill=white,line width = 1.2] (3.0,0) -- (3.3,-0.2) -- (3.3,0.2) -- (3.0,0);
        \node at (0,-0.5) {$C^{j,(1)}$};
        \node at (1.6,-0.5) {$Q^{j,(2)}$};
        \node at (3.2,-0.5) {$Q^{j,(3)}$};
        \node at (0.8,0.2) {$\alpha_1$};
        \node at (2.4,0.2) {$\alpha_2$};
        \node at (0,1) {$k_1$};
        \node at (1.6,1) {$k_2$};
        \node at (3.2,1) {$k_3$};
        \node at (1.6,-1) {$\bf_\rom{1}^{j,n+1}$};
        \end{tikzpicture}
    \end{equation}
\end{itemize}
This finishes the algorithm.

\subsection{Memory requirement and computational complexity}

In this subsection, we compare the proposed low-rank algorithm with the full tensor method in terms of memory requirement and computational complexity.

In terms of memory, the full tensor method requires $\mathcal{O}(N_xN_v^3)$ storage, whereas the low-rank algorithm requires only $\mathcal{O}(r^2N_xN_v)$, assuming $r=\max\{r_1,r_2\}$ is the maximal TT-rank.

In terms of computational complexity, the full tensor method requires at least $\mathcal{O}(N_x N_v^3)$ operations per time step. For the low-rank algorithm, although it consists of five steps, each step shares a similar structure. Specifically, \textbf{steps a and b} compute a tensor train using either \eqref{eq_distribution_propagation} or \eqref{eq_distribution_backward_propagation} and then compute the projection of the resulting tensor train. As will be elaborated below, these two steps together typically cost $\mathcal{O}(rR^2 N_v)$, where $r\leq R\ll N_v$. \textbf{Steps c and d} involve updating a tensor core, whose cost is negligible, followed by either a tensor contraction or a QR decomposition, which costs at most $\mathcal{O}(r^3N_v)$. Therefore, accounting for all spatial grid points, the total computational complexity of the proposed low-rank algorithm is $\mathcal{O}(rR^2N_xN_v)$ per time step.

To better illustrate the
computational cost in \textbf{steps a and b} above, we use \textbf{steps 2ab} as an example (the other steps share a similar cost). First, in practice, the evaluation of $\mathbf{H}^{j,(1)}$ is almost never implemented in a lump sum, but instead is split into a finite sum of tensors, namely, $\bH^{j,(1)} = \mathbf{h}_1^{j,(1)} + \mathbf{h}_2^{j,(1)} + \cdots + \mathbf{h}_P^{j,(1)}$, where we assume that the maximal rank of each tensor $\mathbf{h}_p^{j,(1)}$ is $R$. This can often be identified straightforwardly, depending the structure of the equation and the discretization. Typically, $P=\mathcal{O}(1)$, and $R\geq r$ but is still of the same order as $r$. We then perform the projection for each tensor $\mathbf{h}_p^{j,(1)}$. From \eqref{eq_step2b}, it appears that this step would require $\mathcal{O}(r^3N_v^3)$ operations. However, using the fact that $\mathbf{h}_p^{j,(1)}$ is represented in TT format:
\begin{equation}
    \mathbf{h}^{j,(1)}_{p,k_1k_2k_3} = \sum_{\gamma_1=1}^{R} \sum_{\gamma_2=1}^{R} H_{k_1\gamma_1}^{j,(1)} H_{\gamma_1 k_2\gamma_2}^{j,(2)} H_{\gamma_2 k_3}^{j,(3)},
\end{equation}
the cost can be significantly reduced.
Indeed, \eqref{eq_step2b} can then be rewritten as
\begin{equation}
\begin{aligned}
\tilde{S}_{\alpha_1\alpha_1'}^{j,(1)}
        &=\sum_{k_1=1}^{N_v}\sum_{k_2=1}^{N_v} \sum_{k_3=1}^{N_v} \sum_{\alpha_2=1}^{r_2} 
        P^{j,(1)}_{k_1\alpha_1}
        \sum_{\gamma_1=1}^{R} \sum_{\gamma_2=1}^{R} H_{k_1\gamma_1}^{j,(1)} H_{\gamma_1 k_2\gamma_2}^{j,(2)} H_{\gamma_2 k_3}^{j,(3)} Q_{\alpha_1'k_2\alpha_2}^{j,(2)} Q_{\alpha_2k_3}^{j,(3)} \\
        &= \sum_{k_1=1}^{N_v} P_{k_1,\alpha_1}^{j,(1)}
        \left(\sum_{\gamma_1=1}^{R} H_{k_1\gamma_1}^{j,(1)}
        \left(\sum_{\alpha_2=1}^{r_2} \sum_{k_2=1}^{N_v} Q_{\alpha_1'k_2\alpha_2}^{j,(2)}
        \left(\sum_{\gamma_2=1}^{R} H_{\gamma_1k_2\gamma_2}^{j,(2)}
        \left(\sum_{k_3=1}^{N_v} H_{\gamma_2k_3}^{j,(3)} Q_{\alpha_2k_3}^{j,(3)}\right)\right)\right)\right),
\end{aligned}
\end{equation}
where the expression is evaluated from the innermost parentheses outward.
Diagrammatically, this process can be expressed as 
\begin{equation}
\begin{aligned}
    &\begin{tikzpicture}[baseline=15]
        \draw (0,0) -- (0, 1.4);
        \draw (1.6,0) -- (1.6, 1.4);
        \draw (3.2,0) -- (3.2, 1.4);
        \draw (0,0) -- (3.2,0);
        \draw (0,1.4) -- (3.2,1.4);
        \draw[fill=white,line width = 1.2] (0,1.4) circle (0.2);
        \draw[fill=white,line width = 1.2] (1.6,1.4) circle (0.2);
        \draw[fill=white,line width = 1.2] (3.2,1.4) circle (0.2);
        \draw[fill=white,line width = 1.2] (0.2,0) -- (-0.1,-0.2) -- (-0.1,0.2) -- (0.2,0);
        \draw[fill=white,line width = 1.2] (1.4,0) -- (1.7,-0.2) -- (1.7,0.2) -- (1.4,0);
        \draw[fill=white,line width = 1.2] (3.0,0) -- (3.3,-0.2) -- (3.3,0.2) -- (3.0,0);
        \draw[fill=white,line width = 1.2, draw=white] (0.6,-0.5) -- (0.6,0.5) -- (1.0,0.5) -- (1.0,-0.5) -- (0.6,-0.5);
        \node at (0,1.9) {$H^{j,(1)}$};
        \node at (1.6,1.9) {$H^{j,(2)}$};
        \node at (3.2,1.9) {$H^{j,(3)}$};
        \node at (0,-0.5) {$P^{j,(1)}$};
        \node at (1.6,-0.5) {$Q^{j,(2)}$};
        \node at (3.2,-0.5) {$Q^{j,(3)}$};
        \node at (0.2,0.7) {$k_1$};
        \node at (1.8,0.7) {$k_2$};
        \node at (3.4,0.7) {$k_3$};
        \node at (0.4,0.2) {$\alpha_1$};
        \node at (1.2,0.2) {$\alpha_1'$};
        \node at (2.4,0.2) {$\alpha_2$};
        \node at (0.8,1.6) {$\gamma_1$};
        \node at (2.4,1.6) {$\gamma_2$};
    \end{tikzpicture}
    &=&
    \begin{tikzpicture}[baseline=15]
        \draw (0,0) -- (0, 1.4);
        \draw (1.6,0) -- (1.6, 1.4);
        \draw (3.2,0) -- (3.2, 1.4);
        \draw (0,0) -- (3.2,0);
        \draw (0,1.4) -- (3.2,1.4);
        \draw[fill=white,line width = 1.2] (0,1.4) circle (0.2);
        \draw[fill=white,line width = 1.2] (1.6,1.4) circle (0.2);
        \draw[fill=white,line width = 1.2] (3.2,1.4) circle (0.2);
        \draw[fill=white,line width = 1.2] (0.2,0) -- (-0.1,-0.2) -- (-0.1,0.2) -- (0.2,0);
        \draw[fill=white,line width = 1.2] (1.4,0) -- (1.7,-0.2) -- (1.7,0.2) -- (1.4,0);
        \draw[fill=white,line width = 1.2] (3.0,0) -- (3.3,-0.2) -- (3.3,0.2) -- (3.0,0);
        \draw[fill=white,line width = 1.2, draw=white] (0.6,-0.5) -- (0.6,0.5) -- (1.0,0.5) -- (1.0,-0.5) -- (0.6,-0.5);
        \node at (0,1.9) {$H^{j,(1)}$};
        \node at (1.6,1.9) {$H^{j,(2)}$};
        \node at (0,-0.5) {$P^{j,(1)}$};
        \node at (1.6,-0.5) {$Q^{j,(2)}$};
        \node at (0.2,0.7) {$k_1$};
        \node at (1.8,0.7) {$k_2$};
        \node at (0.4,0.2) {$\alpha_1$};
        \node at (1.2,0.2) {$\alpha_1'$};
        \node at (2.4,0.2) {$\alpha_2$};
        \node at (0.8,1.6) {$\gamma_1$};
        \node at (2.4,1.6) {$\gamma_2$};
        \draw[rounded corners=4pt, fill=white,line width = 1.2] (3.0,-0.2) rectangle (3.4,1.6);
    \end{tikzpicture}
    &=
    \begin{tikzpicture}[baseline=15]
        \draw (0,0) -- (0, 1.4);
        \draw (1.6,0) -- (1.6, 1.4);
        \draw (3.2,0) -- (3.2, 1.4);
        \draw (0,0) -- (3.2,0);
        \draw (0,1.4) -- (3.2,1.4);
        \draw[fill=white,line width = 1.2] (0,1.4) circle (0.2);
        \draw[fill=white,line width = 1.2] (1.6,1.4) circle (0.2);
        \draw[fill=white,line width = 1.2] (3.2,1.4) circle (0.2);
        \draw[fill=white,line width = 1.2] (0.2,0) -- (-0.1,-0.2) -- (-0.1,0.2) -- (0.2,0);
        \draw[fill=white,line width = 1.2] (1.4,0) -- (1.7,-0.2) -- (1.7,0.2) -- (1.4,0);
        \draw[fill=white,line width = 1.2, draw=white] (0.6,-0.5) -- (0.6,0.5) -- (1.0,0.5) -- (1.0,-0.5) -- (0.6,-0.5);
        \node at (0,1.9) {$H^{j,(1)}$};
        \node at (0,-0.5) {$P^{j,(1)}$};
        \node at (1.6,-0.5) {$Q^{j,(2)}$};
        \node at (0.2,0.7) {$k_1$};
        \node at (1.8,0.7) {$k_2$};
        \node at (0.4,0.2) {$\alpha_1$};
        \node at (1.2,0.2) {$\alpha_1'$};
        \node at (2.4,0.2) {$\alpha_2$};
        \node at (0.8,1.6) {$\gamma_1$};
        \draw[rounded corners=4pt, fill=white,line width = 1.2] (2.5,1.2) -- (1.4,1.2) -- (1.4,1.6) -- (3.4,1.6) -- (3.4,-0.2) -- (3.0,-0.2) -- (3.0,1.2) -- (2.0,1.2);
    \end{tikzpicture} \\
    = &
    \begin{tikzpicture}[baseline=15]
        \draw (0,0) -- (0, 1.4);
        \draw (1.6,0) -- (1.6, 1.4);
        \draw (3.2,0) -- (3.2, 1.4);
        \draw (0,0) -- (3.2,0);
        \draw (0,1.4) -- (3.2,1.4);
        \draw[fill=white,line width = 1.2] (0,1.4) circle (0.2);
        \draw[fill=white,line width = 1.2] (1.6,1.4) circle (0.2);
        \draw[fill=white,line width = 1.2] (3.2,1.4) circle (0.2);
        \draw[fill=white,line width = 1.2] (0.2,0) -- (-0.1,-0.2) -- (-0.1,0.2) -- (0.2,0);
        \draw[fill=white,line width = 1.2, draw=white] (0.6,-0.5) -- (0.6,0.5) -- (1.0,0.5) -- (1.0,-0.5) -- (0.6,-0.5);
        \node at (0,1.9) {$H^{j,(1)}$};
        \node at (0,-0.5) {$P^{j,(1)}$};
        \node at (0.2,0.7) {$k_1$};
        \node at (1.8,0.7) {$k_2$};
        \node at (0.4,0.2) {$\alpha_1$};
        \node at (1.2,0.2) {$\alpha_1'$};
        \node at (2.4,0.2) {$\alpha_2$};
        \node at (0.8,1.6) {$\gamma_1$};
        \draw[rounded corners=4pt, fill=white,line width = 1.2] (1.4,-0.2) rectangle (3.4,1.6);
    \end{tikzpicture}
    &= &
    \begin{tikzpicture}[baseline=15]
        \draw (0,0) -- (0, 1.4);
        \draw (1.6,0) -- (1.6, 1.4);
        \draw (3.2,0) -- (3.2, 1.4);
        \draw (0,0) -- (3.2,0);
        \draw (0,1.4) -- (3.2,1.4);
        \draw[fill=white,line width = 1.2] (0,1.4) circle (0.2);
        \draw[fill=white,line width = 1.2] (1.6,1.4) circle (0.2);
        \draw[fill=white,line width = 1.2] (3.2,1.4) circle (0.2);
        \draw[fill=white,line width = 1.2] (0.2,0) -- (-0.1,-0.2) -- (-0.1,0.2) -- (0.2,0);
        \draw[fill=white,line width = 1.2, draw=white] (0.6,-0.5) -- (0.6,0.5) -- (1.0,0.5) -- (1.0,-0.5) -- (0.6,-0.5);
        \node at (0,-0.5) {$P^{j,(1)}$};
        \node at (0.2,0.7) {$k_1$};
        \node at (1.8,0.7) {$k_2$};
        \node at (0.4,0.2) {$\alpha_1$};
        \node at (1.2,0.2) {$\alpha_1'$};
        \node at (2.4,0.2) {$\alpha_2$};
        \draw[rounded corners=4pt, fill=white,line width = 1.2] (1.0,1.2) -- (-0.2,1.2) -- (-0.2,1.6) -- (3.4,1.6) -- (3.4,-0.2) -- (1.4,-0.2) -- (1.4,1.2) -- (0.4,1.2);
    \end{tikzpicture}
    &=\quad 
    \begin{tikzpicture}[baseline=15]
        \draw (0,0) -- (0, 1.4);
        \draw (1.6,0) -- (1.6, 1.4);
        \draw (3.2,0) -- (3.2, 1.4);
        \draw (0,0) -- (3.2,0);
        \draw (0,1.4) -- (3.2,1.4);
        \draw[fill=white,line width = 1.2] (0,1.4) circle (0.2);
        \draw[fill=white,line width = 1.2] (1.6,1.4) circle (0.2);
        \draw[fill=white,line width = 1.2] (3.2,1.4) circle (0.2);
        \draw[fill=white,line width = 1.2, draw=white] (0.6,-0.5) -- (0.6,0.5) -- (1.0,0.5) -- (1.0,-0.5) -- (0.6,-0.5);
        \node at (0.45,0.2) {$\alpha_1$};
        \node at (1.2,0.2) {$\alpha_1'$};
        \draw[rounded corners=4pt, fill=white,line width = 1.2] (1.0,1.2) -- (0.2,1.2) -- (0.2,-0.2) -- (-0.2,-0.2) -- (-0.2,1.6) -- (3.4,1.6) -- (3.4,-0.2) -- (1.4,-0.2) -- (1.4,1.2) -- (0.6,1.2);
    \end{tikzpicture}
\end{aligned}
\end{equation}
The computational costs for each contraction are $\mathcal{O}(rRN_v)$, $\mathcal{O}(rR^2N_v)$, $\mathcal{O}(r^2RN_v)$, $\mathcal{O}({rRN_v})$ and $\mathcal{O}(r^2N_v)$, respectively. Together, this yields the complexity $\mathcal{O}(rR^2N_v)$ as claimed above.

\subsection{Second-order extension}

The above algorithm is first-order accurate in time because it is based on a Lie-Trotter splitting of the projection operator. The method can be extended to second order using Strang splitting \cite{lubich2015time}. We therefore need the following three schemes: 
a half-step forward propagator $\mathcal{G}^j_{\dt/2}$,
a half-step backward propagator $\mathcal{H}^j_{\dt/2}$,
and a full-step forward propagator $\mathcal{G}^j_{\dt}$, where the mappings are defined in \eqref{eq_distribution_propagation}-\eqref{eq_distribution_backward_propagation}, but the underlying time-stepping scheme must be second order. 

Assume that we have $\{\bf^{j,n}_\rom{1}\}_{j=1}^{N_x}$ in form \rombracket{1} at time $t_n$. The following steps compute $\{\bf^{j,n+1}_\rom{1}\}_{j=1}^{N_x}$ in form \rombracket{1} at $t_{n+1}$, which essentially consists of a forward sweep, similar to the first-order scheme, followed by 
a backward sweep, performing the operations in the reverse order. 
\begin{itemize}
    \item Perform \textbf{steps 1, 2, 3, and 4} with half-step propagators $\mathcal{G}_{\dt/2}^j$ and $\mathcal{H}_{\dt/2}^j$.
    \item Perform \textbf{steps 5abc} with the full-step propagator $\mathcal{G}_{\dt}^j$, resulting in $\tilde{\bf}_\rom{5}^j$. 
    Write $\tilde{\bf}_\rom{5}^j$ in form \rombracket{4} by performing the QR decomposition on $\tilde{C}^{j,(3)}$, and denote the resulting tensor by ${\bf}_\rom{4}^j$.
    \item Use $\tilde{\bf}_\rom{5}^j$ as input, perform \textbf{steps 4abc} with the half-step propagator $\mathcal{H}_{\Delta t/2}^j$, resulting in $\tilde{\bf}_\rom{4}^j$. 
    Write $\tilde{\bf}_\rom{4}^j$ in form $\rombracket{3}$ by multiplying $P^{j,(2)}$ and $\tilde{S}^{j,(2)}$, and denote the resulting tensor by $\bf_\rom{3}^j$.
    \item Use $\bf_\rom{3}^j$ as input, perform \textbf{steps 3abc} with the half-step propagator $\mathcal{G}_{\dt/2}^j$, resulting in $\tilde{\bf}_\rom{3}^j$. Write $\tilde{\bf}_\rom{3}^j$ in form \rombracket{2} by performing the QR decomposition on $\tilde{C}^{j,(2)}$, and denote the resulting tensor by $\bf_\rom{2}^j$.
    \item Use $\tilde{\bf}_\rom{3}^j$ as input, perform \textbf{steps 2abc} with the half-step propagator $\mathcal{H}_{\Delta t/2}^j$, resulting in $\tilde{\bf}_\rom{2}^j$. Write $\tilde{\bf}_\rom{2}^j$ in form $\rombracket{1}$ by multiplying $P^{j,(1)}$ with $\tilde{S}^{j,(1)}$, and denote the resulting tensor by $\bf_\rom{1}^j$.
    \item Use $\bf_\rom{1}^j$ as input, perform \textbf{steps 1abc} with the half-step propagator $\mathcal{G}_{\dt/2}^j$, resulting in $\bf_\rom{1}^{j,n+1}$.
\end{itemize}

\section{Application to kinetic equations}
\label{sec_examples}
In this section, we apply the algorithm presented in the previous section to several kinetic equations and demonstrate its performance in both accuracy and efficiency. We begin with the spatially homogeneous case, namely,  equation \eqref{eq_original_equation} without spatial dependence on $x$. We then consider the spatially inhomogeneous case, which includes an extensive discussion of the Vlasov-Amp\`ere-Fokker-Planck equation, a kinetic model widely used to describe plasma dynamics.

\subsection{Spatially homogeneous BGK equation}
\label{sec_spatially_homogeneous_BGK}

We first consider the BGK equation, a simple relaxation-type kinetic model introduced to mimic the full Boltzmann equation \cite{cercignani}. Without spatial dependence, the equation reads:
\begin{equation} \label{eq:BGK}
    \partial_t f(t,\bv) = \eta( \mathcal{M}[f](t,\bv) - f(t,\bv)),
\end{equation}
where $\eta$ is the collision strength and $\mathcal{M}[f]$ is the Maxwellian equilibrium given by
\begin{equation} \label{eq:Max}
    \mathcal{M}[f] = \frac{n}{(2\pi T)^{3/2}} \exp\left(-\frac{\vert \bv - \bu\vert^2}{2T}\right),
\end{equation}
with the density $n$, bulk velocity $\bu$, and temperature $T$ defined by the moments of $f$:
\begin{equation} \label{eq:moments}
n=\int_{\mathbb{R}^3}f\rd{\bv}, \quad \bu=\frac{1}{n}\int_{\mathbb{R}^3}f \bv\rd{\bv}, \quad T=\frac{1}{3n}\int_{\mathbb{R}^3} f|\bv-\bu|^2\rd{\bv}.
\end{equation}
Observing that
\begin{equation}
    \mathcal{M}[f] = \frac{n}{(2\pi T)^{3/2}} 
    \exp\left(-\frac{(v^{(1)}-u^{(1)})^2}{T}\right)
    \exp\left(-\frac{(v^{(2)}-u^{(2)})^2}{T}\right)
    \exp\left(-\frac{(v^{(3)}-u^{(3)})^2}{T}\right),
\end{equation}
it is clear that $\mathcal{M}[f]$ can be represented directly by a tensor train in velocity space with TT-rank $(1,1)$.

For equation \eqref{eq:BGK}, it can be easily verified that $n$, $\bu$, and $T$ remain constant in time, and so does the Maxwellian $\mathcal{M}[f]$; hence we simply denote it by $\mathcal{M}(\bv)$. 

Assume the initial condition is 
\begin{equation}
    f_0(\bv) = \frac{n_1}{(2\pi T_1)^{3/2}}
    \exp\left(
        - \frac{\vert \bv - \bu_1 \vert^2}{2T_1}
    \right)
    + \frac{n_2}{(2\pi T_2)^{3/2}}
    \exp\left(
        - \frac{\vert \bv - \bu_2 \vert^2}{2T_2}
    \right),
\end{equation}
which can again be represented by a tensor train by applying \textbf{operation 2} in \Cref{sec_TT_operators}.
Then \eqref{eq:BGK} can be solved analytically, and the solution is
\begin{equation}
\label{eq_spatial_homo_BGK_exact}
    f_{\text{exact}}(t,\bv) = (1-\e^{-\eta t})\mathcal{M}(\bv) + \e^{-\eta t}f_0(\bv),
\end{equation}
where $\mathcal{M}(\bv)$ is defined in \eqref{eq:Max} with
\begin{equation} \label{eq:Maxpara}
n = n_1 + n_2, \quad \bu = \frac{1}{n}(n_1 \bu_1 + n_2 \bu_2), \quad T = \frac{1}{n}
    \left(
        n_1 T_1 + n_2 T_2 + \frac{n_1 \vert \bu_1 - \bu \vert^2 + n_2 \vert \bu_2 - \bu \vert^2}{3}
    \right).
\end{equation}
Using \textbf{operations 1 and 2}, the exact solution can be represented by a tensor train with TT-rank (3,3).

Since the equation is spatially homogeneous, we need only one tensor train to represent the solution. In this numerical example, we test both the first-order and second-order low-rank algorithms introduced in the previous section. In the first-order algorithm, we use the forward Euler method for \eqref{eq:BGK} in both forward and backward sub-projection steps; that is, the mapping in  \eqref{eq_distribution_propagation}-\eqref{eq_distribution_backward_propagation} are given by
\begin{equation}
\mathcal{G}_{\Delta t}(\bf) = \bf +  \Delta t\eta \left(\bM - \bf\right), \quad 
\mathcal{H}_{\Delta t}(\bf) = \bf -  \Delta t\eta \left(\bM - \bf\right),
\end{equation}
where $\bf$ and $\bM$ denote the TT representations of the solution $f(t,\bv)$ and the Maxwellian $\mathcal{M}(\bv)$, respectively. In the second-order algorithm, we use the Heun's method for \eqref{eq:BGK} in both the forward and backward sub-projection steps, and the mappings $\mathcal{G}_{\dt}$ and $\mathcal{H}_{\dt}$ are given by
\begin{equation}
\begin{split}
    &\mathcal{G}_{\dt}(\bf) = \frac{1}{2} \left(\bf+\bf^{\circ}+\dt\eta(\bM-\bf^\circ)\right),
    \text{~with~}
     \bf^\circ = \bf+\dt \eta \left(\bM-\bf\right), \\
    &\mathcal{H}_{\dt}(\bf) = \frac{1}{2} \left(\bf+\bf^{\circ}-\dt\eta(\bM-\bf^\circ)\right),
    \text{~with~}
    \bf^\circ = \bf-\dt \eta \left(\bM-\bf\right).
\end{split}
\end{equation}
Note that the computation of $\mathcal{G}_{\dt}(\bf)$, $ \mathcal{H}_{\dt}(\bf)$ can be completed in TT form given that $\mathbf{f}$ is represented by a tensor train.

In the numerical test, we choose the following parameter set
\begin{equation}
\begin{split}
    &n_1 = \frac{1}{2}, \quad 
    \bu_1 = [-1,2,0]^\T, \quad
    T_1 = 1, \quad
    n_2 = \frac{1}{2}, \quad
    \bu_2 = [3,-3,2], \quad
    T_2 = 1,
\end{split}
\end{equation}
with collision strength $\eta = 1$. Then $n = 1$, $\bu = \left[1,-\frac{1}{2},1\right]^\T$, and $T = \frac{19}{4}$ according to \eqref{eq:Maxpara}.

\begin{figure}[htp!]
    \centering
    \includegraphics[width=0.4\linewidth]{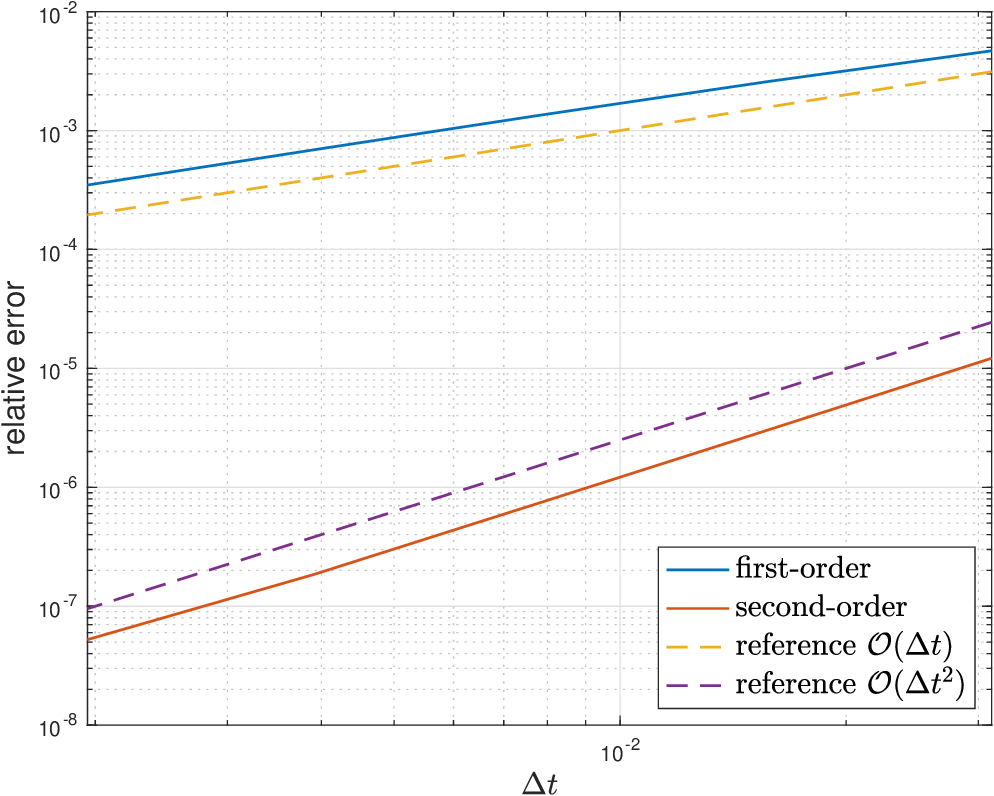}
    \caption{Numerical error for the spatially homogeneous BGK equation.}
    \label{fig_spatial_homo_Boltzmann_BGK}
\end{figure}

We truncate the velocity domain to $[v_{\min},v_{\max}]^3 = [-8,8]^3$ and fix $N_v = 256$ (so $\Delta v=(v_{\max}-v_{\min})/N_v=1/16$). The grid points in each velocity dimension are given by
\begin{equation}
\label{eq_velocity_discretization}
    v_k = v_{\min} + \left(k-\frac{1}{2}\right)\dv, \quad k = 1,\cdots,N_v.
\end{equation}
The TT-rank in the entire simulation is fixed as $(5,5)$. We solve the equation up to time $t=5$ with time steps $\dt=\frac{1}{32},\frac{1}{64},\frac{1}{128}, \frac{1}{256}, \frac{1}{512}$.
Since the exact solution is known, the relative error is defined as $ \displaystyle
    \frac{\Vert \bf_{\mathrm{TT}} - \bf_{\mathrm{exact}} \Vert_F}{\Vert \bf_{\mathrm{exact}} \Vert_F},
$
where $\bf_{\mathrm{TT}}$ is the numerical solution, and the Frobenius norm of a tensor train $\bA$ is defined as $
 \Vert \bA \Vert_F = \left(\sum_{k_1,k_2,k_3=1}^{N_v}\bA_{k_1k_2k_3}^2\right)^{1/2}$. The relative errors for both the first-order and the second-order low-rank algorithms with different time steps $\dt$ are plotted in \Cref{fig_spatial_homo_Boltzmann_BGK}, which clearly demonstrate the expected orders of accuracy.

\subsection{Heat equation}
We next consider the three-dimensional heat equation. Although it is not typically classified as a kinetic equation, this example is useful for assessing the effectiveness of our method on problems  involving diffusive terms, which will be relevant for the Fokker-Planck equation considered below. 
The equation reads
\begin{equation} \label{eq:heat}
    \partial_t f(t,\bv) = \eta \Delta_\bv f(t,\bv),
\end{equation}
with $\eta$ being the diffusion coefficient. Assume the initial condition is
\begin{equation}
    f_0(\bv) 
    = A_1 \exp(-\beta_1 \vert\bv-\bu_1\vert^2)
    + A_2 \exp(-\beta_2 \vert\bv-\bu_2\vert^2).
\end{equation}
Then \eqref{eq:heat} can be solved analytically and the solution is given by
\begin{equation}
    f_{\text{exact}}(t,\bv) 
    = \frac{A_1}{(1+4\eta \beta_1 t)^{3/2}} \exp\left({-\frac{\beta_1 \vert \bv - \bu_1 \vert^2}{1+4\eta\beta_1 t}}\right)
    + \frac{A_2}{(1+4\eta \beta_2 t)^{3/2}} \exp\left({-\frac{\beta_2 \vert \bv - \bu_2 \vert^2}{1+4\eta\beta_2 t}}\right).
\end{equation}

We truncate the velocity domain to $[v_{\min},v_{\max}]^3$ and choose the grid points as in \eqref{eq_velocity_discretization}.
We discretize the diffusion operator $\eta \Delta_{\bv}f$ at the grid point $(v_{k_1},v_{k_2},v_{k_3})$ as
\begin{equation}
\label{eq_discrete_diffusion}
	(\mathbf{L} \bf)_{k_1k_2k_3}
	= \frac{\eta}{\dv} \left(\mathfrak{F}_{k_1+\frac{1}{2}, k_2k_3} - \mathfrak{F}_{k_1-\frac{1}{2},k_2k_3}
 + \mathfrak{F}_{k_1, k_2+\frac{1}{2},k_3} - \mathfrak{F}_{k_1,k_2-\frac{1}{2},k_3} 
 + \mathfrak{F}_{k_1k_2,k_3+\frac{1}{2}} - \mathfrak{F}_{k_1k_2,k_3-\frac{1}{2}}\right),
\end{equation}
where the fluxes $\mathfrak{F}$ at the interior grid points $k_1=1,\cdots,N_v-1$ are given by
\begin{equation}
	\mathfrak{F}_{k_1+\frac{1}{2},k_2k_3} = \frac{\bf_{k_1+1,k_2k_3}-\bf_{k_1k_2k_3}}{\dv},
\end{equation}
and are zero at the boundaries
$
	\mathfrak{F}_{\frac{1}{2},k_2k_3} = \mathfrak{F}_{N_v+\frac{1}{2},k_2k_3} 
	= 0
$.
The definitions in the other two velocity dimensions are analogous.
In the implementation, $\bf_{k_1+1,k_2,k_3}$ and $\bf_{k_1,k_2,k_3}$ are two tensor trains with a single varying core and their difference can be obtained using \textbf{operations 1 and 3} in \Cref{sec_TT_operators}.
The discrete diffusion operator $\mathbf{Lf}$, which involves \textbf{operations 1 and 2}, again has a TT form. 

We then consider the first-order low-rank algorithm with forward Euler method in both forward and backward sub-projection steps, i.e., the mappings in \eqref{eq_distribution_propagation}-\eqref{eq_distribution_backward_propagation} are given by
\begin{equation}
    \mathcal{G}_{\Delta t}(\bf) = \bf + \dt \mathbf{L}\bf,
    \quad
    \mathcal{H}_{\Delta t}(\bf) = \bf - \dt \mathbf{L}\bf.
\end{equation}

\begin{figure}[htp!]
    \centering
    \includegraphics[width=0.4\linewidth]{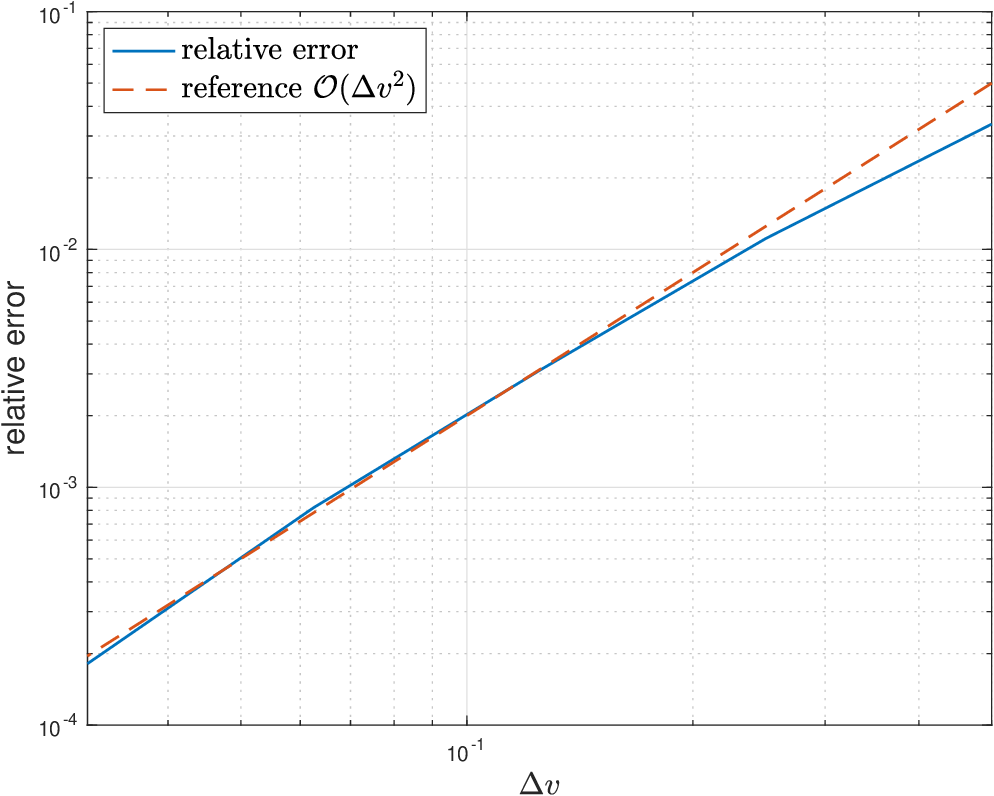}
    \caption{Numerical error for the heat equation.}
    \label{fig_heat_equation}
\end{figure}

In the numerical test, we choose the following parameters
\begin{equation}
\begin{aligned}
    &A_1 = \frac{1}{3}, \quad\bu_1 = [1,2,-1]^T, \quad\beta_1 = 1,\quad
    A_2 = \frac{2}{3}, \quad\bu_2 = [3,-1,-2]^T, \quad\beta_2 = \frac{3}{2},
\end{aligned}
\end{equation}
and the diffusion coefficient $\eta =1$.
We solve the equation up to time $t = 5$ with a fixed TT-rank $(5,5)$.
The velocity domain is chosen as  $[v_{\min},v_{\max}]^3 = [-16,16]^3$ with $\dv = \frac{1}{2}, \frac{1}{4},\frac{1}{8}, \frac{1}{16}, \frac{1}{32}$.
To fulfill the CFL condition of the heat equation, we choose the time step $ \dt=\dfrac{\dv^2}{12\eta}$.
We again compute the relative error 
$\dfrac{
	\Vert \bf_{\mathrm{TT}} - \bf_{\mathrm{exact}} \Vert_F }
	{\Vert \bf_{\mathrm{exact}}\Vert_F}$. As shown in \Cref{fig_heat_equation}, the error exhibits the expected second-order accuracy in velocity, or first-order accuracy in time.

\subsection{Linear Fokker-Planck equation}

The third example is the linear Fokker-Planck equation given by 
\begin{equation} \label{eq:FP}
    \pdt f(t,\bv) = \nabla_\bv \cdot \left(
        {M}(\bv)\nabla_\bv \left(\frac{f}{{M}(\bv)}\right)\right)
    ,
\end{equation}
where the function $M$ is fixed as $\displaystyle
    M(\bv) = \exp\left(-\frac{\vert \bv \vert^2}{2}\right)$. 
Note that the right-hand side of \eqref{eq:FP} can be written equivalently as $\nabla_{\bv}\cdot(\nabla_{\bv}f+\bv f)$, which is the familiar drift-diffusion form.

If the initial condition is chosen as
\begin{equation}
    f_0(\bv) = \frac{1}{(2\pi(1-\e^{-1}))^{3/2}} \exp\left(-\frac{\vert \bv \vert^2}{2-2\e^{-1}}\right),
\end{equation}
then \eqref{eq:FP} admits an analytical solution 
\begin{equation}
    f_{\text{exact}}(t,\bv) = \frac{1}{(2\pi(1-\e^{-2t-1}))^{3/2}}
    \exp\left(-\frac{\vert \bv \vert^2}{2-2\e^{-2t-1}}\right).
\end{equation}

We choose the same grid points as in the previous heat equation and discretize the Fokker-Planck operator at $(v_{k_1},v_{k_2},v_{k_3})$ as
\begin{equation}
\label{eq_Fokker_Planck_discretization}
\begin{aligned}
    (\mathbf{Q}\bf)_{k_1k_2k_3}
    =& \frac{1}{\dv} \left(
        \mathfrak{F}_{k_1+\frac{1}{2},k_2k_3} 
        -  \mathfrak{F}_{k_1-\frac{1}{2},k_2k_3}
        +  \mathfrak{F}_{k_1,k_2+\frac{1}{2},k_3} 
        -  \mathfrak{F}_{k_1,k_2-\frac{1}{2},k_3}
        +  \mathfrak{F}_{k_1k_2,k_3+\frac{1}{2}} 
        -  \mathfrak{F}_{k_1k_2,k_3-\frac{1}{2}}
    \right),
\end{aligned}
\end{equation}
where the fluxes at the interior grid points $k_1=1,\cdots,N_v-1$ are given by
\begin{equation}
\label{eq_Fokker_Planck_fluxes}
\begin{aligned}
    \mathfrak{F}_{k_1+\frac{1}{2},k_2k_3}
=\frac{M_{k_1k_2k_3}+M_{k_1+1,k_2k_3}}{2\dv}
    \left(
        \frac{f_{k_1+1,k_2k_3}}{M_{k_1+1,k_2k_3}}
        - \frac{f_{k_1k_2k_3}}{M_{k_1k_2k_3}}
    \right),
\end{aligned}
\end{equation}
and are zero at the boundaries $\mathfrak{F}_{\frac{1}{2},k_2k_3}=\mathfrak{F}_{N_v+\frac{1}{2},k_2k_3}=0$. The definition in the other two velocity dimensions are analogous. 
In addition to the operations mentioned in the discussion of \eqref{eq_discrete_diffusion}, the computation of $\mathbf{Qf}$ involves \textbf{operation 4 and 5} in \Cref{sec_TT_operators}.

We then consider the first-order low-rank algorithm with forward Euler method in both forward and backward sub-projection steps, i.e., the mappings in \eqref{eq_distribution_propagation}-\eqref{eq_distribution_backward_propagation} are given by 
\begin{equation}
    \mathcal{G}_{\dt}(\bf) = \bf + \dt \mathbf{Q}\bf,
    \quad
    \mathcal{H}_{\dt}(\bf) = \bf - \dt \mathbf{Q}\bf.
\end{equation}

\begin{figure}[htp!]
    \centering
    \includegraphics[width=0.4\linewidth]{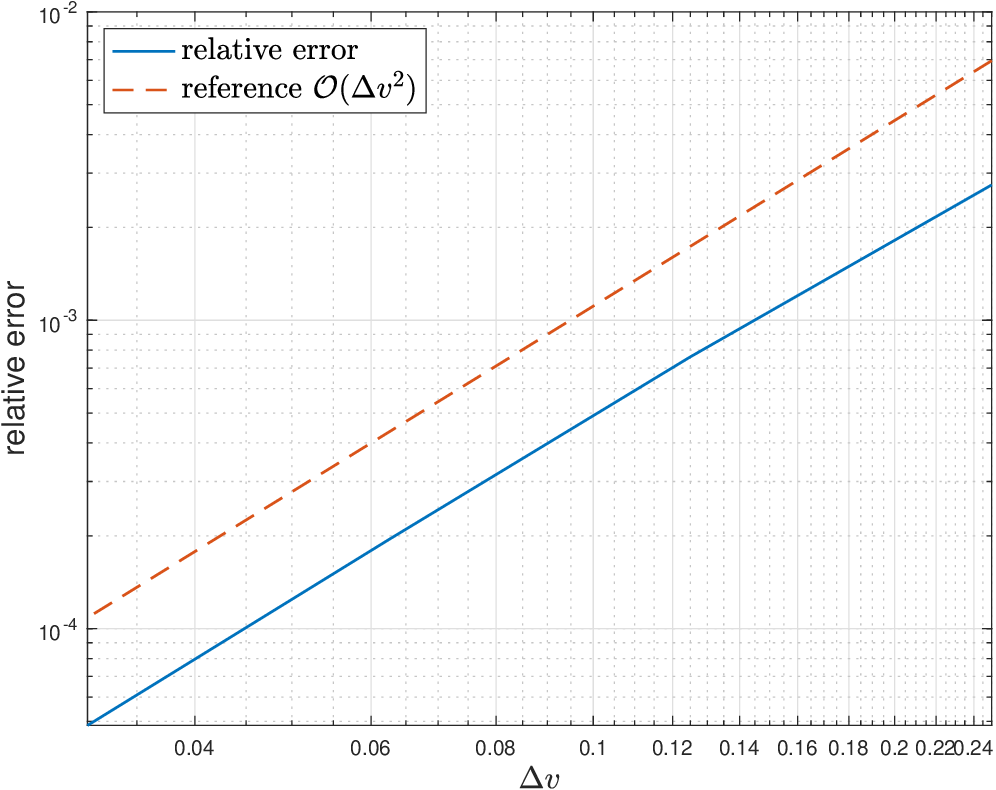}
    \caption{Numerical error for the linear Fokker-Planck equation.}
    \label{fig_spatial_homo_Fokker_Planck}
\end{figure}

We solve the equation up to time $t=1$ with a fixed TT-rank $(5,5)$. 
The velocity domain is chosen as $[v_{\min},v_{\max}]^3 = [-8,8]^3$ with $\dv = \frac{1}{4},\frac{1}{8},\frac{1}{16},\frac{1}{32}$. The time step is chosen as $\dt = \frac{\dv^2}{12}$. The relative error is shown in \Cref{fig_spatial_homo_Fokker_Planck}, 
which clearly demonstrates the second-order accuracy in $\dv$, or first-order in $\dt$.

\subsection{Spatially inhomogeneous kinetic equation}
\label{sec_inhomogeneous_kinetic_equation}

In this section, we consider the spatially inhomogeneous kinetic equation in the 1D3V setting:
\begin{equation} \label{eq:inhomo_kinetic}
    \pdt f(t,x,\bv) + v^{(1)} \pdx f(t,x,\bv) = \mathcal{Q}[f](t,x,\bv),
\end{equation}
where the collision operator $\mathcal{Q}[f]$ is either the BGK operator 
\begin{equation} \label{eq:BGK_oper}
\mathcal{Q}[f](t,x,\bv)=\eta(\mathcal{M}[f]-f),
\end{equation}
or the kinetic Fokker-Planck operator
\begin{equation} \label{eq:FP_oper}
\mathcal{Q}[f](t,x,\bv)=\eta\nabla_\bv \cdot \left(
        \mathcal{M}[f]\nabla_\bv \left(\frac{f}{\mathcal{M}[f]}\right)\right).
\end{equation}
Here the Maxwellian distribution $\mathcal{M}[f]$ is still given by \eqref{eq:Max}. However, due to the presence of the transport term, $n$, $\bu$, and $T$ in \eqref{eq:moments} depend on both time $t$ and spatial variable $x$, and therefore need to be updated at every time step and at every spatial point. 

For equation \eqref{eq:inhomo_kinetic} with the BGK operator \eqref{eq:BGK_oper}, we consider the case where the collision strength $\eta$ can be strong, and therefore use a first-order IMEX scheme:
\begin{equation}
\label{eq_IMEX_scheme}
\frac{\bf^{j,n+1}-\bf^{j,n}}{\dt} + \left(v^{(1)}\mathbf{D}_x^{\mathrm{up}} \bf\right)^{j,n} = \eta \left(\bM^{j,n+1} - \bf^{j,n+1}\right),
\end{equation}
where $\left(v^{(1)}\mathbf{D}_x^{\mathrm{up}}\bf\right)^{j,n}$ is the second-order upwind scheme for the transport term:
\begin{equation}
\label{eq_upwind_adv}
\begin{aligned}
    &\left(v^{(1)}\mathbf{D}_x^{\mathrm{up}}\bf \right)^{j,n}
    = 
        \frac{v_{k_1}^+}{2\dx} 
        \left(
        3\bf^{j,n} - 4\bf^{j-1,n} + \bf^{j-2,n}
        \right)
        - \frac{v_{k_1}^-}{2\dx} 
        \left(
        -\bf^{j+2,n} + 4 \bf^{j+1,n} - 3 \bf^{j,n}
        \right),
\end{aligned}
\end{equation}
with $v_{k_1}^+ = \max\{v_{k_1},0\}$ and $v_{k_1}^- = \max\{-v_{k_1},0\}$ and periodic boundary condition.
When $\bf^{j,n}$ is in the TT form, $(v^{(1)}\mathbf{D}_x^{\mathrm{up}}\bf)^{j,n}$ can be also written in TT form using \textbf{operations 1, 2 and 5} in \Cref{sec_TT_operators}.

In scheme \eqref{eq_IMEX_scheme}, $\bM^{j,n+1}$ appears implicitly, but there is a standard trick to deal with this (c.f. \cite{HJL17}). By taking the discrete moments $\sum_{k_1,k_2,k_3=1}^{N_v} \cdot \ (1,\bv_k, |\bv_k|^2)^T\Delta v^3$, $k=(k_1,k_2,k_3)$ on both sides of the scheme and using that the BGK operator is conservative, one obtains
\begin{equation} \label{eq_moment}
\frac{\mathbf{U}^{j,n+1}-\mathbf{U}^{j,n}}{\dt} + \sum_{k_1,k_2,k_3=1}^{N_v}\left(v^{(1)}\mathbf{D}_x^{\mathrm{up}} \bf\right)^{j,n}(1,\bv_k,|\bv_k|^2)^T \Delta v^3 = 0,
\end{equation}
where $\mathbf{U}:=(n,n \bu, n|\bu|^2+3nT)^T=\sum_{k_1,k_2,k_3=1}^{N_v} {\bf}(1,\bv_k, |\bv_k|^2)^T\Delta v^3$. In this way, the macroscopic quantities $n^{j,n+1}$, $\bu^{j,n+1}$, and $T^{j,n+1}$ can be obtained first, and hence $\bM^{j,n+1}$ is known. Therefore, at the beginning of each time step, we first solve \eqref{eq_moment} to obtain $\bM^{j,n+1}$, which is then used in all substeps of the low-rank algorithm. In the implementation, we evaluate the moments from the TT representations of $\bf^{j,n}$ and $(v^{(1)}\mathbf{D}_x^{\mathrm{up}\bf})^{j,n}$ using \textbf{operation 6} in \Cref{sec_TT_operators}. 

We then consider the first-order low-rank algorithm for the BGK equation. If we use the above first-order IMEX scheme in both forward and backward sub-projection steps, the mappings in \eqref{eq_distribution_propagation}-\eqref{eq_distribution_backward_propagation} are given by 
\begin{equation}
\label{eq_IMEX_IMEX}
\begin{split}
&\mathcal{G}^j_{\Delta t}(\bf^{1},\cdots,\bf^{N_x}) = \frac{1}{1+\dt\eta}
\left({\bf^{j} - \dt \left(v^{(1)}\mathbf{D}_x^{\mathrm{up}}\bf\right)^{j}}+\dt\eta \bM^{j,n+1}\right), \\
& \mathcal{H}^j_{\Delta t}(\bf^{1},\cdots,\bf^{N_x}) = \frac{1}{1-\dt\eta}
\left({\bf^{j} + \dt \left(v^{(1)}\mathbf{D}_x^{\mathrm{up}}\bf\right)^{j}} - \dt\eta \bM^{j,n+1}\right).
\end{split}
\end{equation}
Another choice is to use IMEX in forward substeps and forward Euler in backward substeps, where $\mathcal{G}^j_{\Delta t}$ remains the same and $\mathcal{H}^j_{\Delta t}$ becomes
\begin{equation}
\mathcal{H}_{\Delta t}^j(\bf^{1},\cdots,\bf^{N_x}) = 
\bf^{j} + \dt \left(v^{(1)}\mathbf{D}_x^{\mathrm{up}}\bf\right)^{j} - \dt\eta \left(\bM^{j,n+1} - \bf^{j}\right).
\end{equation}

\begin{figure}[t]
     \centering
     \begin{subfigure}[b]{0.32\textwidth}
         \centering
         \includegraphics[width=\textwidth]{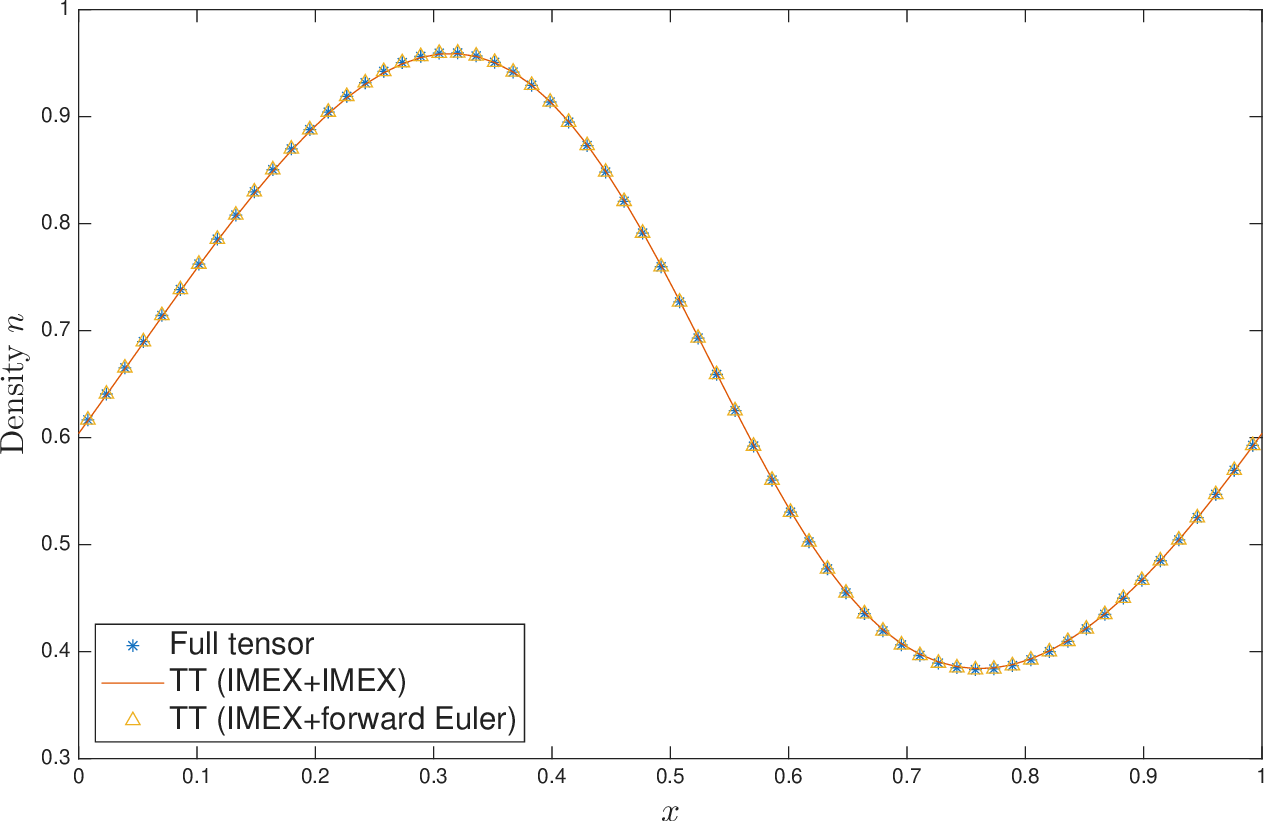}
         \caption{Density at $t=0.1$.}
     \end{subfigure}
     \hfill
     \begin{subfigure}[b]{0.32\textwidth}
         \centering
         \includegraphics[width=\textwidth]{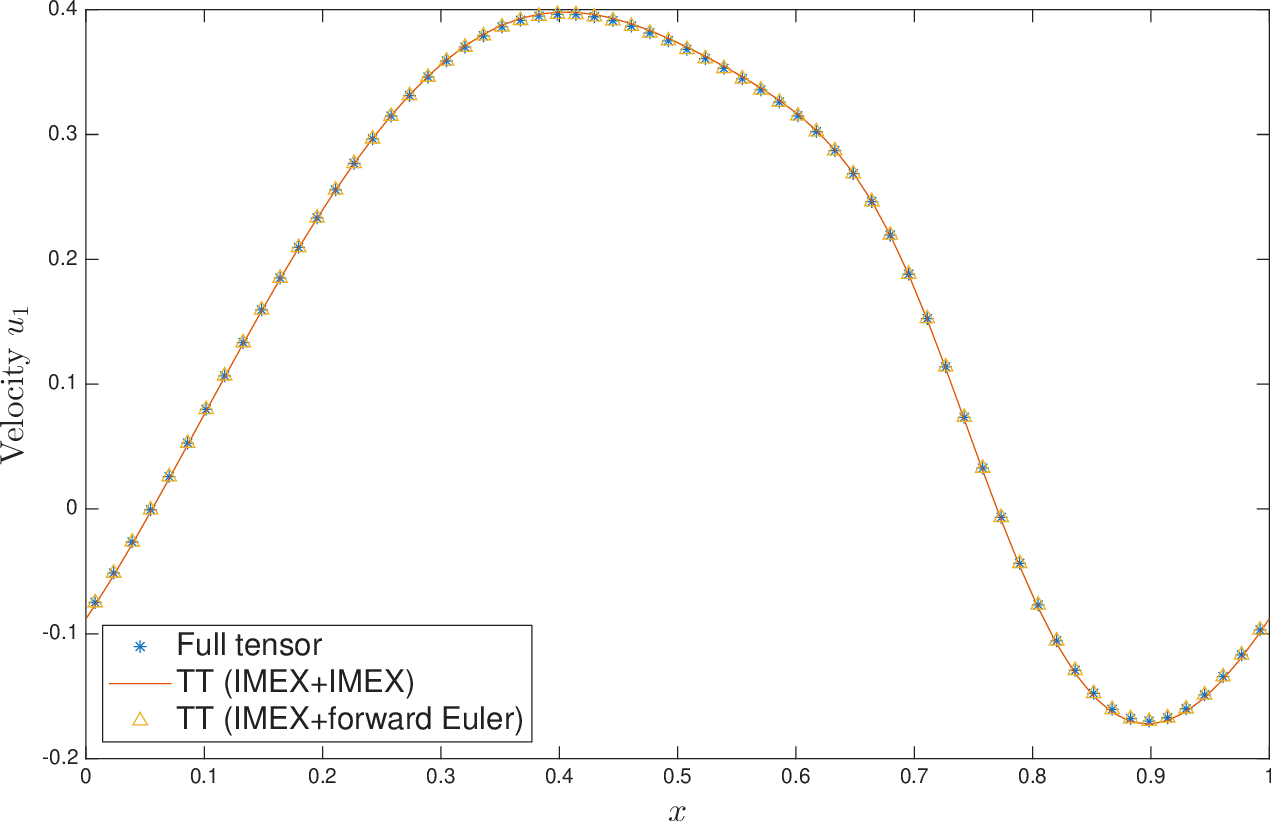}
         \caption{Velocity at $t=0.1$.}
     \end{subfigure}
     \hfill
     \begin{subfigure}[b]{0.32\textwidth}
         \centering
         \includegraphics[width=\textwidth]{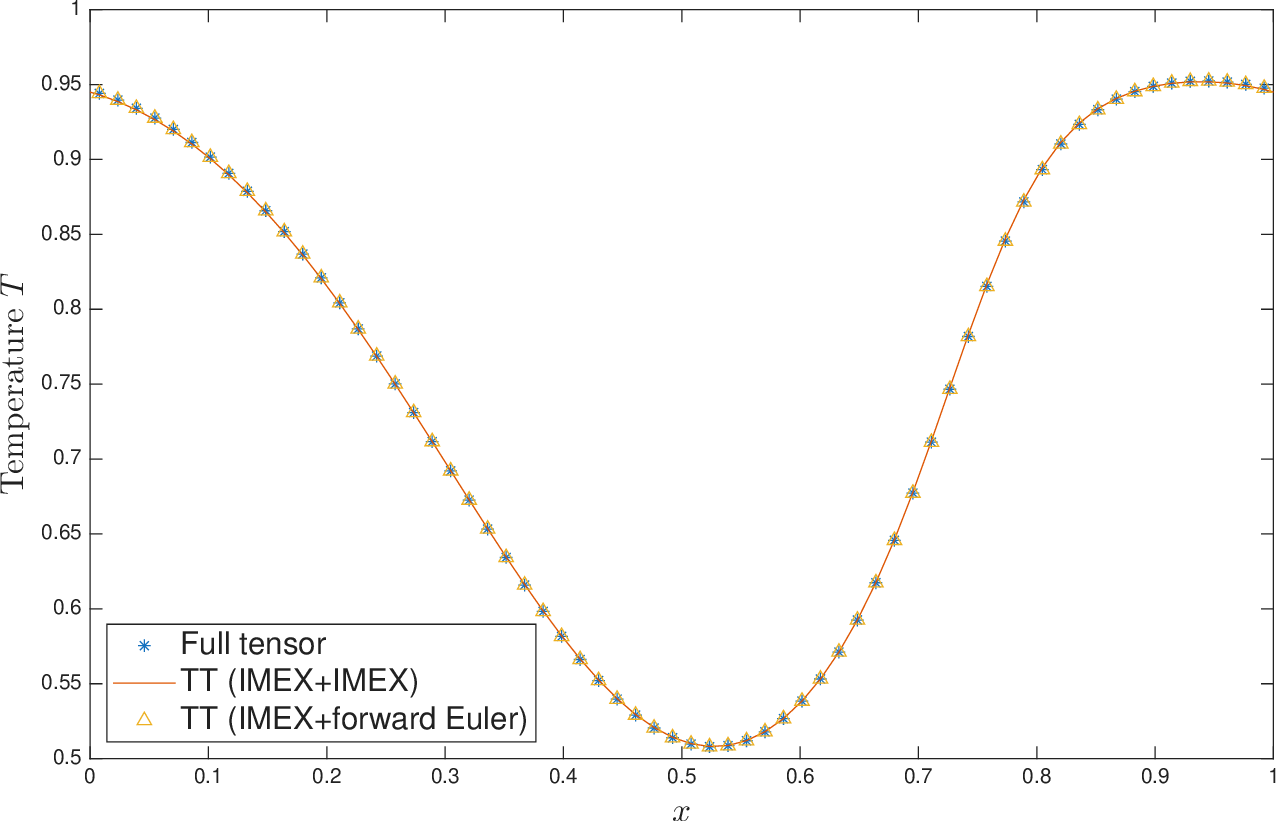}
         \caption{Temperature at $t=0.1$.}
     \end{subfigure}
        \caption{Density, bulk velocity, and temperature for the stiff spatially inhomogeneous BGK equation.}
        \label{fig_inhomo_BGK}
\end{figure}

In the numerical test, we consider
the initial condition given by
\begin{equation}
    f_0(x,\bv) = \frac{n_0(x)}{(2\pi T_0(x))^{3/2}} \exp\left(
        - \frac{\vert \bv - \bu_0 \vert^2}{2T_0(x)}
    \right),
\end{equation}
with
\begin{equation}
    n_0(x) = \frac{2+\sin(2\pi x)}{3}, \quad
    \bu_0 = (0.2,0,0)^T, \quad
    T_0(x) = \frac{3+\cos(2\pi x)}{4}.
\end{equation}
The spatial domain is chosen as $[0,1]$ with $N_x=64$, and the grid points are given by 
\begin{equation}
x_j = \left(j-\frac{1}{2}\right)\dx, \quad j=1,\cdots,N_x.
\end{equation}
The velocity domain is truncated to $[v_{\min},v_{\max}]^3 = [-6,6]^3$ with $N_v = 64$, and the grid points are given by
\eqref{eq_velocity_discretization}. Due to the IMEX treatment, we are able to consider a very stiff problem with $\eta = 10^5$. We set the time step to $\dt = 0.001$ and sovle the equation up to time $t=0.1$.
The TT-rank in the entire simulation is fixed as $(5,5)$. The macroscopic quantities $n$, $u^{(1)}$ (the first component of $\bu$), and $T$ at the final time are shown in \Cref{fig_inhomo_BGK}. For comparison, the results obtained using the  full tensor numerical scheme \eqref{eq_IMEX_scheme} are also included as a reference. The results of three methods agree quite well. 

For equation \eqref{eq:inhomo_kinetic} with the kinetic Fokker-Planck operator \eqref{eq:FP_oper}, we consider a fully explicit scheme, since handling the stiff collision operator would require additional techniques in the TT format that are beyond the scope of the current work. As a result, we consider a relatively small collision strength, $\eta=1$.

In the first-order low-rank algorithm, we use the forward Euler time-stepping in both forward and backward sub-projections steps, i.e., the mappings in \eqref{eq_distribution_propagation}-\eqref{eq_distribution_backward_propagation} are given by 
\begin{equation}
\begin{split}
&\mathcal{G}_{\Delta t}^j(\bf^{1},\cdots,\bf^{N_x})
= \bf^{j} - \dt \left(v^{(1)}\mathbf{D}_x^{\mathrm{up}} \bf\right)^j + \dt \eta (\mathbf{Q}\bf)^j, \\
&\mathcal{H}_{\Delta t}^j(\bf^{1},\cdots,\bf^{N_x})
= \bf^{j} + \dt \left(v^{(1)}\mathbf{D}_x^{\mathrm{up}} \bf\right)^j - \dt \eta (\mathbf{Q}\bf)^j,
\end{split}
\end{equation}
where the transport operator $\left(v^{(1)}\mathbf{D}_x^{\mathrm{up}} \bf\right)^j$ is given by \eqref{eq_upwind_adv} and the Fokker-Planck operator is discretized as in \eqref{eq_Fokker_Planck_discretization}-\eqref{eq_Fokker_Planck_fluxes}, except that the Maxwellian used is $\mathbf{M}^{j,n}$.
Specifically, at the beginning of each time step, we first evaluate the moments $\mathbf{U}^{j,n}$ based on the distribution $\bf^{j,n}$.
The Maxwellian $\mathbf{M}^{j,n}$ is then constructed using $\mathbf{U}^{j,n}$ and used in all substeps of the low-rank algorithm.

\begin{figure}[t]
     \centering
     \begin{subfigure}[b]{0.32\textwidth}
         \centering
         \includegraphics[width=\textwidth]{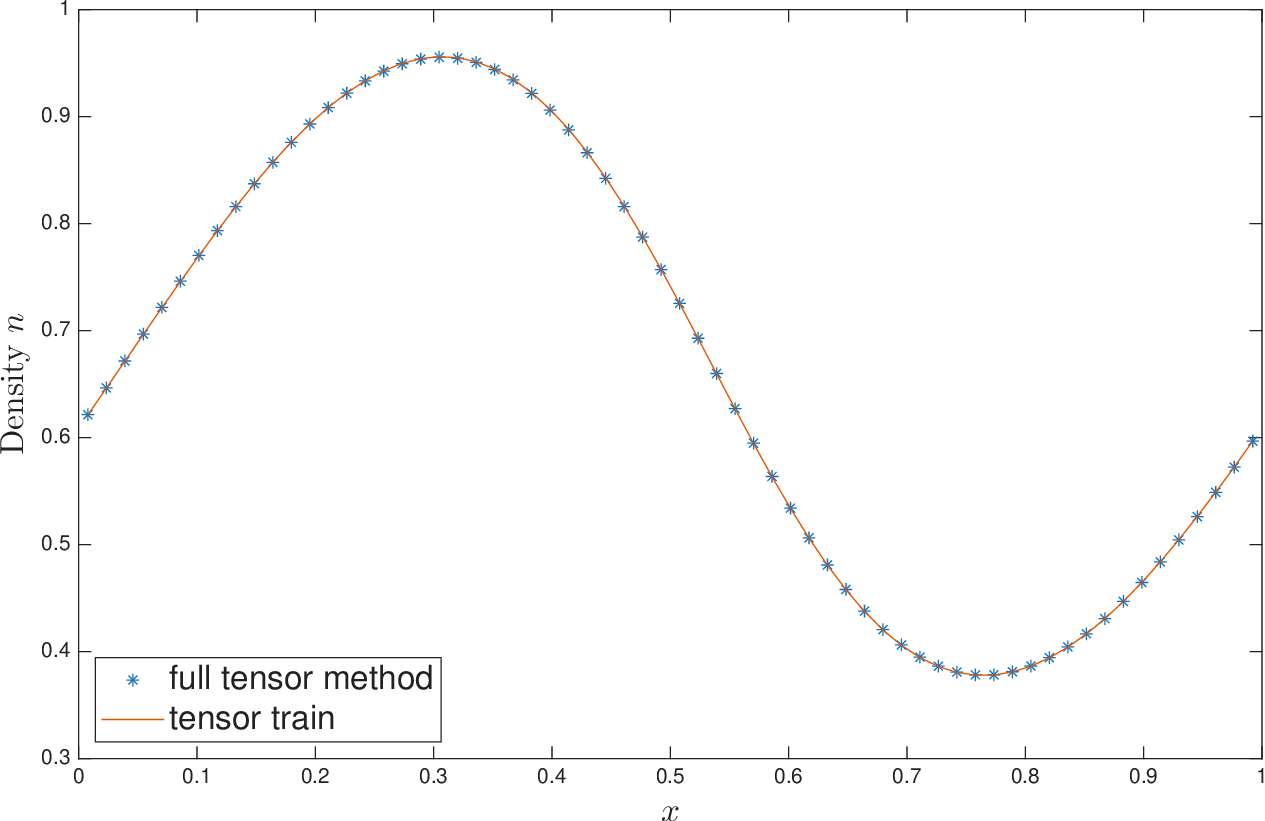}
         \caption{Density at $t=0.1$.}
     \end{subfigure}
     \hfill
     \begin{subfigure}[b]{0.32\textwidth}
         \centering
         \includegraphics[width=\textwidth]{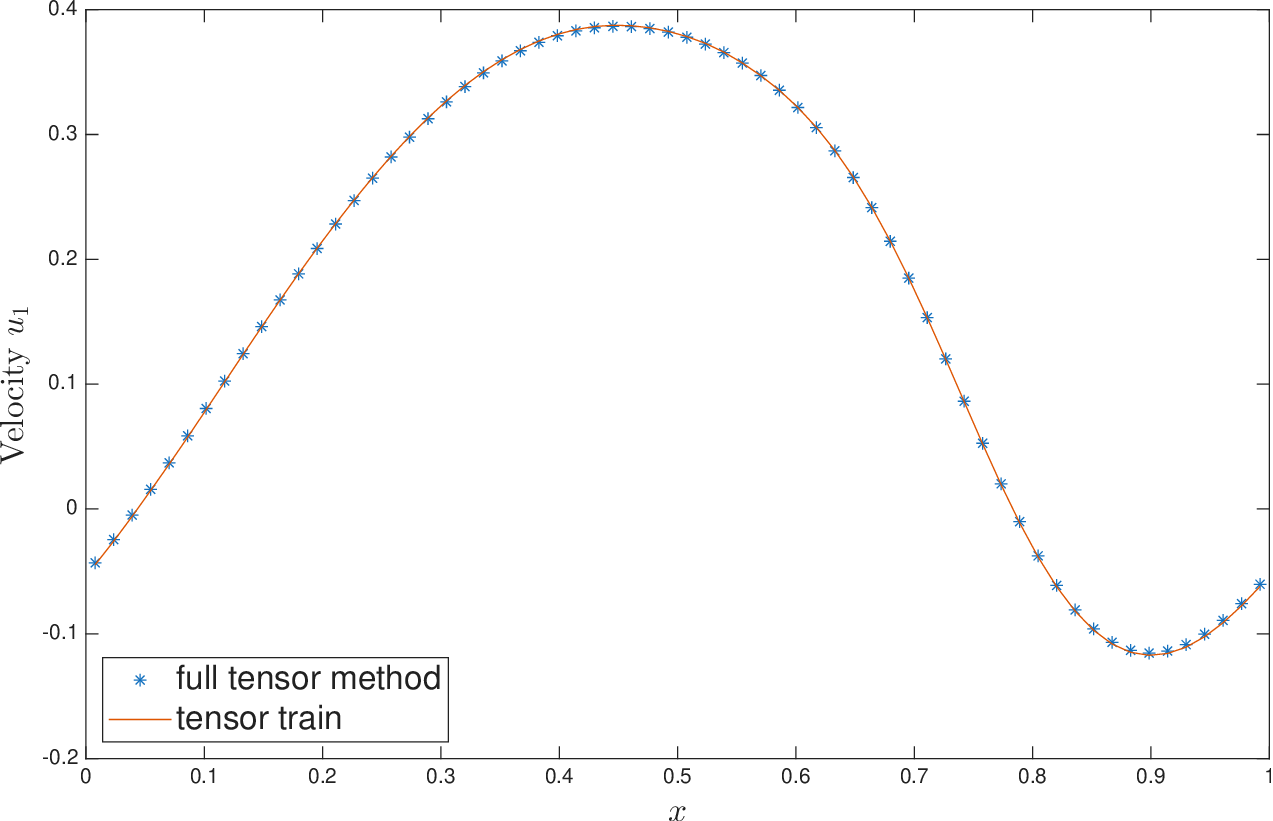}
         \caption{Velocity at $t=0.1$.}
     \end{subfigure}
     \hfill
     \begin{subfigure}[b]{0.32\textwidth}
         \centering
         \includegraphics[width=\textwidth]{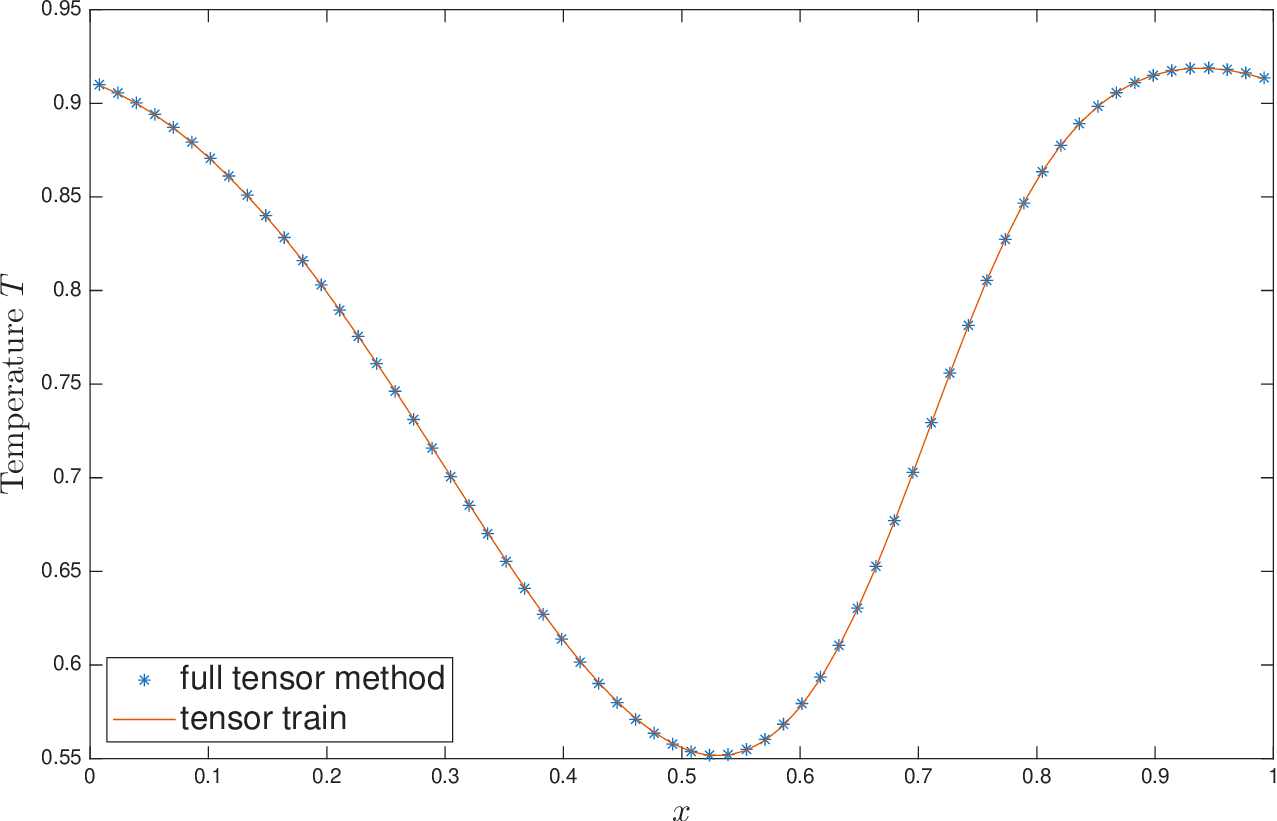}
         \caption{Temperature at $t=0.1$.}
     \end{subfigure}
        \caption{Density, bulk velocity, and temperature for the spatially inhomogeneous Fokker-Planck equation.}
        \label{fig_inhomo_FP}
\end{figure}

We use the same spatial and velocity discretization as in the spatially inhomogeneous BGK equation, and choose the time step as $\dt = 0.1 \min\left(\dfrac{\dx}{\max\vert v^{(1)} \vert}, \dfrac{\dv^2}{6\eta}\right)$. 
The equation is again solved to $t=0.1$ and the macroscopic quantities at $t=0.1$ are shown in \Cref{fig_inhomo_FP}.
The TT-rank in the entire simulation is fixed as (5,5).
The full tensor results are also included for reference. 

To conclude this subsection, we compare the full-tensor method and the proposed low-rank method in terms of memory requirements and computational cost. 
First, the full-tensor method requires storing a tensor of size $64\times 64 \times 64 \times 64$, which amounts to approximately 134 MB of memory.
In contrast, the low-rank method uses only about 1.18 MB of memory. This indicates that when higher-resolution solutions are needed, the full-tensor approach is likely to become computationally intractable, whereas the TT representation offers much greater flexibility and efficiency. Second, under the same experimental conditions\footnote{The experiments were conducted on a MacBook Pro equipped with an Apple M4 CPU.}, the simulation times (from $t=0$ to $t=0.1$) for both methods are reported in \Cref{table_simulation_time}.
For the BGK equation, the ``low-rank method'' in \Cref{table_simulation_time} refers to the IMEX–IMEX scheme; the IMEX–forward Euler scheme requires a comparable amount of simulation time.
The comparison shows that our method achieves accurate results in significantly less time. 
\begin{table}[ht]
    \centering
    \begin{tabular}{|c|c|c|}
        \hline
        time used (in seconds) & BGK & Fokker-Planck \\ \hline
        low-rank method & 7.7979 & 31.4565 \\ \hline
        full tensor method & 33.6206 & 254.2530 \\ \hline
    \end{tabular}
    \caption{Computational time of the low-rank method and full tensor method.}
    \label{table_simulation_time}
\end{table}

\subsection{Vlasov-Amp\`ere-Fokker-Planck (VAFP) equation}
\label{sec_VA_DFP}
In this section, we apply our method to the Vlasov-Fokker-Planck equation coupled with Ampère's law. This is a widely used kinetic model for plasma dynamics.

The full system reads as follows:
\begin{equation}
\label{eq_vlasov_ampere}
    \pdt f + \bv \cdot \gradx f - \bE \cdot \gradv f = \eta\, T\,\gradv \cdot \left(\mathcal{M}[f]\gradv\left(\frac{f}{\mathcal{M}[f]}\right)\right),
\end{equation}
where $f(t,\bx,\bv)$ is the distribution function of electrons depending on time $t$, position $\bx$, and velocity $\bv$. $\bE$ is the electric field determined by the Gauss's law
\begin{equation}
    \gradx \cdot \bE(t,\bx) = \rho(t,\bx) - \rho_i,
\end{equation}
where $\rho(t,\bx) = -\int_{\mathbb{R}^3} f \dd \bv$ is the charge density and $\rho_i$ is a uniform background density satisfying
$
     \int_{\Omega_\bx} \left( \rho(t,\bx) - \rho_i\right) \dd \bx = 0.
$
The electric field $\bE$ also follows the Amp\`ere's law
\begin{equation}
\label{eq_ampere_law}
    \pdt \bE(t,\bx) = - \bJ(t,\bx),
\end{equation}
where 
$
    \bJ(t,\bx) = -\int_{\mathbb{R}^3} \bv f \dd \bv
$
is the current density. It can be shown that in the continuous case, the Gauss's law and the Amp\`ere's law are equivalent. 

The right-hand side of \eqref{eq_vlasov_ampere} is the Fokker-Planck operator with the Maxwellian and moments defined by \eqref{eq:Max}-\eqref{eq:moments}. Note that it can be equivalently written as $\eta \nabla_{\bv}\cdot \left(T\nabla_{\bv} f+(\bv-\bu)\right)$, which is the more familiar Dougherty operator \cite{Doughterty64} often appearing in the physics literature.

\subsubsection{Numerical discretization of the VAFP equation}

Following the previous sections, we limit our discussion of the VAFP equation to the 1D3V setting, that is, $f=f(t,x,\bv)$, $\bv=(v^{(1)},v^{(2)},v^{(3)})$. The equation is then reduced to
\begin{equation}
\partial_t f + v^{(1)} \partial_x f - E^{(1)} \partial_{v^{(1)}}f
= \eta\, T\,\gradv \cdot \left(\mathcal{M}[f]\gradv\left(\frac{f}{\mathcal{M}[f]}\right)\right),
\end{equation}
where $E^{(1)}$ is the first component of $\bE$. We initialize the electric field by solving the Gauss's law and then evolve it in time by solving the Amp\`ere's law.

For both transport terms in $x$ and $v^{(1)}$, we use the second-order upwind scheme. That is, $v^{(1)}\partial_x f$ is discretized using \eqref{eq_upwind_adv} with periodic boundary condition. $E^{(1)}\partial_{v^{(1)}}f$ is discretized as
\begin{equation}
\begin{aligned}
&\left(E^{(1)}\mathbf{D}_{v^{(1)}}^{\mathrm{up}}\bf \right)^{j}_{k_1k_2k_3}
    = 
        \frac{(E^{(1)}_j)^+ (-f_{k_1+2,k_2k_3}^j+4f_{k_1+1,k_2k_3}^j-3f_{k_1k_2k_3}^j) 
    - (E_{j}^{(1)})^- (3f_{k_1k_2k_3}^j-4f_{k_1-1,k_2k_3}^j+f_{k_1-2,k_2k_3}^j)}{2\dv}
\end{aligned}
\end{equation}
with $(E_j^{(1)})^+ = \max\{(E^{(1)})^j,0\}$ and $(E_j^{(1)})^- = \max\{-(E^{(1)})^j,0\}$.
We apply the zero boundary condition in the $v^{(1)}$ direction:
\begin{equation}
    f_{0,k_2k_3}^j = f_{-1,k_2k_3}^j = f_{N_v+1,k_2k_3}^j = f_{N_v+2,k_2k_3}^j = 0.
\end{equation}

At each time step, we first compute the density $n^{j,n}$, bulk velocity $\bu^{j,n}$, temperature $T^{j,n}$ and current density $(J^{(1)})^{j,n}$ (the first component of $\bJ$) using $\bf^{j,n}$.
With these macroscopic quantities, we are able to compute $(E^{(1)})^{j,n+1}$  by solving \eqref{eq_ampere_law} with the forward Euler method and construct the Maxwellian $\bM^{j,n}$.
The Fokker-Planck operator is then discretized the same as in \Cref{sec_inhomogeneous_kinetic_equation} using $\bM^{j,n}$. 

In the first-order low-rank algorithm, we use the forward Euler time-stepping in both forward and
backward sub-projections steps, i.e., the mappings in \eqref{eq_distribution_propagation}-\eqref{eq_distribution_backward_propagation} are given by 
\begin{equation}
\begin{split}
&\mathcal{G}_{\Delta t}^j(\bf^1,\cdots,\bf^{N_x})
= \bf^j - \dt \left(v^{(1)} \mathbf{D}_x^{\mathrm{up}}\bf\right)^j
+ \dt \left(E^{(1)} \mathbf{D}_{v^{(1)}}^\mathrm{up} \bf\right)^j  + \dt \eta (\mathbf{Q}\bf)^j, \\
&\mathcal{H}_{\Delta t}^j(\bf^1,\cdots,\bf^{N_x})
= \bf^j + \dt \left(v^{(1)}\mathbf{D}_x^{\mathrm{up}}\bf\right)^j
- \dt \left(E^{(1)} \mathbf{D}_{v^{(1)}}^\mathrm{up} \bf\right)^j  - \dt \eta (\mathbf{Q}\bf)^j.
\end{split}
\end{equation}

\subsubsection{Numerical results of the VAFP equation}

We consider two benchmark tests for the VAFP equation: the linear Landau damping and the two-stream instability. \\

\noindent \textbf{Linear Landau damping}
\begin{figure}
    \centering
    \includegraphics[width=0.6\linewidth]{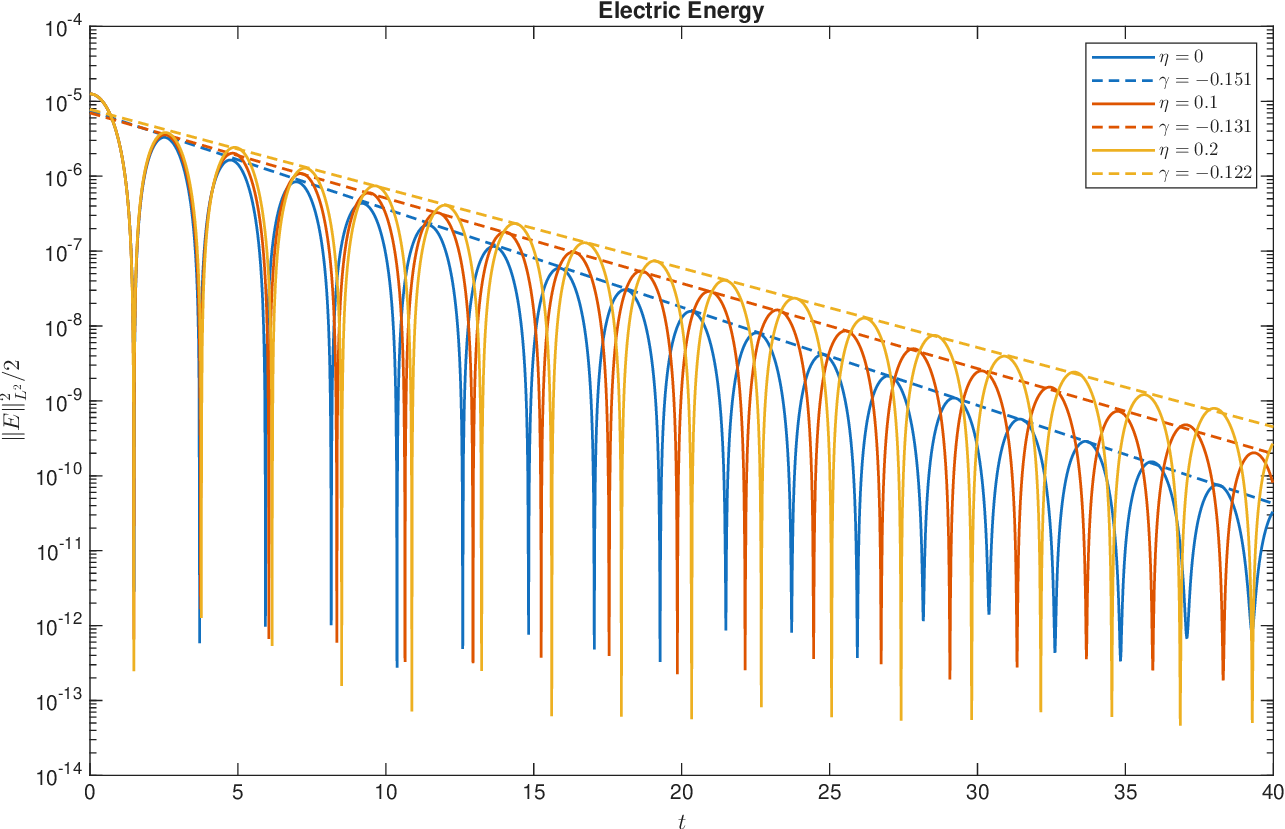}
    \caption{Linear Landau damping. Evolution of electric energy and damping rates for different collision strengths.}
    \label{fig_LLD_electric_energy}
\end{figure}
\begin{figure}
     \centering
     \begin{subfigure}[b]{0.3\textwidth}
         \centering
         \includegraphics[width=\textwidth]{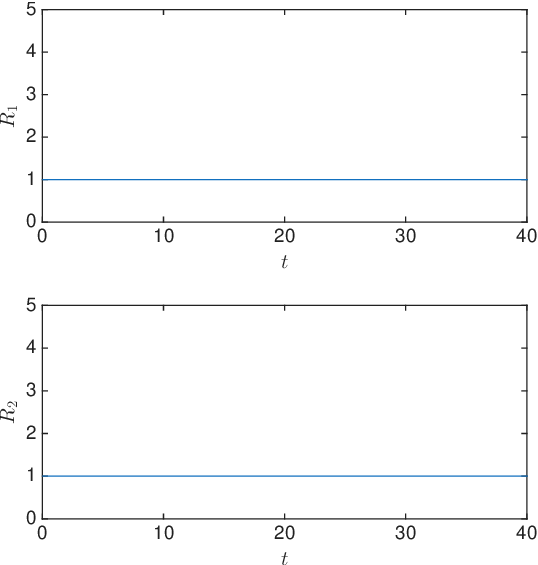}
         \caption{$\eta=0$}
     \end{subfigure}
     \hspace{-1mm}
     \begin{subfigure}[b]{0.3\textwidth}
         \centering
         \includegraphics[width=\textwidth]{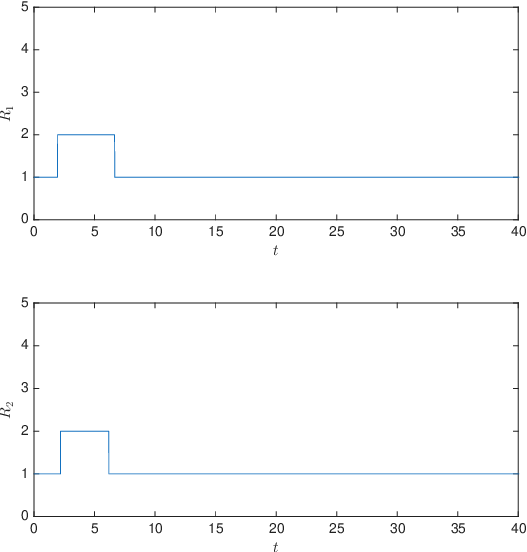}
         \caption{$\eta =0.1$}
     \end{subfigure}
     \hspace{-1mm}
     \begin{subfigure}[b]{0.3\textwidth}
         \centering
         \includegraphics[width=\textwidth]{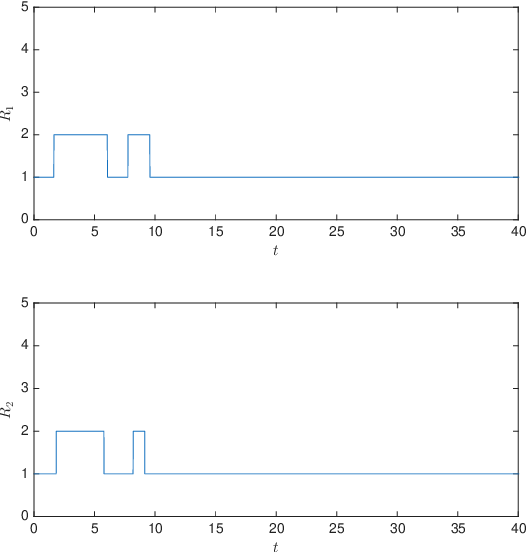}
         \caption{$\eta = 0.2$}
     \end{subfigure}
        \caption{Linear Landau damping. Evolution of effective ranks for different collision strengths. 
        }
        \label{fig_LLD_rank}
\end{figure}
We first consider the linear Landau damping with initial condition 
\begin{equation}
    f(0,x,\bv) = \frac{1}{(2\pi)^{3/2}} (1+A\cos(\kappa x)) \e^{-v_1^2/2}\e^{-v_2^2/2}\e^{-v_3^2/2}.
\end{equation}
In this case, we can compute the initial electric field explicitly:
\begin{equation}
    E^{(1)}(0,x)= - \frac{A}{\kappa} \sin(\kappa x).
\end{equation}
The physical parameters are $A=0.001$ and $\kappa = 0.5$.
We choose the spatial domain as
$x\in[0,2\pi/\kappa]=[0,4\pi]$ with $N_x=128$ and periodic boundary condition.
The velocity domain is truncated to $[v_{\min},v_{\max}]^3=[-9,9]^3$ with $N_v=128$.
We choose the time step as 
\begin{equation}
\label{eq_time_step_CFL}
    \dt = 0.1 \min\left\{\frac{\dx}{\max |v^{(1)}|}, \frac{\dv}{\max \vert E^{(1)}(0,x)\vert}, \frac{\dv^2}{6 \eta} \right\}.
\end{equation}
The TT-rank is fixed as $(5,5)$ during the simulation.
We evaluate the effect of different collision strengths by choosing $\eta =0,0.1, 0.2$.
The electric energy is defined as
\begin{equation}
    \mathcal{E}(t) = \frac{1}{2}\int_{0}^{2\pi/\kappa} \left(E^{(1)}(t,x)\right)^2 \dd x
    \approx \frac{1}{2} \sum_{j=1}^{N_x} ((E^{(1)})^j)^2\dx,
\end{equation}
whose evolution is shown in \Cref{fig_LLD_electric_energy}.
We observed that the damping rate (in absolute value) decreases as the collision strength increases.
In the collisionless case, the damping rate is in good agreement with the linear theory prediction of $-0.153$.
Additionally, we track the \emph{effective rank} of the solution during the simulation.
For a tensor train $\bf^j$ given in \eqref{eq_TT_decomposition_forms1}, we define two effective ranks.
We compute the singular values $\sigma_1,\cdots,\sigma_{r_1}$ of the matrix $S^{j,(1)}$ when $\bf^j$ is in form \rombracket{2} and the first effective rank is defined as
\begin{equation}
\mathfrak{r}_1(\bf^j) = \max\{r; \ \sigma_r \geqslant \delta \sigma_1\}.
\end{equation}
Similarly, we can also compute the singular values of the matrix $S^{j,(2)}$ when $\bf^j$ is in form \rombracket{4} and define another effective rank $\mathfrak{r}_2$.
The final effective tensor ranks of the solution are defined as
\begin{equation}
    {R}_1 = \max_{j=1,\cdots,N_x} \mathfrak
    {r}_1(\bf^j),\quad
    {R}_2 = \max_{j=1,\cdots,N_x} \mathfrak
    {r}_2(\bf^j).
\end{equation}
The effective tensor ranks of the numerical solution, computed with threshold $\delta = 10^{-5}$, are presented in \Cref{fig_LLD_rank}.
For the case without collisions ($\eta=0$), our simulation maintains a TT-rank of $(1,1)$ throughout the simulation, as expected.
Even with collisions,
the effective ranks of the solution remain low, indicating that a small TT-rank is sufficient for this example.\\

\noindent \textbf{Two-stream instability}
\begin{figure}
    \centering
    \includegraphics[width=0.5\linewidth]{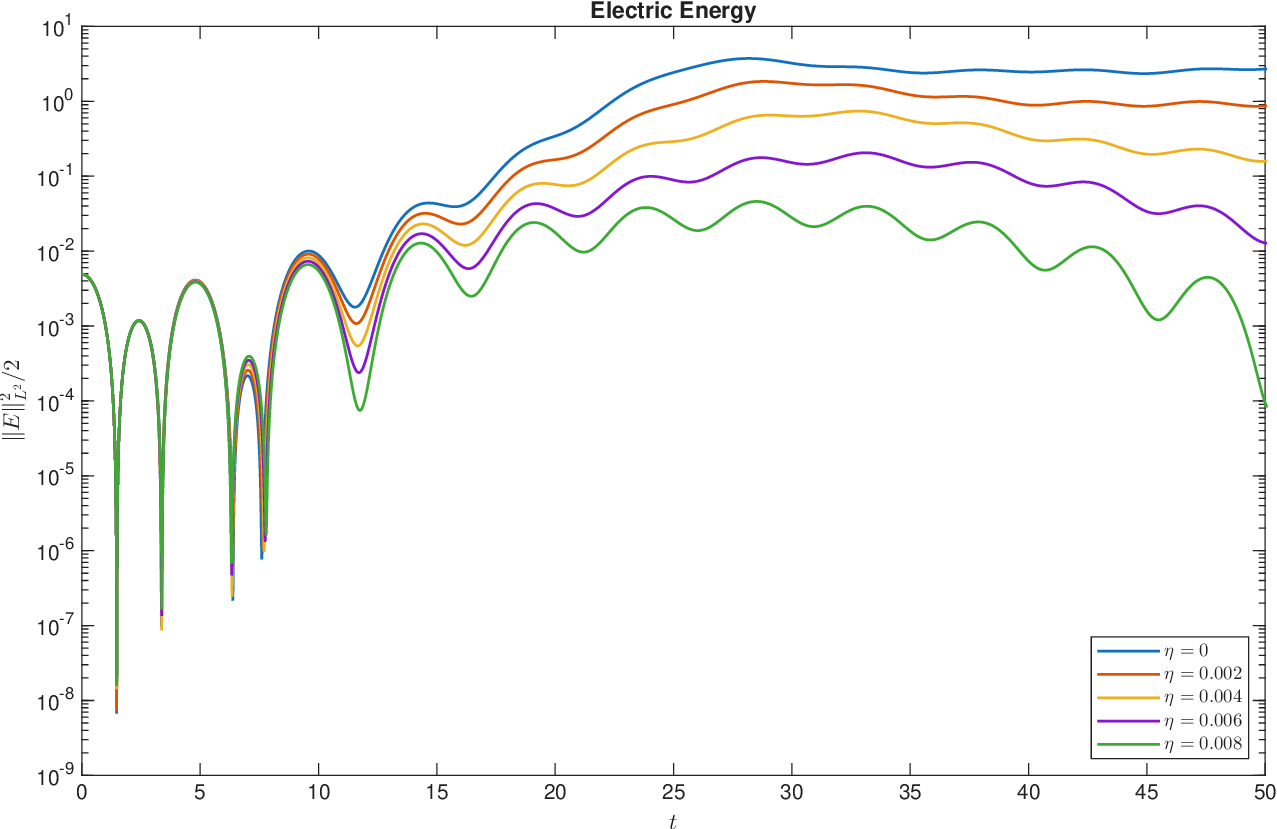}
    \caption{Two-stream instability. Evolution of electric energy for different collision strengths.}
    \label{fig_TSI_energy}
\end{figure}
\begin{figure}
     \centering
     \begin{subfigure}[b]{0.22\textwidth}
         \centering
         \includegraphics[width=\textwidth]{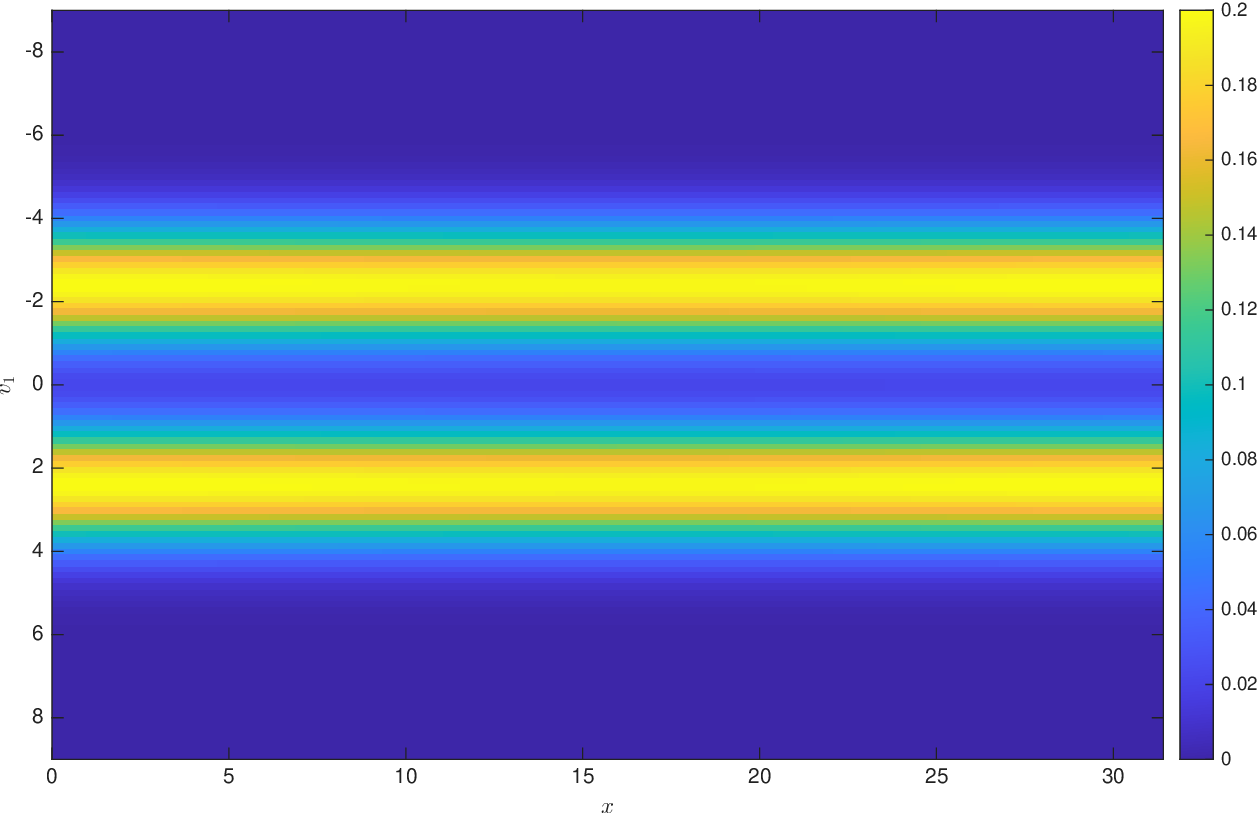}
         \caption{$\eta=0,t=0$}
     \end{subfigure}
     \begin{subfigure}[b]{0.22\textwidth}
         \centering
         \includegraphics[width=\textwidth]{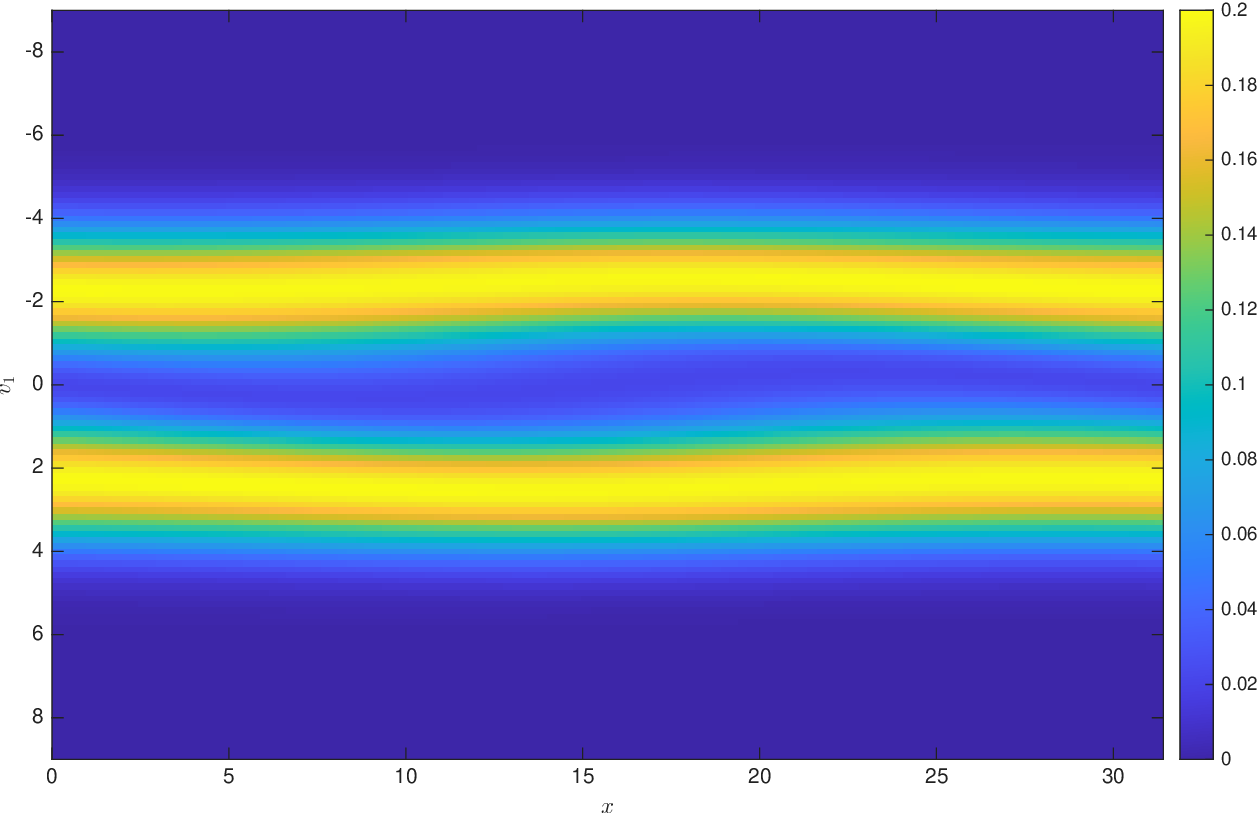}
         \caption{$\eta =0,t=15$}
     \end{subfigure}
     \begin{subfigure}[b]{0.22\textwidth}
         \centering
         \includegraphics[width=\textwidth]{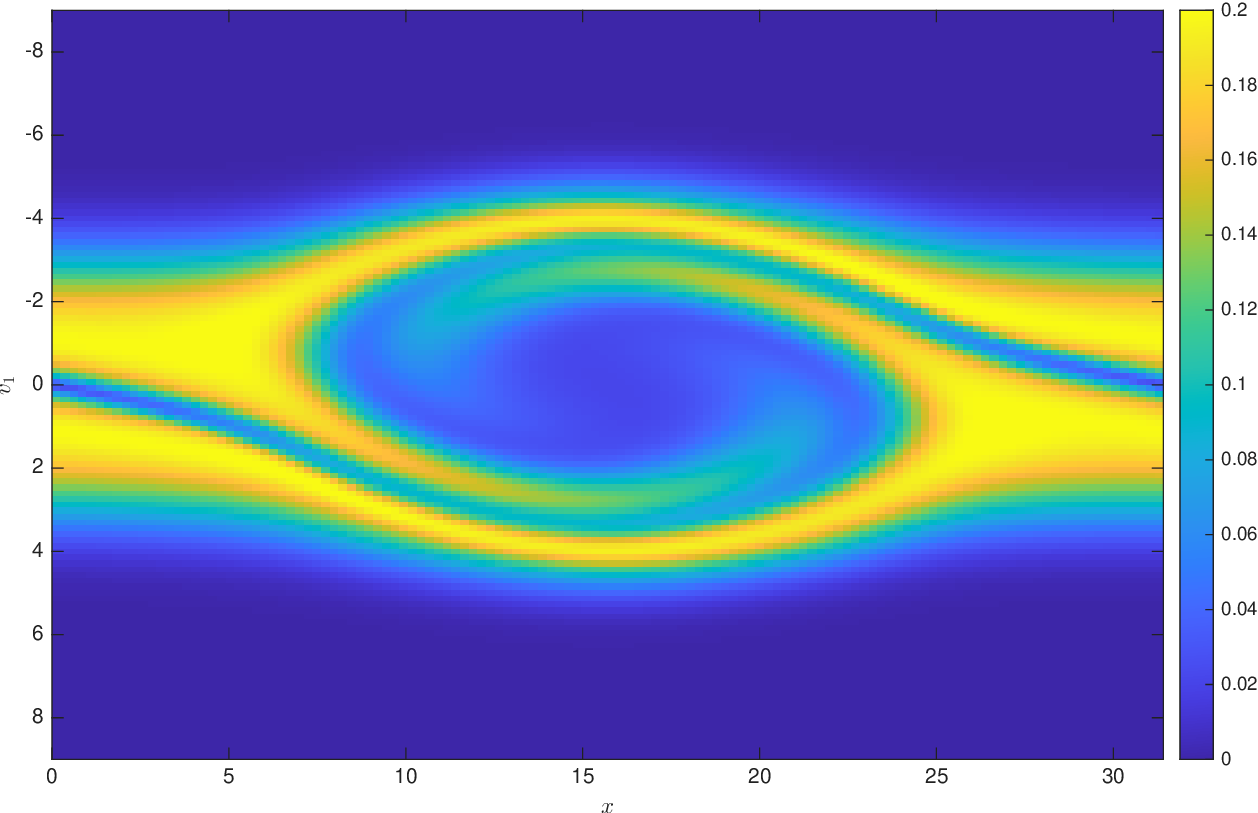}
         \caption{$\eta = 0,t=30$}
     \end{subfigure}
     \begin{subfigure}[b]{0.22\textwidth}
         \centering
         \includegraphics[width=\textwidth]{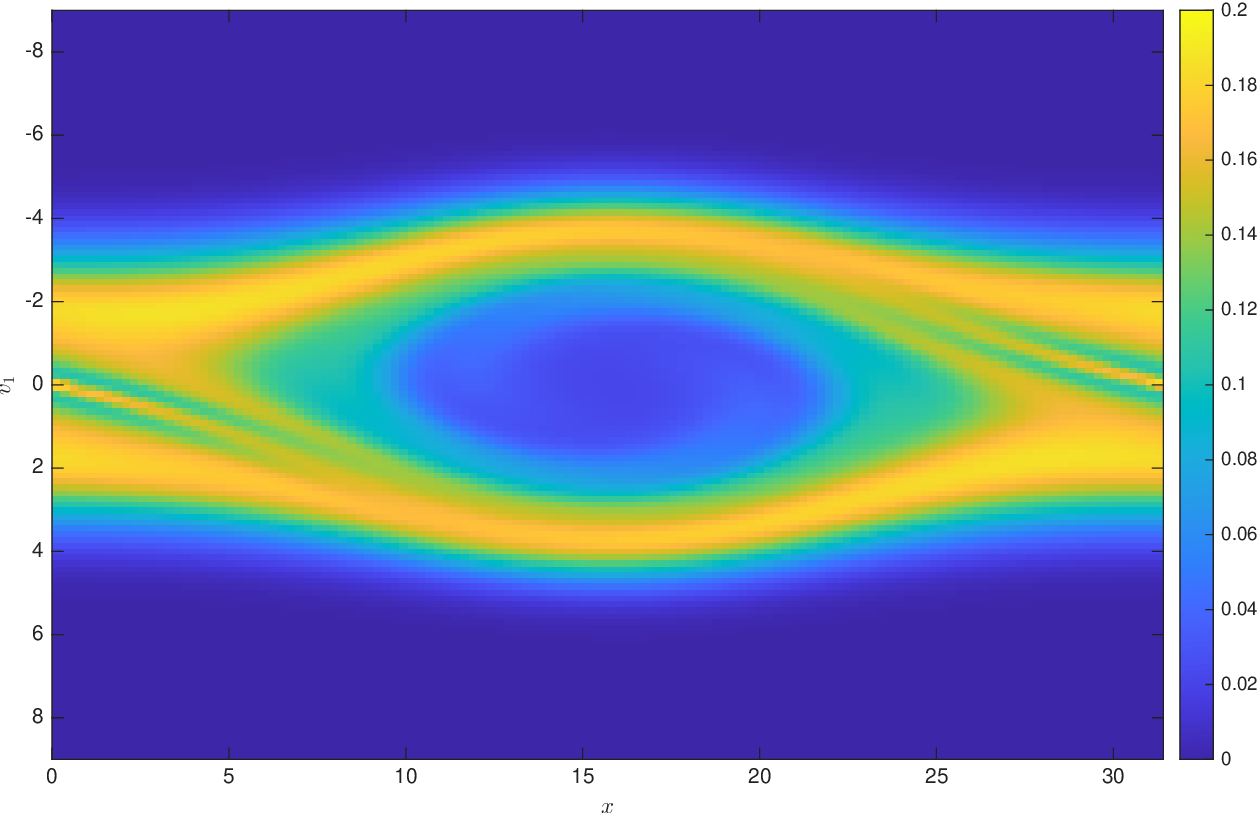}
         \caption{$\eta =0,t=45$}
     \end{subfigure} \\
     \begin{subfigure}[b]{0.22\textwidth}
         \centering
         \includegraphics[width=\textwidth]{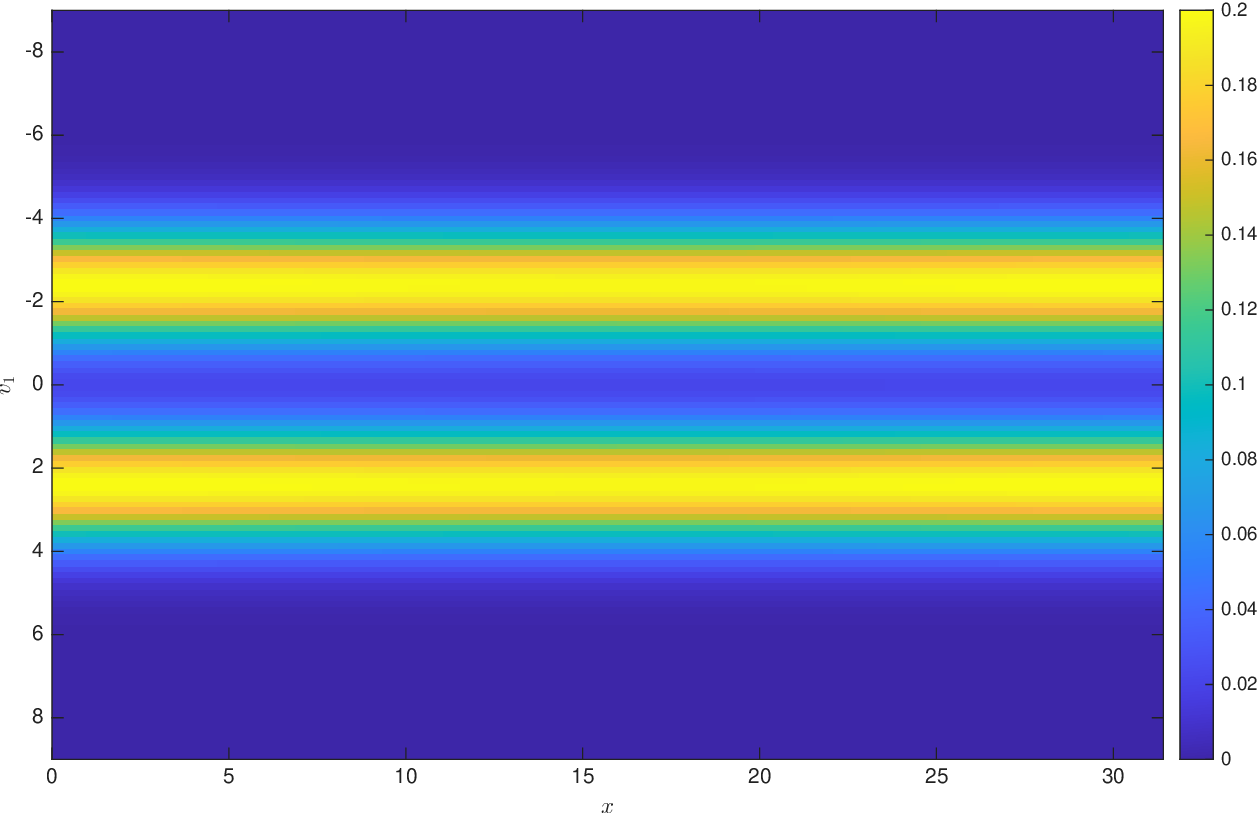}
         \caption{$\eta= 0.002,t=0$}
     \end{subfigure}
     \begin{subfigure}[b]{0.22\textwidth}
         \centering
         \includegraphics[width=\textwidth]{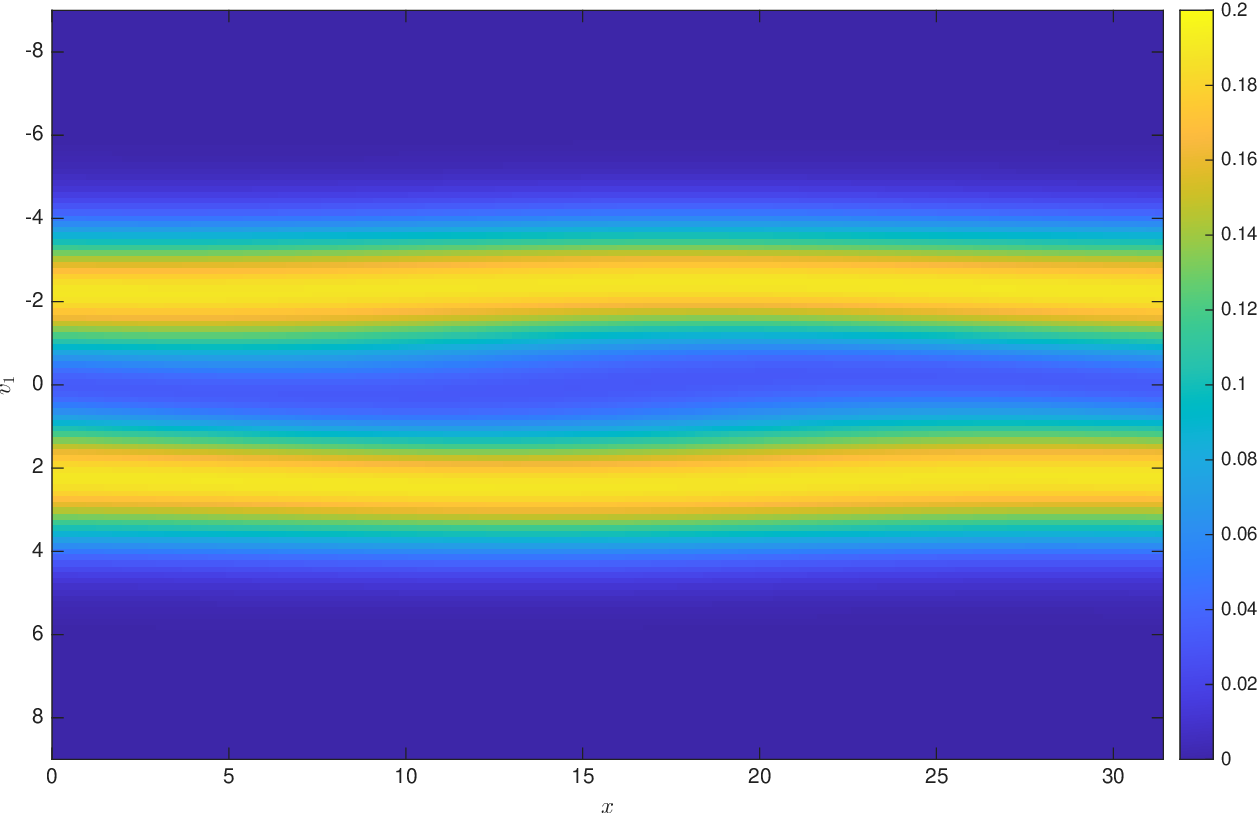}
         \caption{$\eta = 0.002,t=15$}
     \end{subfigure}
     \begin{subfigure}[b]{0.22\textwidth}
         \centering
         \includegraphics[width=\textwidth]{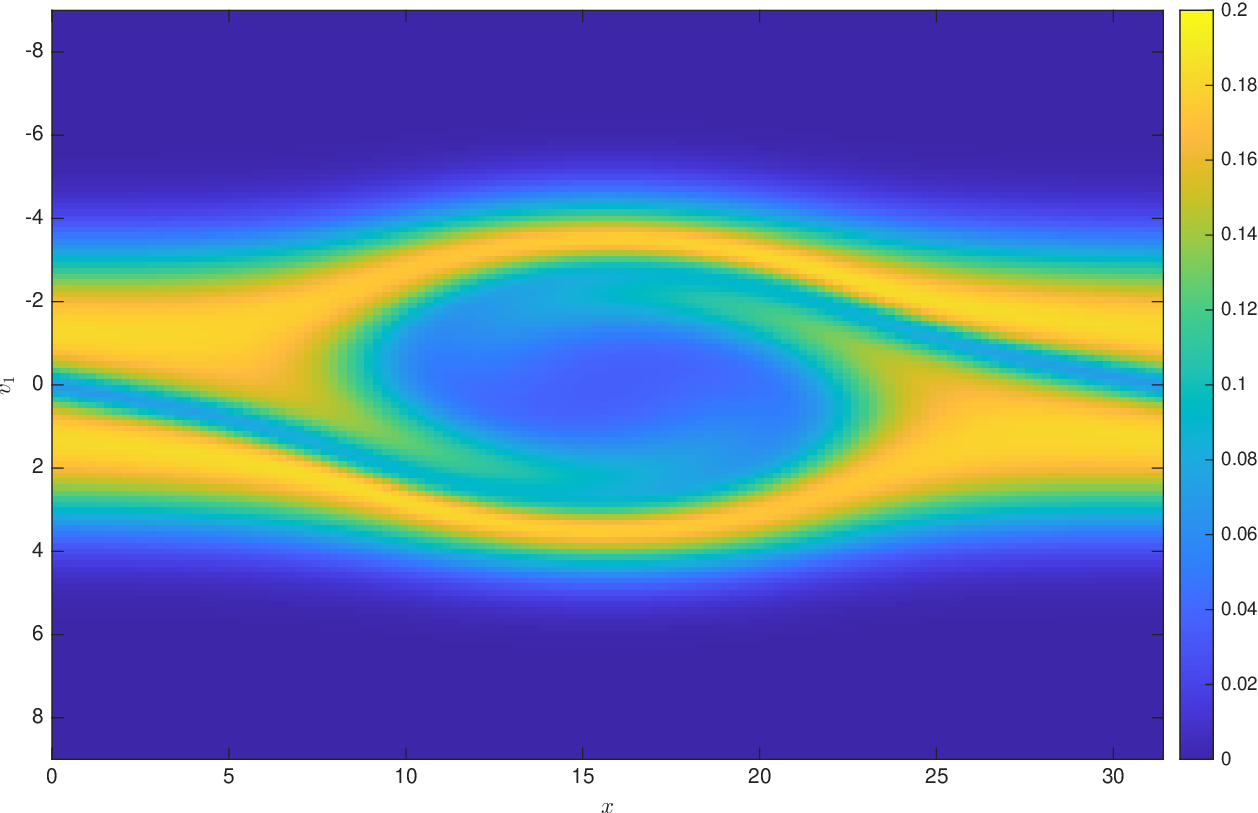}
         \caption{$\eta = 0.002,t=30$}
     \end{subfigure}
     \begin{subfigure}[b]{0.22\textwidth}
         \centering
         \includegraphics[width=\textwidth]{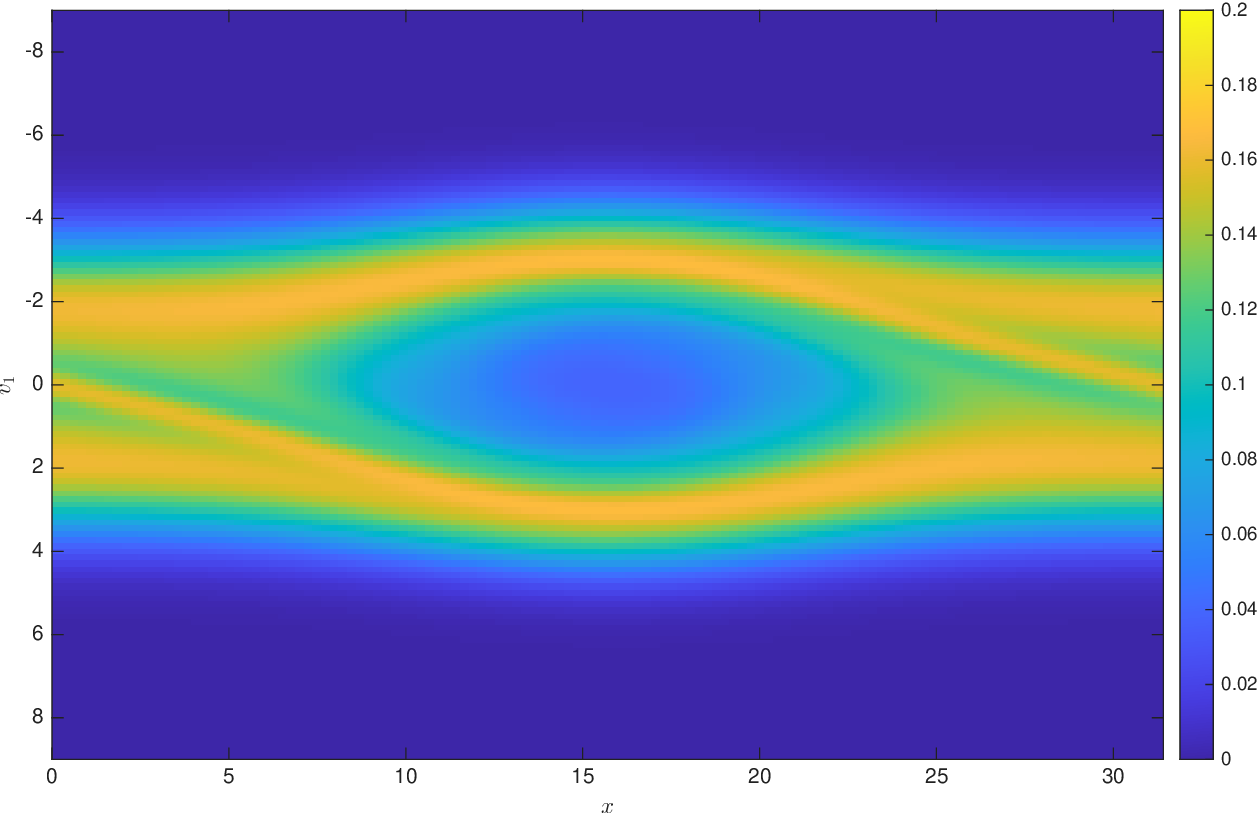}
         \caption{$\eta =0.002,t=45$}
     \end{subfigure}\\
     \begin{subfigure}[b]{0.22\textwidth}
         \centering
         \includegraphics[width=\textwidth]{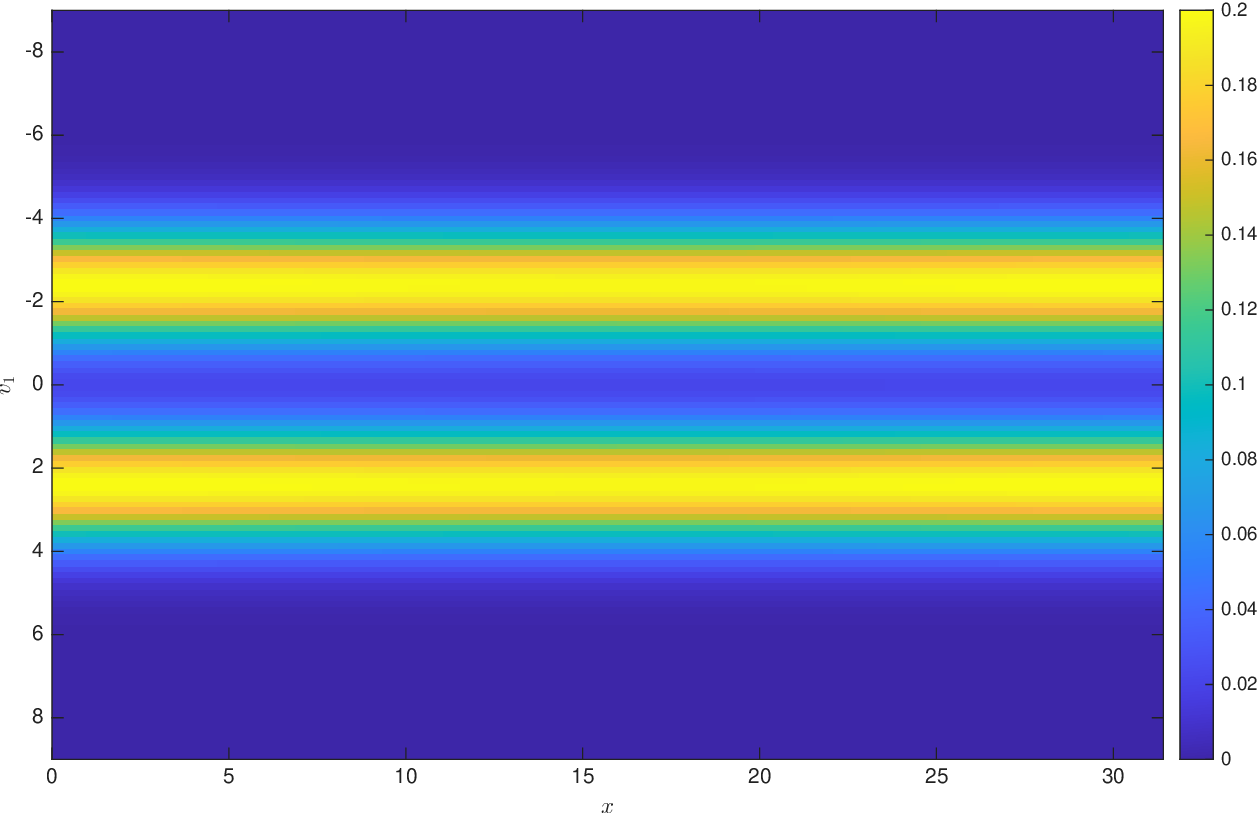}
         \caption{$\eta= 0.004,t=0$}
     \end{subfigure}
     \begin{subfigure}[b]{0.22\textwidth}
         \centering
         \includegraphics[width=\textwidth]{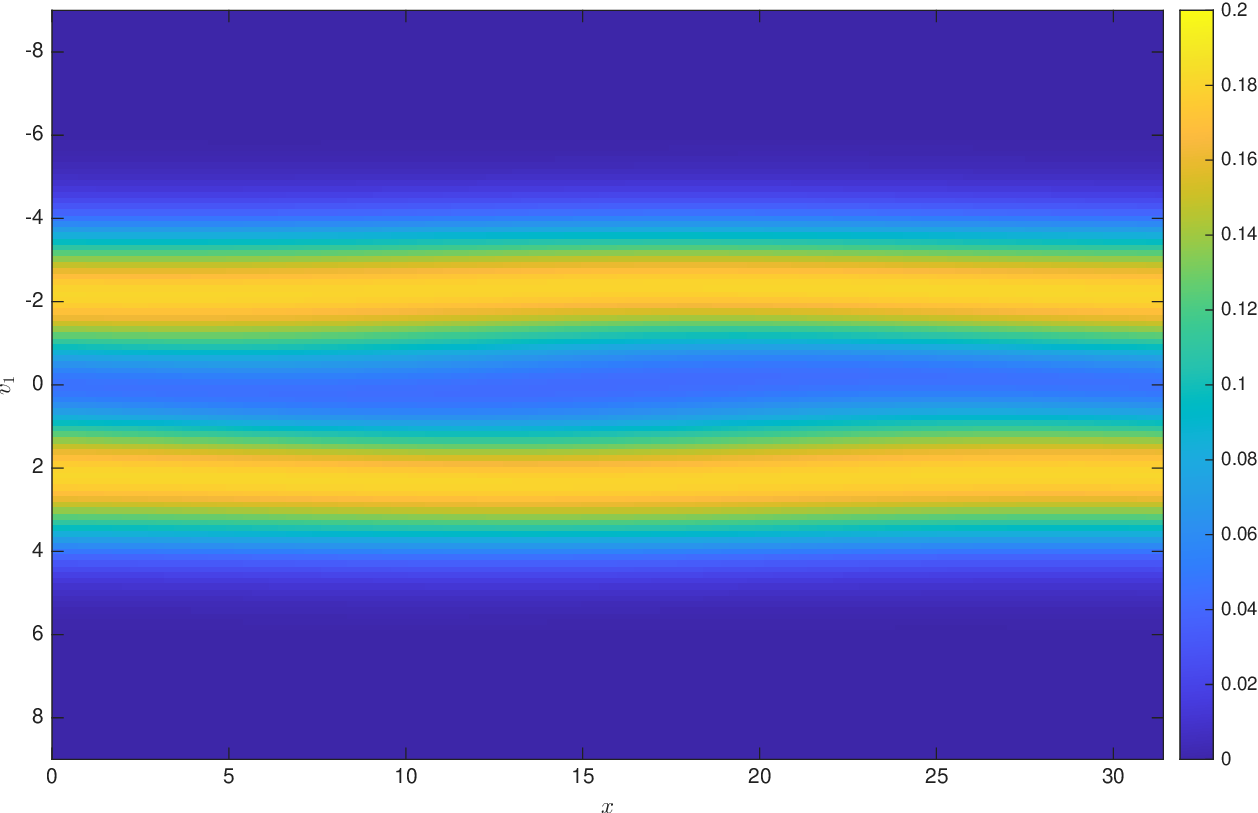}
         \caption{$\eta = 0.004,t=15$}
     \end{subfigure}
     \begin{subfigure}[b]{0.22\textwidth}
         \centering
         \includegraphics[width=\textwidth]{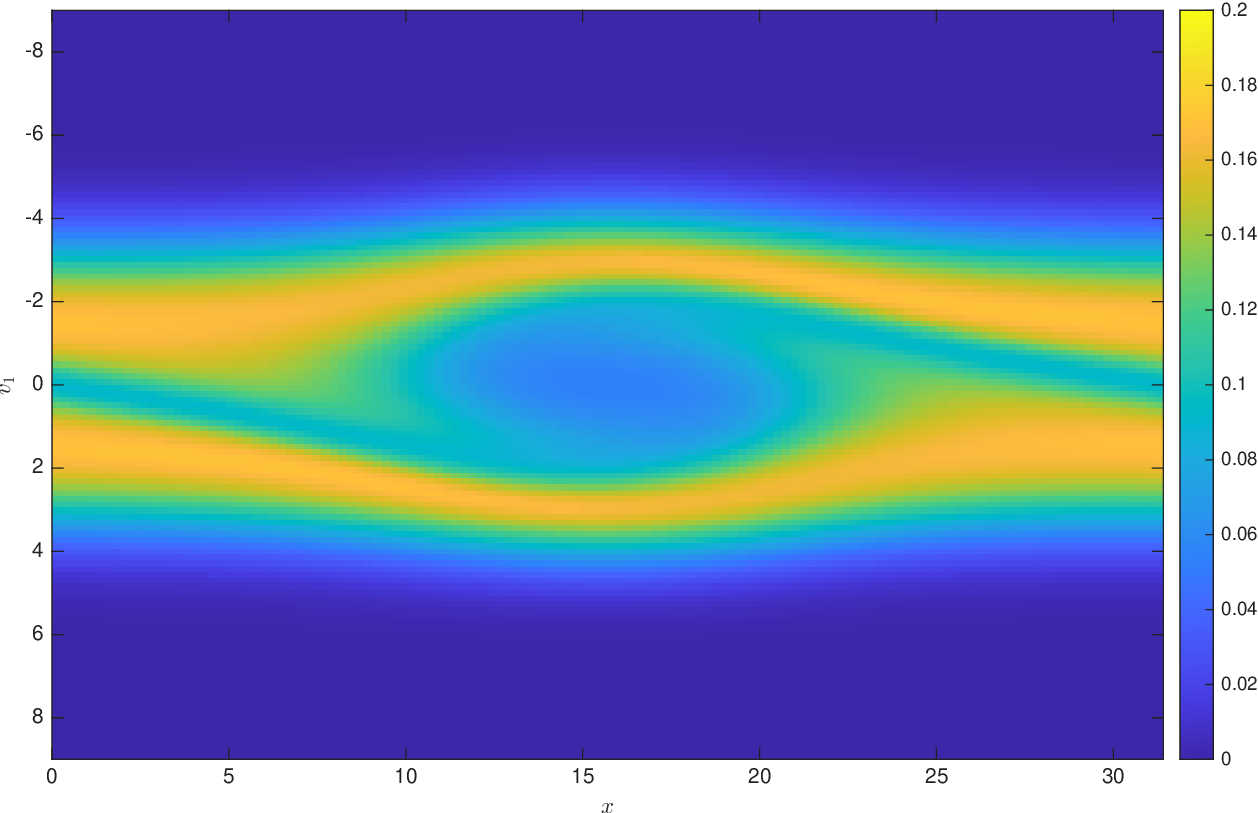}
         \caption{$\eta = 0.004,t=30$}
     \end{subfigure}
     \begin{subfigure}[b]{0.22\textwidth}
         \centering
         \includegraphics[width=\textwidth]{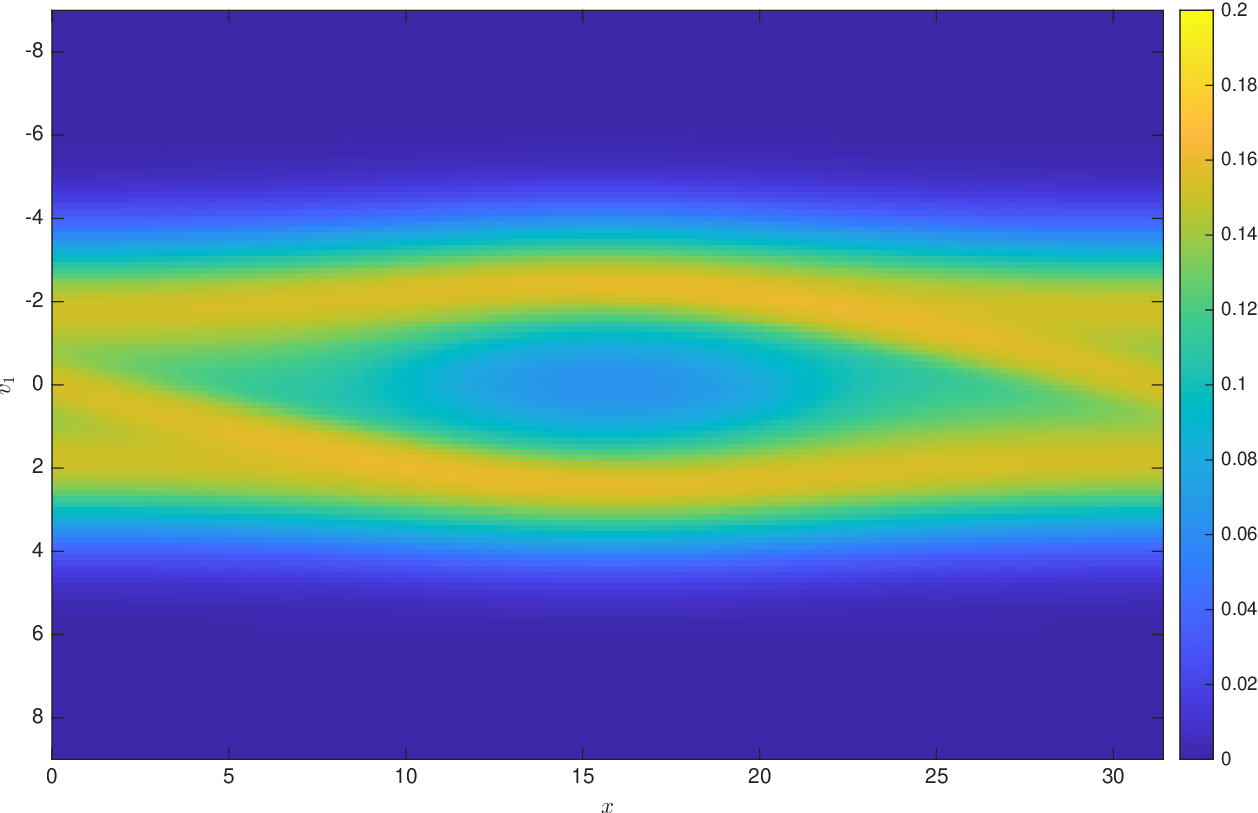}
         \caption{$\eta =0.004,t=45$}
     \end{subfigure}\\
     \begin{subfigure}[b]{0.22\textwidth}
         \centering
         \includegraphics[width=\textwidth]{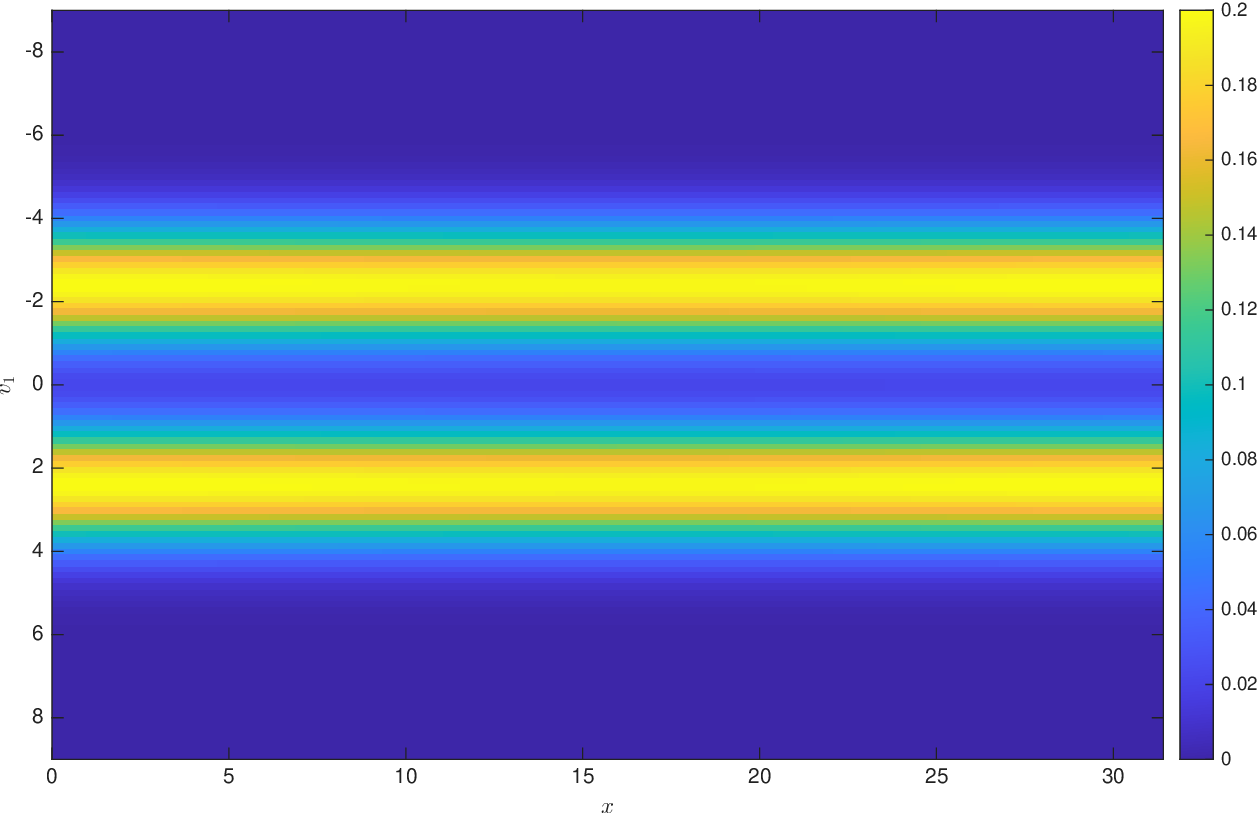}
         \caption{$\eta= 0.006,t=0$}
     \end{subfigure}
     \begin{subfigure}[b]{0.22\textwidth}
         \centering
         \includegraphics[width=\textwidth]{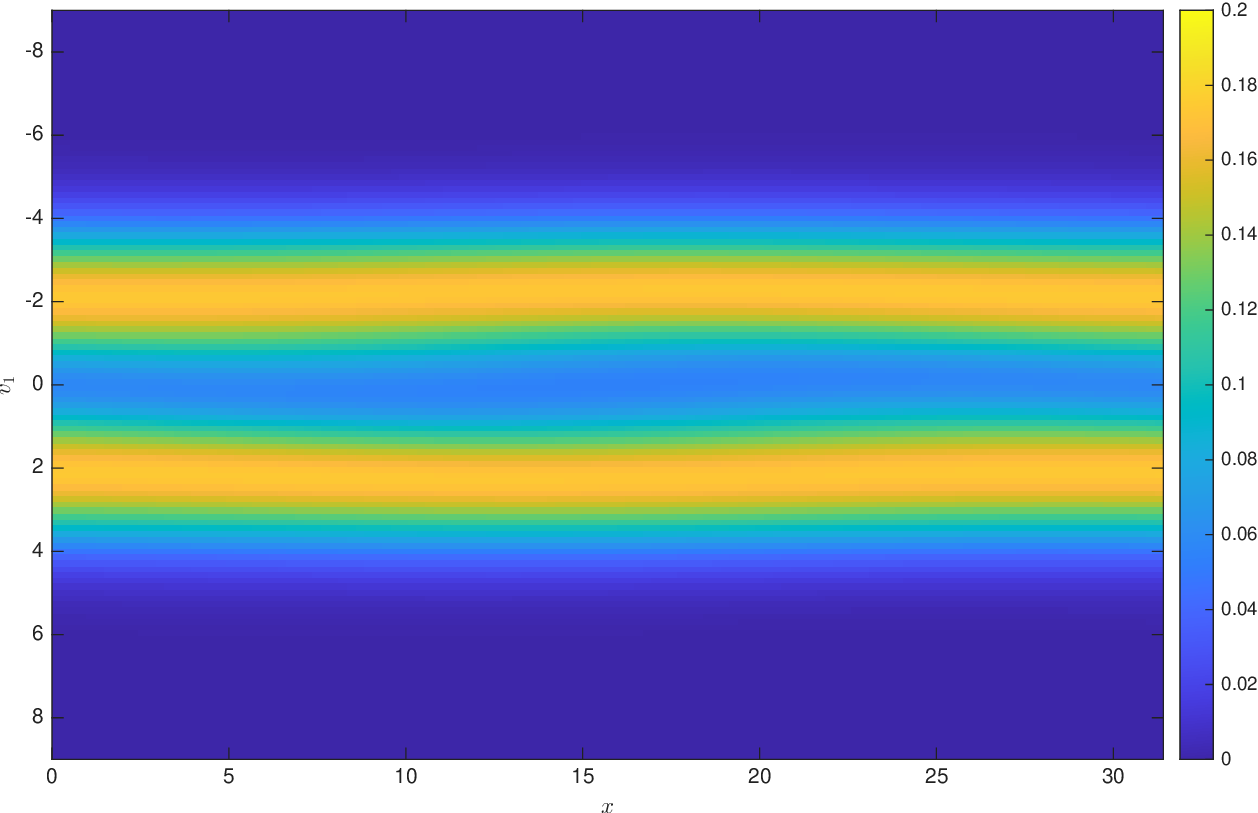}
         \caption{$\eta = 0.006,t=15$}
     \end{subfigure}
     \begin{subfigure}[b]{0.22\textwidth}
         \centering
         \includegraphics[width=\textwidth]{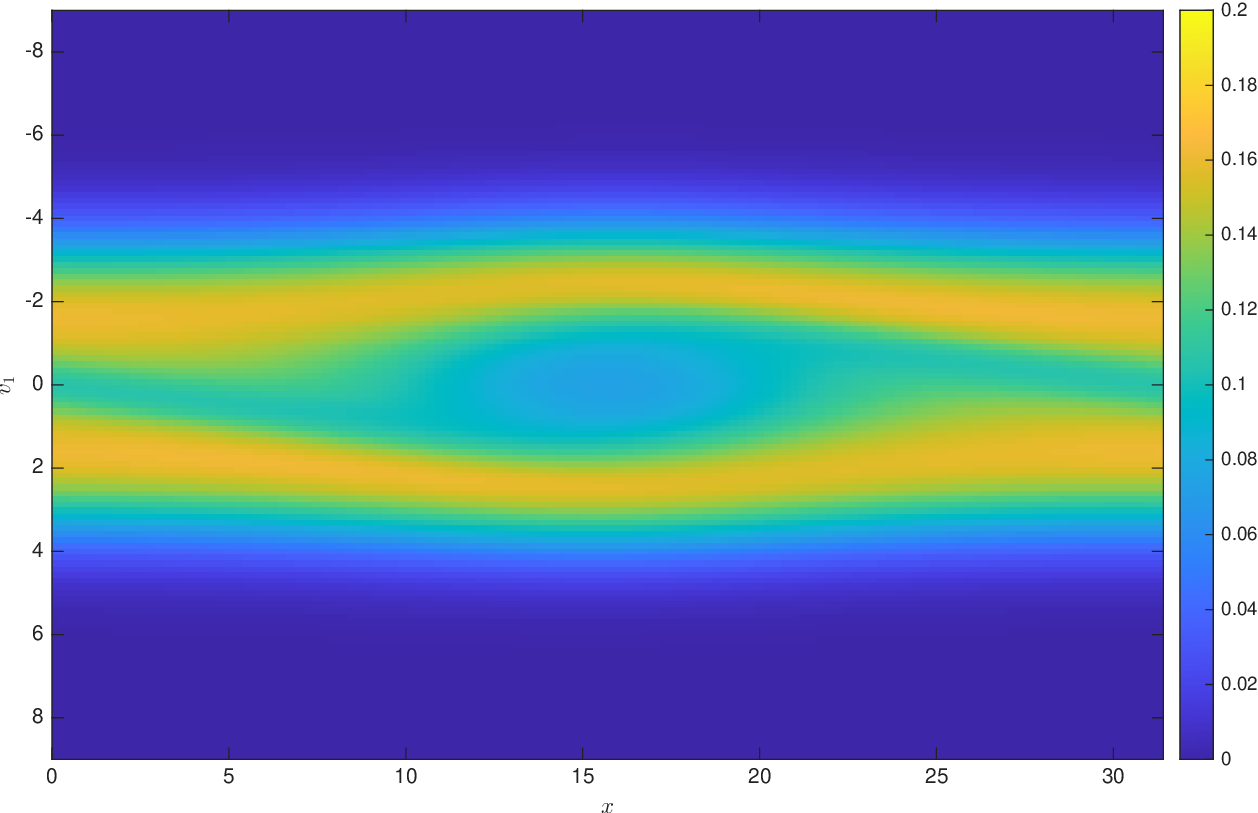}
         \caption{$\eta = 0.006,t=30$}
     \end{subfigure}
     \begin{subfigure}[b]{0.22\textwidth}
         \centering
         \includegraphics[width=\textwidth]{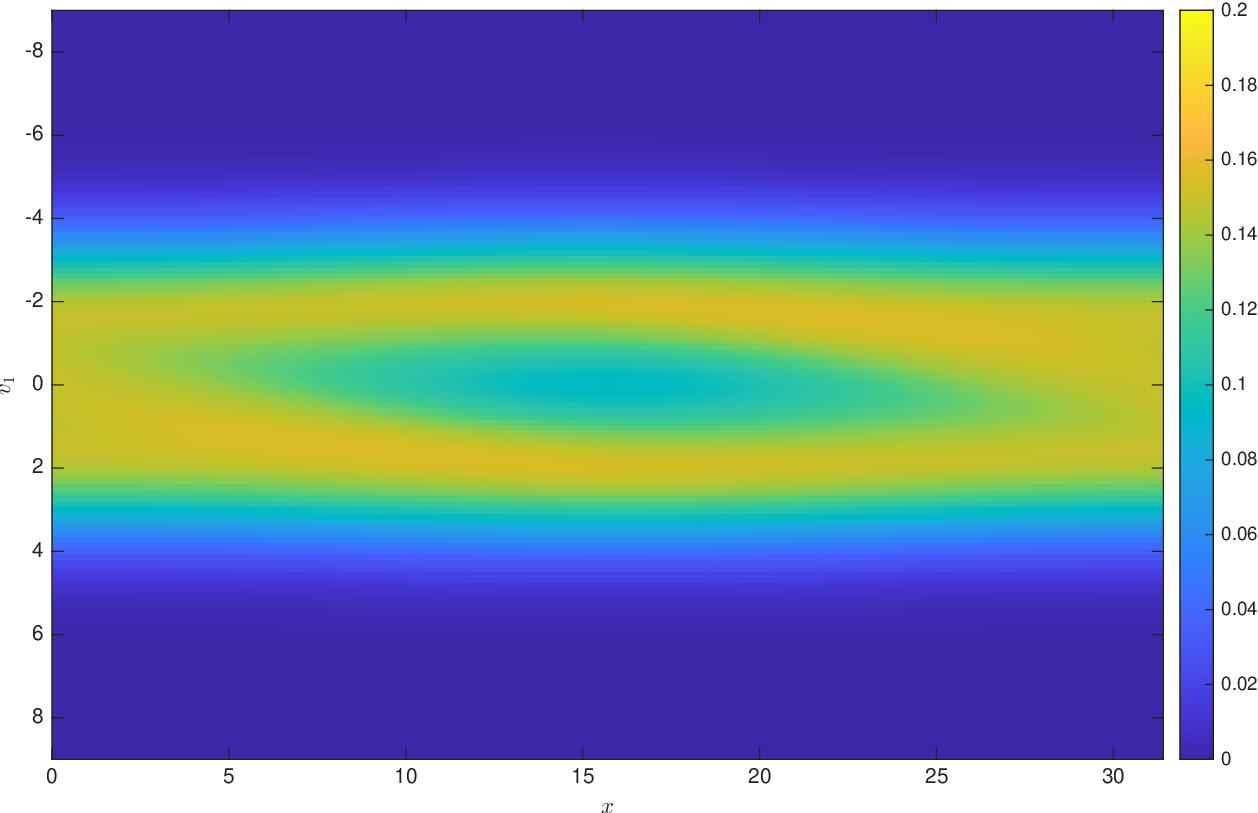}
         \caption{$\eta =0.006,t=45$}
     \end{subfigure}\\
     \begin{subfigure}[b]{0.22\textwidth}
         \centering
         \includegraphics[width=\textwidth]{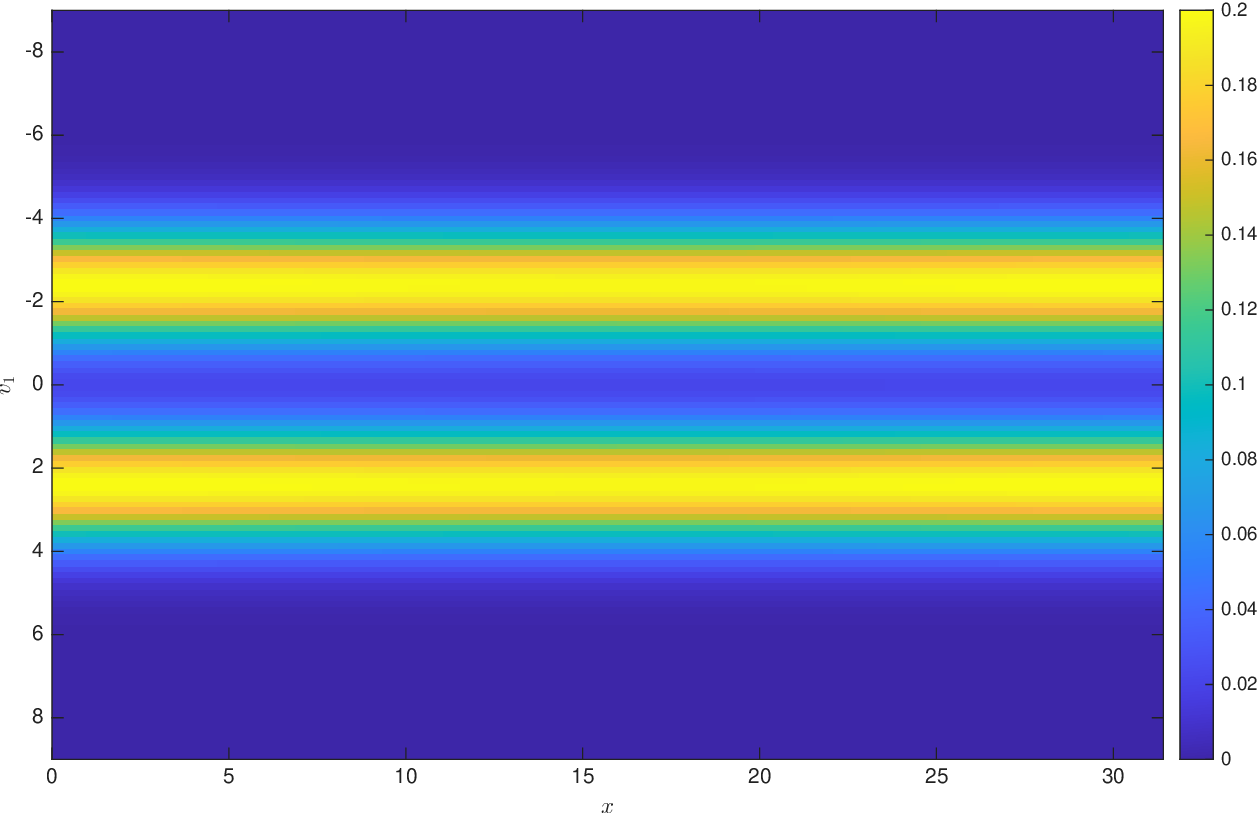}
         \caption{$\eta= 0.008,t=0$}
     \end{subfigure}
     \begin{subfigure}[b]{0.22\textwidth}
         \centering
         \includegraphics[width=\textwidth]{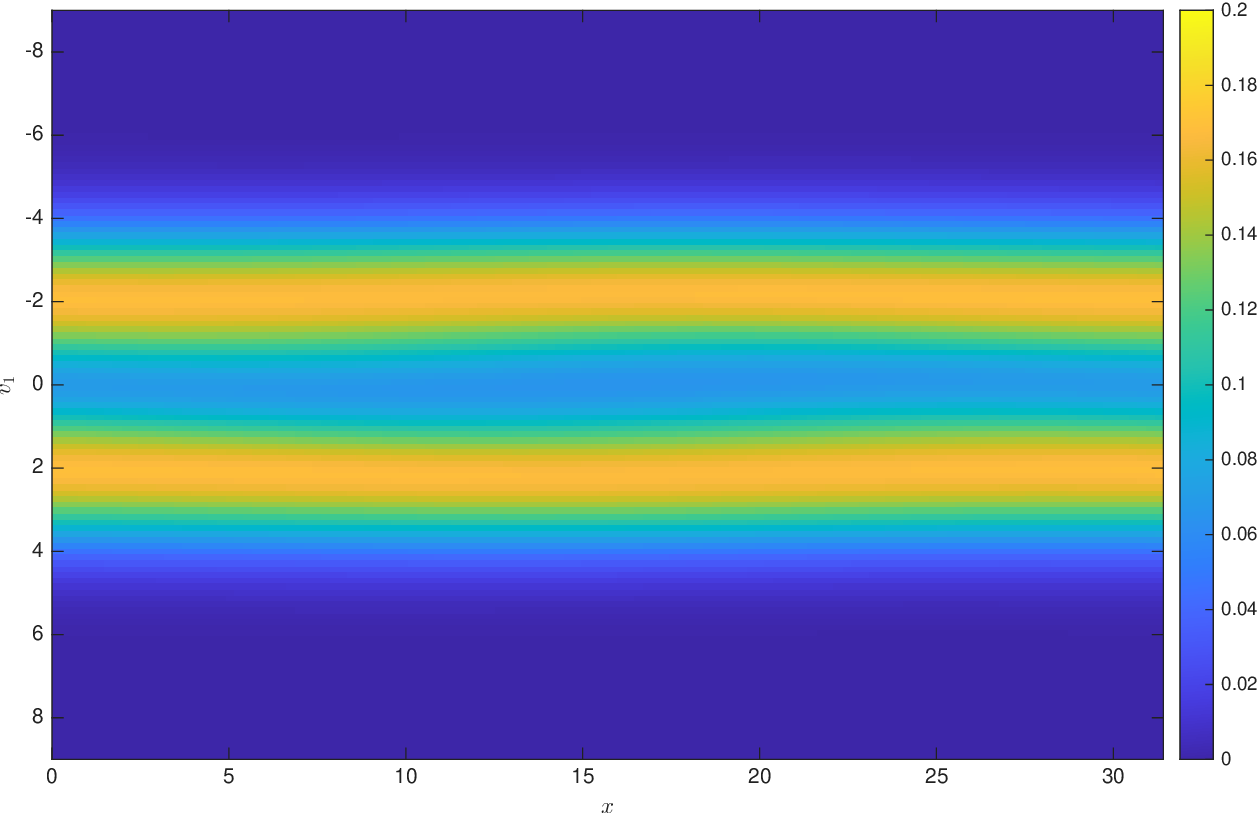}
         \caption{$\eta = 0.008,t=15$}
     \end{subfigure}
     \begin{subfigure}[b]{0.22\textwidth}
         \centering
         \includegraphics[width=\textwidth]{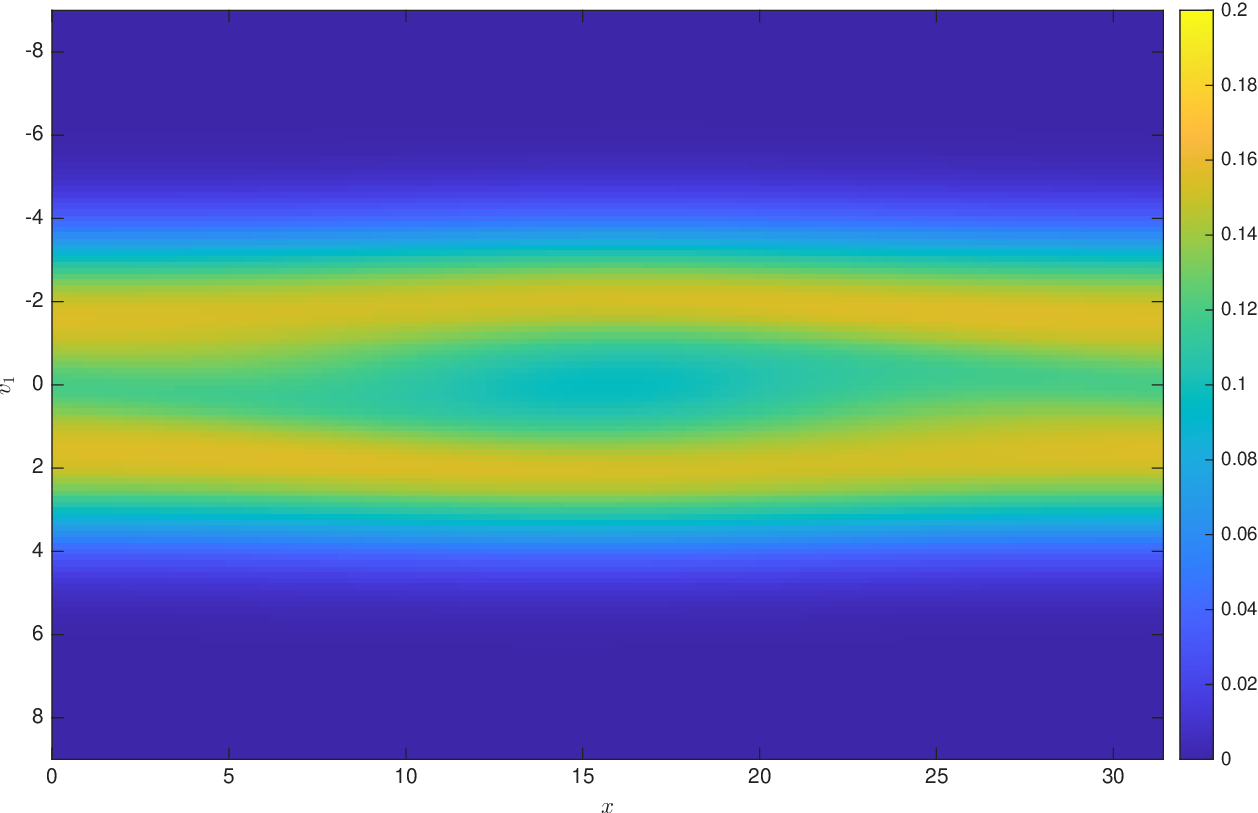}
         \caption{$\eta = 0.008,t=30$}
     \end{subfigure}
     \begin{subfigure}[b]{0.22\textwidth}
         \centering
         \includegraphics[width=\textwidth]{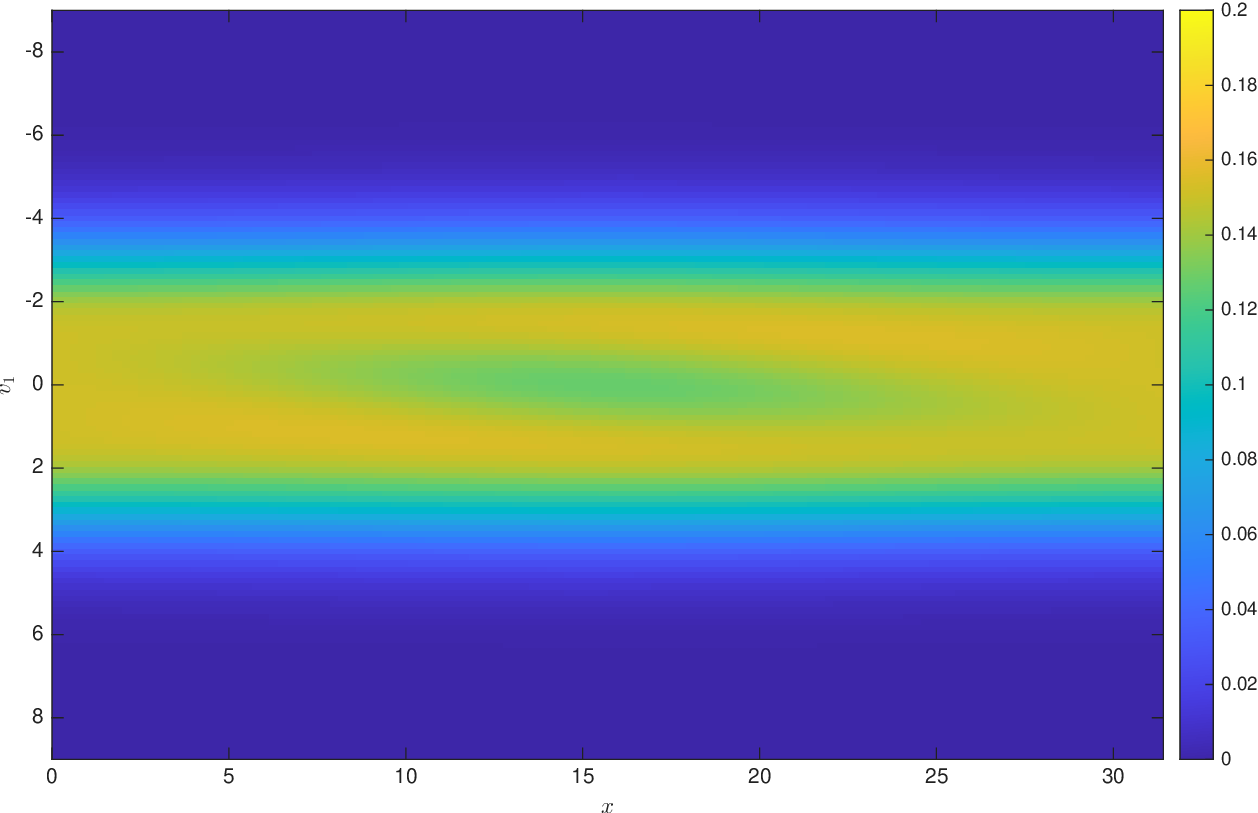}
         \caption{$\eta =0.008,t=45$}
     \end{subfigure}
    \caption{Two-stream instability. Phase plots for different collision strengths.}
    \label{fig_TSI_phase_plot}
\end{figure}
\begin{figure}
     \centering
     \begin{subfigure}[b]{0.24\textwidth}
         \centering
         \includegraphics[width=\textwidth]{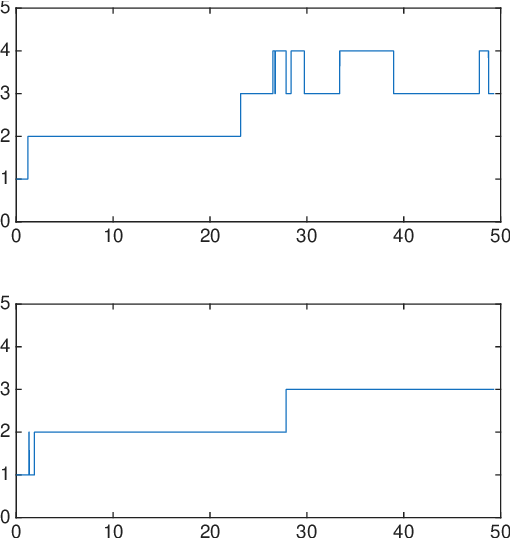}
         \caption{$\eta =0.002$}
     \end{subfigure}
     \begin{subfigure}[b]{0.24\textwidth}
         \centering
         \includegraphics[width=\textwidth]{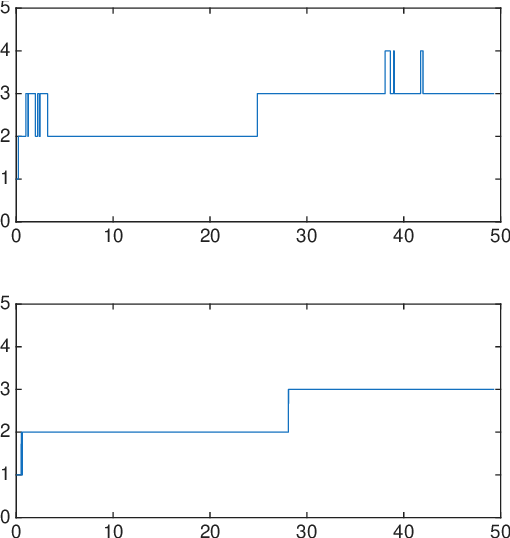}
         \caption{$\eta = 0.004$}
     \end{subfigure}
     \begin{subfigure}[b]{0.24\textwidth}
         \centering
         \includegraphics[width=\textwidth]{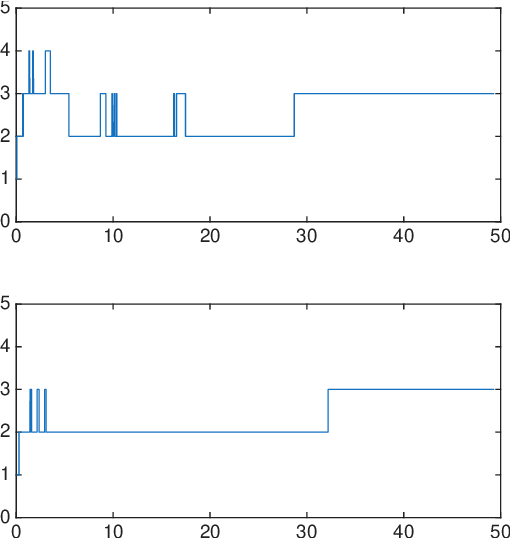}
         \caption{$\eta =0.006$}
     \end{subfigure}
     \begin{subfigure}[b]{0.24\textwidth}
         \centering
         \includegraphics[width=\textwidth]{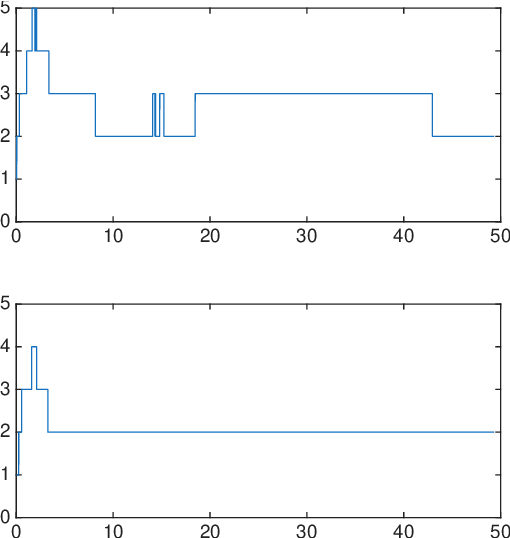}
         \caption{$\eta = 0.008$}
     \end{subfigure}
        \caption{Two-stream instability. Evolution of effective ranks for different collision strengths.}
        \label{fig_TSI_rank}
\end{figure}
We also carry out numerical experiments for the two-stream instability.
The initial condition is given by
\begin{equation}
    f_0(x,\bv) = \frac{1}{2(2\pi)^{3/2}}
    \left(
        1+A\cos(\kappa x)
    \right) \left(
        \e^{-(v_1-v^*)^2/2}
        + \e^{-(v_1+v^*)^2/2}
    \right)
    \e^{-v_2^2/2}
    \e^{-v_3^2/2}.
\end{equation}
The initial electric field $E^{(1)}$ can be computed explicitly:
\begin{equation}
    E^{(1)}(0,x) = -\frac{A}{\kappa} \sin(\kappa x).
\end{equation}
We choose the following parameters $A=0.005$, $\kappa = 0.2$ and $v^*=2.4$.
The spatial domain is $x\in[0,2\pi/\kappa]=[0,10\pi]$ with $N_x=128$ and periodic boundary condition.
The velocity domain is truncated to $[v_{\min},v_{\max}]^3=[-9,9]^3$ with $N_v = 128$. The time step is chosen as in \eqref{eq_time_step_CFL}.
The TT-rank is fixed as $(5,5)$ during the simulation.
We choose different collision strengths as $\eta = 0, 0.002,0.004,0.006,0.008$.
The evolution of the electric energy is shown in \Cref{fig_TSI_energy}, and 
the evolution of the phase space is plotted in \Cref{fig_TSI_phase_plot}. The effective TT-ranks are shown in  \Cref{fig_TSI_rank}, where we omit the case $\eta = 0$ since the effective rank remains 1 throughout the entire simulation.
The results show a clear vortex structure when collisions are absent. 
When collisions are included, the vortex structure gradually smears out as time evolves.
Stronger collisions tend to drive the solution closer to the Maxwellian.
Throughout the simulation, the effectively ranks again remain low, indicating the efficiency of the TT representation.

\section{Conclusion}
\label{sec_future_work}
In this paper, we presented a dynamical tensor-train method applied to a large class of kinetic equations, in which the velocity space is discretized using tensor trains, while the spatial variable is treated as a parameter.
Since the local equilibrium of kinetic equations admit a TT-rank of $(1,1)$, we expect that this discretization enables the use of relatively small TT-ranks when the system is close to equilibrium. A series of numerical examples including the spatially homogeneous and inhomogeneous cases confirmed the efficiency and accuracy of the method.

In addition to the numerical examples presented in this paper, we outline several flexible aspects of the proposed method that are not implemented here but will be considered in future work.
\begin{itemize}
    \item \emph{Different bond dimensions within a tensor train.} In all experiments, the TT-ranks of the tensor trains are chosen as $(r,r)$.
    However, it is not necessary to keep the bond dimensions between the first two and the last two modes identical. 
    In general, one may select TT-ranks of the form $(r_1,r_2)$ with $r_1\not=r_2$. 
    Allowing for different bond dimensions can be advantageous, as it provides greater flexibility to adapt to possible anisotropies in the solution and may further reduce computational cost without sacrificing accuracy.
    \item \emph{Domain decomposition.} In our method, tensor trains at different spatial locations only interact through the evaluation of $\mathcal{G}^j_{\Delta t}$ and $\mathcal{H}^j_{\Delta t}$.
    In most cases, there is no restriction on the TT-ranks of these tensor trains. 
    In other words, different TT-ranks can be chosen independently at different spatial grid points.
    This flexibility enables adaptive rank selection across the spatial domain, allowing the method to allocate higher ranks where the solution exhibits more complexity while keeping ranks small in regions close to equilibrium, thereby improving efficiency.
    \item \emph{Rank adaptivity.} Like most dynamical low-rank approaches, our method allows for rank-adaptivity during the simulation.
    In practice, the ranks can be dynamically increased or decreased through truncation strategies or error-based criteria, ensuring that the representation remains both accurate and efficient.
\end{itemize}

In most numerical examples of this paper (except for the stiff spatially inhomogeneous BGK equation), we employ the explicit time-stepping method in the low-rank algorithm, chosen for its simplicity and ease of implementation.
Since implicit and IMEX methods are generally more stable, another interesting direction is to explore efficient implicit implementation within the framework of dynamical tensor trains.


\section*{Acknowledgement}
This work was partially supported by DOE grant DE-SC0023164, NSF grants DMS-2409858 and IIS-2433957, and DoD MURI grant FA9550-24-1-0254.

\bibliographystyle{abbrv}
\bibliography{myBib_abbr}

\end{document}